\documentclass[poms,final,nonblindrev]{poms1V1}

\OneAndAHalfSpacedXI

\usepackage{graphicx}
\usepackage{subfigure, epsfig}
\usepackage{natbib}

 \bibpunct[, ]{(}{)}{,}{a}{}{,}
 \def\bibfont{\small}

\TheoremsNumberedThrough
\ECRepeatTheorems

\EquationsNumberedThrough

\usepackage{hyperref}
    \hypersetup{
        colorlinks   = true,
        citecolor    = blue,
        linkcolor=blue
    }
\usepackage{adjustbox}
\usepackage{placeins}
\usepackage{array, makecell}
\usepackage{color, colortbl}
\definecolor{Gray}{gray}{0.9}
\usepackage{setspace,lipsum}
\usepackage{style}
\usepackage{xcolor}
\usepackage{arydshln}
\usepackage{pdflscape}
\usepackage{rotating}

\newcommand{\srevision}[1]{{\color{black}{#1}}}

\begin{document}

\RUNAUTHOR{S. Dey, A.K. Kurbanzade, E. S. Gel, J. Mihaljevic, S. Mehrotra}
\RUNTITLE{Review of Pandemic Vaccine Supply Chain Modeling}
\TITLE{Optimization Modeling for Pandemic Vaccine Supply Chain Management: A Review \textcolor{black}{and Future Research Opportunities}}

\ARTICLEAUTHORS{
\AUTHOR{Shibshankar Dey {$^{a,b}$} $\bullet$ Ali Kaan Kurbanzade {$^{a,b}$} $\bullet$ Esma S. Gel {$^{c}$} \\ Joseph Mihaljevic {$^{d}$} $\bullet$ Sanjay Mehrotra {$^{a,b,\ast}$}}
\AFF{$^a$ Department of Industrial Engineering and Management Sciences, Northwestern University, Evanston, IL, USA \\ 
$^b$ Center for Engineering and Health, Northwestern University, Feinberg School of Medicine, Chicago, IL, USA \\
$^c$ Department of Supply Chain Management and Analytics, University of Nebraska-Lincoln, Lincoln, NB, USA \\
$^d$ School of Informatics, Computing, and Cyber Systems, Northern Arizona University, Flagstaff, AZ, USA \\
$^\ast$Corresponsing author: \EMAIL{mehrotra@northwestern.edu}}
}

\ABSTRACT{
During various stages of the COVID-19 pandemic, countries implemented diverse vaccine management approaches, influenced by variations in infrastructure and socio-economic conditions. This article provides a comprehensive overview of optimization models developed by the research community throughout the COVID-19 era, aimed at enhancing vaccine distribution and establishing a standardized framework for future pandemic preparedness. These models address critical issues such as site selection, inventory management, allocation strategies, distribution logistics, and route optimization encountered during the COVID-19 crisis. A unified framework is employed to describe the models, emphasizing their integration with epidemiological models to facilitate a holistic understanding. \srevision{This article also summarizes evolving nature of literature, relevant research gaps, and authors' perspectives for model selection. Finally, future research scopes are detailed both in the context of modeling and solutions approaches.}
}

\HISTORY{The initial draft was first submitted in May, 2023. \textcolor{black}{The revised version is submitted in November, 2023.}}

\maketitle

\vspace{-1cm}

\section{Introduction}

The COVID-19 pandemic emerged as the largest global health crisis in recent history, with an staggering incidence of over 750 million reported cases and death toll approaching seven million individuals worldwide~\citep{WHO-COVID-19}. Beyond the profound health implications, the pandemic has presented an array of multifaceted challenges encompassing social, technological, economic, environmental, political, legal, and ethical domains \citep{sampath2021pandemics,taskinsoy2020covid,klemevs2020minimising,shah2021prevalence,saberi2020ethical}. In response, numerous researchers have delved into crucial issues such as resource allocation for essential care resources like ventilators and hospital beds/capacity \citep{mehrotra2020model,bozkir2023capacity,yang2021design}, as well as investigating the impact of community mitigation strategies like social distancing and mask-wearing \citep{keskinocak2020impact,gel2020covid} during the initial stages of the pandemic. However, with the rapid development of highly effective vaccines, the research efforts have progressively shifted towards studying operational dilemmas related to vaccine allocation and distribution.

The COVID-19 pandemic underscored the significance of timely and widespread vaccine deployment, particularly when highly effective vaccines are shown to prevent infections and hospitalizations to achieve herd immunity against novel pathogens~\citep{macintyre2022modelling}. However, the operational challenges associated with rapidly developed vaccines, such as those created and approved for COVID-19, necessitate careful consideration to ensure efficient, effective, and equitable deployment~\citep{huang2012models}. Primarily, vaccine availability is typically limited, especially in the early stages of production, requiring prioritized allocation of scarce doses among various subgroups and geographical regions. Additionally, the diverse range of vaccines developed by different manufacturers offer varying levels and types of immunity against different pathogen strains. These vaccines also present different storage requirements (shelf life, doses per vial, cold storage, etc.) and dosing protocols (single or multiple shots), thereby complicating transportation, allocation, and distribution logistics. From a global perspective, addressing these challenges requires tailored solutions in different regions due to factors such as varying vaccine availability, healthcare infrastructure, socioeconomic pressures, education levels, and risk perceptions among citizens. Moreover, different countries may prioritize subpopulations and employ conflicting economic measures, leading to the adoption of diverse vaccine deployment and delivery strategies~\citep{saadi2021models}.

In response to the various operational challenges encountered during the COVID-19 pandemic, the Operations Research/Management Science community has generated a significant body of research comprising modeling and optimization studies, introducing novel dimensions and capabilities to address operational problems in this domain. Some of these optimization models incorporate dynamic considerations of disease transmission and infection stages, enabling accurate estimation of population-level vaccination requirements. This literature review focuses on optimization modeling for pandemic vaccine deployment, offering a comprehensive overview of these papers and the current state-of-the-art in this field. The primary objectives of this review are to (i) advocate for the adoption of such models in decision-making processes by public health officials in future pandemics and (ii) identify research gaps that hinder ongoing pandemic preparedness efforts undertaken by our nation and various global agencies.

Our search methodology involved identifying relevant articles published or made available online between the years 2020 and 2023, employing various combinations of keywords listed in the first, second, and third columns of Table~\ref{table:keywords}. Through this approach, we identified a \textcolor{black}{total of 109 articles} that met the inclusion criteria. Additionally, to ensure comprehensive coverage of pertinent topics, we incorporated a selection of highly cited articles published between 2011 and 2019, as they offered valuable insights into issues of interest such as equitable vaccine allocation, vial size, and the vaccine supply chain. \textcolor{black}{Figure \ref{fig:PRISMA} in the Appendix outlines the paper selection and classification process as a flow chart.}

\begin{table}[htp] \setlength{\tabcolsep}{8pt}
\begin{center}
\caption{List of keywords used in the initial search.}
\begin{normalsize}
\renewcommand{\arraystretch}{1}
\begin{tabular}{lll}
\hline
\rowcolor{Gray}
\textbf{Keyword Set 1} & \textbf{Keyword Set 2} & \textbf{Keyword Set 3} \\
\hline
coronavirus vaccine & allocation & mathematical modelling \\
COVID-19 vaccine & distribution & mathematical programming \\
influenza vaccine & equity & network \\
pandemic vaccine & facility location problem & optimization \\
vaccination & inventory & optimization model \\
& supply chain management & simulation \\
\hline
\end{tabular}
\end{normalsize}
\label{table:keywords}
\end{center}
\end{table}

The articles covered in this review can be categorized into three main groups (i) Optimization articles, which encompass resource allocation, location-allocation, inventory management, supply chain management, and routing problems, employing at least one mathematical optimization model. Some of these articles incorporate simulation and compartmental epidemiological dynamics models, such as Susceptible-Infected-Recovered (SIR) models, integrated with the optimization models. (ii) Simulation-based articles, which focus on prioritization, management of multiple vaccine, and inventory control, utilizing simulation modeling to compare various scenarios or interventions of interest. These articles do not include an optimization model but rely on simulation-based approaches. (iii) Other articles, involving descriptive analysis, conceptualization, vaccine hesitancy, and waste management. These articles do not involve optimization or simulation models but offer valuable insights to guide future research in this field. \textcolor{black}{In total, our review covers 77 optimization model-based articles, 17 simulation-based articles, and 15 articles employing alternative approaches. Notably, 76 out of the 94 articles involving optimization and simulation modeling specifically pertain to COVID-19 vaccine studies,} highlighting the significant recent interest in pandemic vaccine management prompted by the COVID-19 pandemic.

The subsequent sections of this paper are organized as follows: Section \S \ref{sec:framework} presents the vaccine supply chain framework along with the nomenclature adopted throughout the paper. In Sections \S\ref{sec:model}, and \S\ref{sec: pure facility location model}, we propose a comprehensive optimization modeling framework for COVID-19 vaccine allocation and distribution management. These sections respectively address supply chain and location problems using mixed-integer linear programming techniques. Section \S \ref{sec: uncertainty-scm} reviews the optimization technique found in vaccine logistics literature in presence of uncertainties. Additionally, Section \S \ref{sec:eq-and-sus} provides an overview of how equity and sustainability are incorporated into the models. In Section \S \ref{sec:SEIRD models}, we outline a generalized epidemiological modeling framework utilizing an extended version of the compartmental SIR model called DELPHI-V and present a framework for integrating optimization and epidemiological models. \textcolor{black}{In Section \S \ref{sec:Methods}, we provide a classification of methodological considerations and highlight the important insights from case studies. In Section \S \ref{sec:conclusion}, we offer an overall discussion of the literature as well as our own perspectives for model selection. In Section \S \ref{sec:Critic}, we provide a critique and suggestions for future research directions.} {Appendix~\S\ref{sec: vehicle routing}} covers routing problems. Finally, {Appendix \S \ref{sec:sim-and-others}} covers simulation-based and other relevant articles in alignment with the optimization-based articles.

\newpage

\section{Vaccine Supply Chain Framework and Model Notations} \label{sec:framework}

\begin{figure}
\begin{center}
\includegraphics[scale=0.03]{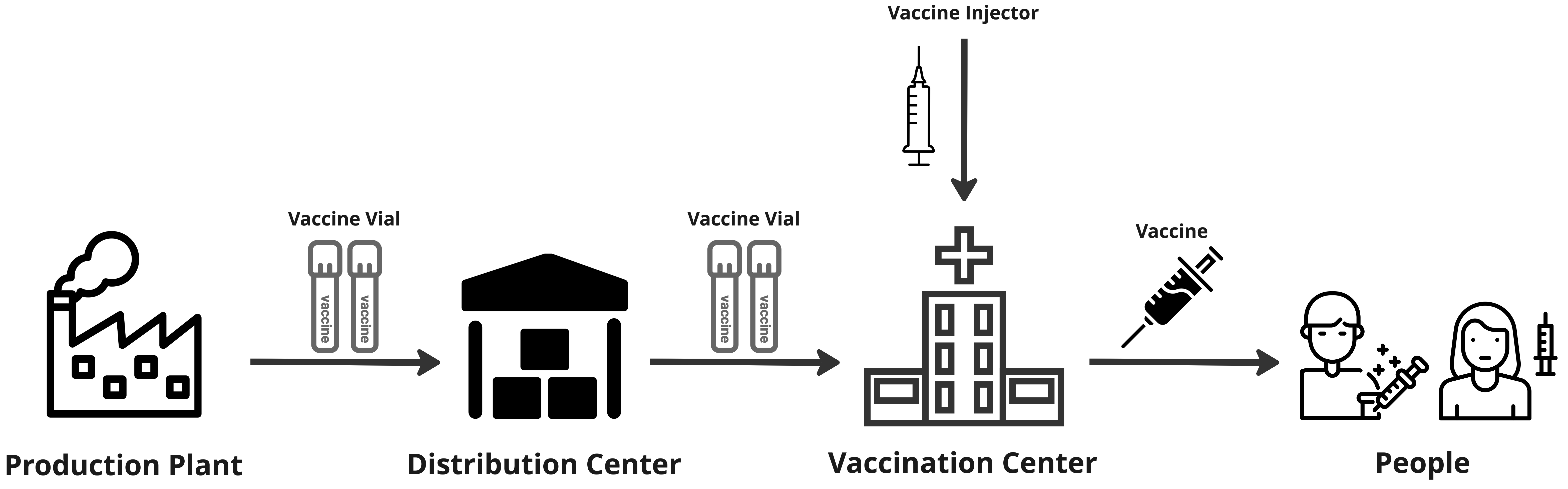}
\caption{Vaccine Supply Chain} \medskip
\label{Img:Vaccine-SC}
\end{center}
\end{figure}

The vaccine supply chain comprises of four primary echelons: manufacturers (M), distribution centers (DC), vaccination centers (VC), and population demand (PS), as depicted in Figure~\ref{Img:Vaccine-SC}. Production facilities can be located domestically or internationally, catering to multiple distribution centers (DCs) or warehouses at various levels, such as national, state, city, regional, or district levels. The fundamental unit of inventory is the vaccine vial, which is typically transported in refrigerated trucks to numerous distribution centers (DCs) that serve geographically dispersed vaccination centers (VCs). These vaccination centers (VCs) encompass hospitals, community health centers, pharmacies, and mobile vaccination clinics. Depending on vaccine characteristics and availability, certain subpopulations may be prioritized for vaccination based on factors like age, occupation, or risk level. It is important to note that vaccine vials can contain multiple doses, necessitating the consideration of shelf-life and the duration for which an open vial remains viable. In this paper, we assume that each vial contains a single vaccine dose. Furthermore, the paper employs four distinct time frames: ordering time, receiving time, opening time, and vaccination time. The primary types of decisions addressed in the surveyed papers are as follows:
\begin{itemize}
\item \textit{Location:} Decisions regarding locating facilities, such as mass vaccination centers. 
\item \textit{Inventory:} Decisions on managing various aspects of vaccine inventory. 
\item \textit{Allocation:} Joint decisions on vaccine allocation and assignment of demand nodes to locations. 
\item \textit{Distribution:} Decisions regarding distribution of vaccines to demand nodes.
\item \textit{Routing:} Decisions regarding the management of vehicles/mobile vaccination clinics visiting depots, vaccination centers (VCs), etc. 
\end{itemize} 

Tables~\ref{Tab:Index Set description}, \ref{Tab:parameter description}, and \ref{Tab:Decision Variable Description} provide the index sets, parameters, and decision variables employed in our overarching optimization model. We strongly suggest the reader to refer to these tables when reading the technical expressions appearing in this review, since they are not necessarily described in the main text. Note that, in the above notation convention, the subscripts refer to the indexes. Letters $M$, $D$, $m$, $d$, $v$, etc. are used as superscripts in many variables or parameters to highlight their association with the manufacturer (M), distribution center (DC), and vaccination center (VC), respectively. The manufacturer set $\mathcal{M}$ can also be defined as $\{1, 2, \cdots, |\mathcal{M}|\}$ where $|\cdot|$ is its cardinality. All sets are defined using calligraphic letters. As necessary, we also interchangeably use $[\cdot]$ to represent an index set with the entry `$\cdot$' representing the last entry of that set starting from~$1$. Decision variable types are indicated with letters B (binary), C (continuous) and I (integer). I/C refers to decision variables that should have been integer but were treated as continuous to reduce complexity.  

Table~\ref{table:optimization-based-classification} classifies the articles based on the specific type of decision being modeled, such as location, inventory, allocation, distribution, or routing. It also indicates whether a compartmental Susceptible-Infected-Recovered (SIR) type epidemiological model was utilized to capture the dynamic nature of the infectious disease being studied. It is worth noting that papers incorporating a SIR-type model often consider demand uncertainty within the modeling framework. 

Optimization articles are further grouped into five categories based on the decisions studied: (i) resource allocation, (ii) location-allocation, (iii) inventory management, (iv) supply chain management, and (v) routing. Further details on each paper in each category  (model type, uncertainty type, commodity, period, echelons, objective, case study, and solution methodology) are compiled in Tables~\ref{table:resource-allocation-articles} through~\ref{table:routing-articles} in Appendix. \textcolor{black}{Table \ref{table:classification} provides an overall summary for all of the optimization articles. 
}

\begin{table}[htp] \setlength{\tabcolsep}{8pt}
\begin{center}
\caption{Summary of articles based on decision levels, model type, commodity, period and case study for each category.}
\begin{adjustbox}{width=\textwidth}
\renewcommand{\arraystretch}{1}
\begin{tabular}{llllllllllllllllll}
\hline
\rowcolor{Gray}
 & & \multicolumn{5}{c}{\textbf{Decision Levels}} & & \multicolumn{2}{c}{\textbf{Model Type}} & \multicolumn{2}{c}{\textbf{Commodity Type}} & & \multicolumn{2}{c}{\textbf{Period}} & \multicolumn{2}{c}{\textbf{Case Study}} & \textbf{Solution Method} \\
\hline
\rowcolor{Gray}
\textbf{Category} & \textbf{Total} & Location & Inventory & Allocation & Distribution & Routing & SIR & Deterministic & Uncertain & Single & Multiple & Multiple Dose & Single & Multiple & Real Life & Synthetic & Commercial Solver \\ 
\hline
\rowcolor{Gray}
\textbf{Supply Chain Models} \\
\hline
Resource Allocation & 23 & 0 (0\%) & 0 (0\%) & 23 (100\%) & 10 (43\%) & 0 (0\%) & 13 (57\%) & 9 (39\%) & 14 (61\%) & 17 (74\%) & 6 (26\%) & 2 (9\%) & 10 (43\%) & 13 (57\%) & 18 (78\%) & 5 (22\%) & 4 (17\%) \\ 
Inventory Management & 13 & 0 (0\%) & 13 (100\%) & 9 (69\%) & 7 (54\%) & 0 (0\%) & 1 (8\%) & 4 (31\%) & 9 (69\%) & 8 (62\%) & 5 (38\%) & 4 (31\%) & 2 (15\%) & 11 (85\%) & 10 (77\%) & 3 (23\%) & 1 (8\%) \\ 
Supply Chain Management & 20 & 20 (100\%) & 18 (90\%) & 20 (100\%) & 16 (80\%) & 2 (10\%) & 0 (0\%) & 9 (45\%) & 11 (55\%) & 13 (65\%) & 8 (40\%) & 0 (0\%) & 7 (35\%) & 14 (70\%) & 19 (95\%) & 1 (5\%) & 4 (20\%) \\ 
\hline
\rowcolor{Gray}
\textbf{Location and Routing Models} \\
\hline
Location-Allocation & 18 & 18 (100\%) & 0 (0\%) & 16 (89\%) & 5 (28\%) & 0 (0\%) & 1 (6\%) & 12 (67\%) & 6 (33\%) & 15 (83\%) & 2 (11\%) & 0 (0\%) & 14 (78\%) & 3 (17\%) & 16 (89\%) & 2 (11\%) & 7 (39\%) \\ 
Routing & 3 & 2 (67\%) & 0 (0\%) & 0 (0\%) & 3 (100\%) & 3 (100\%) & 1 (33\%) & 1 (33\%) & 2 (67\%) & 3 (100\%) & 0 (0\%) & 1 (33\%) & 2 (67\%) & 1 (33\%) & 3 (100\%) & 0 (0\%) & 2 (67\%) \\
\hline
\textbf{Total} & 77 & 40 (52\%) & 31 (40\%) & 68 (88\%) & 41 (53\%) & 5 (6\%) & 16 (21\%) & 35 (45\%) & 42 (55\%) & 56 (73\%) & 21 (27\%) & 7 (9\%) & 35 (45\%) & 42 (55\%) & 66 (86\%) & 11 (14\%) & 18 (23\%) \\ 
\hline
\end{tabular}
\end{adjustbox}
\label{table:classification}
\end{center}
\end{table}

\begin{table} \setlength{\tabcolsep}{4pt}
\caption{Sets}
\begin{center}
\renewcommand{\arraystretch}{1}
\begin{normalsize}
\begin{tabular}{cl}
\hline
\rowcolor{Gray}
\textbf{Notation}&\textbf{Description}\\
\hline
$\mathcal{M}$ & Set of manufacturers (M)/suppliers/vendors, indexed by $i$. \\
$\mathcal{N}$ & Set of distribution centers (DC), indexed by $j$. \\
$\mathcal{V}$ & Set of vaccination centers (VC), indexed by $k$. \\
$\mathcal{V}_j$ & Set of vaccination centers (VC) that can be served by distribution center (DC) $j$, indexed by $k$. \\
$\mathcal{G}$ & Set of risk-based population groups, indexed by $g$. \\
$\mathcal{P}$ & Set of different types of vaccines, indexed by $p$. \\
$\mathcal{P}_i$ & Set of vaccine types that can only come from manufacturer (M) $i$, indexed by $p$. \\
$\mathcal{T}$ & Set of discrete time periods, indexed by $t$. \\
$\mathcal{H}$ & Set of vehicles, indexed by $h$. \\
$\mathcal{R}$ & Set of regions, indexed by $r$ or $r'$. \\
$\mathcal{S}$ & Set of population sites (PS), indexed by $s$. \\
$\mathcal{Q}_s$ & Set of coverage levels to serve population sites (PS)/villages $s$, indexed by $q$. \\
$\mathcal{O}$ & Set of outreach centers/teams, indexed by $\ell$ . \\
$\Omega$ & Set of scenarios representing uncertainty, indexed by $\omega$. \\
$\mathcal{B}$ & Ambiguity set in distributionally robust optimization (DRO). \\
\hline
\end{tabular}
\end{normalsize}
\end{center}
\label{Tab:Index Set description}
\end{table}

\begin{table} \setlength{\tabcolsep}{4pt}
\vspace{-.3cm}
\caption{Parameters}
\vspace{-0.5cm}
\begin{center}
\begin{adjustbox}{width=\textwidth}
\renewcommand{\arraystretch}{0.9}
\begin{tabular}{cll}
\hline
\rowcolor{Gray}
\textbf{Notation}&\textbf{Description}&\textbf{Equations}\\
\hline
\rowcolor{Gray}
\multicolumn{2}{l}{\textbf{Capacity}} \\
\hdashline
$C^M_{ipt}$ & Production capacity of manufacturer $i$ at time period $t$ for type $p$ of vaccine. & \ref{eq:capacity_man0}\\
$C^D_j$ & Capacity of distribution center (DC) $j$. & \ref{eq: inventory-capacity} \\
$C^V_k (C^V) $ & Capacity of vaccination center (VC) $k$ (Overall capacity/availability over all VCs). & \ref{eq: p-median and max coverage mixture},\ref{eq: herd-immunity} \\
$C^F$ & Fleet size availability. & \ref{eq: fleet size capacity} \\
$S^D_j$ & Safety stock level of distribution center (DC) $j$. & \ref{eq: inventory-capacity} \\ 
$C^O_{\ell}$ & Capacity of outreach center (OC) $\ell$. & \ref{eq: step-wise detail coverage1} \\
$C^H_h$ & Capacity of vehicle $h$. & \ref{eq: vehicle-routing} \\
$C^{D(c)}_j$ & Cold refrigeration capacity of distribution center (DC) $j$. & \ref{eq: Capacity with Forcing} \\
$C^{D(vc)}_j$ & Very cold refrigeration capacity of distribution center (DC) $j$. & \ref{eq: Capacity with Forcing2} \\
$C^{D(uc)}_j$ & Ultra cold refrigeration capacity of distribution center (DC) $j$. & \ref{eq: Capacity with Forcing2},\ref{eq: Capacity with Forcing3} \\
\hdashline
\rowcolor{Gray}
\multicolumn{2}{l}{\textbf{Cost}} \\
\hdashline
$c^{SH}_{k p}$ & Shortage cost for vaccine $p$ at vaccination center (VC) $k$. & \textit{In Text} \\
$c^{EI}$ & Environmental impact cost for establishing a facility. & \ref{eq:env obj} \\
$c^{CR}$ & Amount of carbon emmission by transshipping a unit of vaccine. & \ref{eq:env obj} \\
$c_{k}$ & Cost of locating vaccination center (VC) at $k$. & \ref{eq: p-median2},\ref{eq: step-wise coverage} \\
$\mathsf{b}$ & Total budget. & \ref{eq: p-median2},\ref{eq: step-wise coverage} \\
\hdashline
\rowcolor{Gray}
\multicolumn{2}{l}{\textbf{Demand}} \\
\hdashline
$d_{g k p t}$ & Demand of vaccine type $p$ among population group $g$ at location $k$ in time period $t$. & \ref{eq:OPen Vial Wastage in backward time},\ref{eq:OPen Vial Wastage in backward time2},\ref{eq:priority group focus},\ref{eq:priority group focus2},\ref{eq:priority group focus3}, \\
& & \ref{eq:priority group focus4},\ref{eq: equity based objective},\ref{eq: equity based objective 2} \\
$d_{kpt}$ & Demand for vaccine type $p$ at vaccination center (VC) $k$ at time period $t$. & \ref{eq: Inventory_balance_VC},\ref{eq:shortage-demand} \\
$d_{r g}$ & Demand of region $r$ for risk group $g$. & \ref{eq: vehicle-routing},\ref{eq: equity based objective3.1},\ref{eq: equity based objective3.2} \\
$d^D_j$ & Demand (average demand rate) of distribution center (DC) $j$. & \ref{eq:uncertain-DC} \\
$\gamma$ & Minimum demand satisfaction ratio. & \ref{eq: equity based objective},\ref{eq: equity based objective 2} \\
$d^{P}_{s}$ & Demand of population site (PS) $s$. & \ref{eq: p-median and max coverage mixture} \\
$\bar{P}_{s}$ & Coverage coefficient at population site (PS) $s$. & \ref{eq: p-median2},\ref{eq: step-wise coverage} \\
$\theta_q$ & Coverage fraction attainable when the nearest serving center is between $\mathfrak{y}_{q - 1}$ and $\mathfrak{y}_{q}$. & \ref{eq: step-wise coverage} \\
$\hat{P}_r (\hat{P}_s)$ & Population size of region $r$ (of population site (PS) $s$). & \ref{eq: herd-immunity}\\
\hdashline
\rowcolor{Gray}
\multicolumn{2}{l}{\textbf{Inventory}} \\
\hdashline
$l_j$ & Lead time to deliver order to distribution center (DC) $j$. & \ref{eq:uncertain-DC} \\
$c^{OC}_{i j}$ & Ordering cost to order from manufacturer (M) $i$ by distribution center (DC) $j$. & \ref{eq:uncertain-DC} \\
$c^{HC}_{j}$ & Holding cost to hold inventory at facility $j$. & \ref{eq:uncertain-DC} \\
\hdashline
\rowcolor{Gray}
\multicolumn{2}{l}{\textbf{Logistics}} \\
\hdashline
$\delta^{(vc)}_i (\delta^{(vc)}_p) $ & \textbf{1}, if vaccine coming from manufacturer $i$ (type $p$) needs very cold or ultra cold refrigeration. & \ref{eq: Capacity with Forcing},\ref{eq: Capacity with Forcing2},\ref{eq: Capacity with Forcing3},\ref{eq: Forcing individually second echelon} \\
$\delta^{(uc)}_i (\delta^{(uc)}_p) $ & \textbf{1}, if vaccine coming from manufacturer $i$ (type $p$) needs ultra cold refrigeration. & \ref{eq: Capacity with Forcing2},\ref{eq: Capacity with Forcing3},\ref{eq: Forcing individually second echelon} \\
 $a_{k j}$ & \textbf{1}, if vaccination center (VC) $k$ is feasible to be assigned to distribution center (DC) $j$. & \ref{eq:DC-VC-linking} \\
 $a_{s k}$ & \textbf{1}, if population site (PS) $s$ is feasible to be assigned to vaccination center (VC) $k$. & \ref{eq: p-median and max coverage mixture} \\
 $\eta_{a b}$ & Distance from node $a$ to node $b$ (general index notation $a$ and $ b$) . & \ref{eq: p-median and max coverage mixture},\ref{eq: p-median2},\ref{eq: step-wise coverage},\ref{eq: step-wise detail coverage1},\ref{eq:env obj} \\
$\tau_{rr'}$ & Travel time from region $r$ to region $r'$. & \ref{eq: vehicle-routing} \\
$\tau^S_{r}$ & Service time of vehicle (mobile vaccination clinic) at region $r$. & \ref{eq: vehicle-routing} \\
 $\eta_\text{max}$ & Maximum allowable distance. & \ref{eq: p-median2},\ref{eq: step-wise detail coverage1} \\
  $\mathfrak{y}_q$ & Distance at coverage level $q$ & \ref{eq: step-wise detail coverage1} \\
 $e^v_{r g}$ & Estimated equitable quantity that should be allocated to region $r$ for priority group $g$. & \ref{eq: equity based objective3.2} \\
\hdashline
\rowcolor{Gray}
\multicolumn{2}{l}{\textbf{Others}} \\
\hdashline
$\varsigma$ & Weight such that $0 \leq \varsigma \leq 1$. & \ref{eq: equity based objective3.1} \\
$\mathbf{M}$ & Big-$M$ value. \\
$\check{\mathbf{M}}$ & Context specific penalty term (not necessarily a large value) or importance factor. & \ref{eq: robustify by std},\ref{eq: equity based objective3.1} \\
\hdashline
\rowcolor{Gray}
\multicolumn{2}{l}{\textbf{Vaccine}} \\
\hdashline
$\lambda^{sl}$ & Vaccine vial shelf life. & \ref{eq: Shelf life of vaccine vial 1},\ref{eq: Shelf life of vaccine vial 2},\ref{eq: Shelf life of vaccine vial 4},\ref{eq: Shelf life of vaccine vial 3},\ref{eq: Shelf life of vaccine vial 5} \\
$\tau$ & Vaccine safe-use-life after opening the vial. & \ref{eq:OPen Vial Wastage in forward time},\ref{eq:OPen Vial Wastage in forward time2},\ref{eq:OPen Vial Wastage in backward time},\ref{eq:OPen Vial Wastage in backward time2} \\
$\iota_{\nu}$ & Number of vaccines contained in a vial size $\nu \in \mathcal{V}$. & \ref{eq: vial vs dosage balancing} \\
$\iota_{\nu p}$ & Number of vaccines contained in a vial size $\nu \in \mathcal{V}$, for group $p$. & \ref{eq:OPen Vial Wastage in forward time},\ref{eq:OPen Vial Wastage in forward time2} \\
\hdashline
\rowcolor{Gray}
\multicolumn{2}{l}{\textbf{Workforce}} \\
\hdashline
$C^{HW}$ & Number of people a health worker can serve within time $t$. & \ref{eq: healthcare personnel availability} \\
$C^{EH}_{k t}$ & Existing workforce in vaccination center (VC) $k$ in time $t$. & \ref{eq: healthcare personnel availability} \\
\hdashline
\rowcolor{Gray}
\multicolumn{2}{l}{\textbf{Uncertainty Modeling}} \\
\hdashline
$\alpha$ & Significance level such that $0\leq \alpha \leq 1$. & \ref{eq:uncertain-DC},\ref{eq: CC normal} \\
$\mathsf{Z}_{1 - \alpha}$ & $\mathsf{Z}$ score of standard normal distribution at $1-\alpha$ confidence level. & \ref{eq:uncertain-DC},\ref{eq: CC normal}  \\
$\mu (\sigma)$ & Mean (standard deviation) estimate of an uncertain parameter. & \ref{eq:uncertain-DC},\ref{eq:DROvsSP},\ref{eq: CC normal}\\
$\epsilon$ & Ambiguity in empirical estimate of a parameter vaccination center (VC) $k$ in time $t$. & \ref{eq:DROvsSP} \\
$\underline{\epsilon}, \overline{\epsilon}$ & Factor to scale an ambiguously estimated parameter such that $0 \leq \underline{\epsilon} \leq 1 \leq \overline{\epsilon}$. & \ref{eq:DROvsSP} \\
$\mathsf{p}_{\omega}$ & Probability that a scenario $\omega \in \Omega$ would happen. &  \ref{eq: robustify by std} \\
\hline
\end{tabular}
\end{adjustbox}
\end{center}
\label{Tab:parameter description}
\end{table}

\begin{table} \setlength{\tabcolsep}{4pt}
\vspace{-.3cm}
\caption{Decision Variables}
\vspace{-0.5cm}
\begin{center}
\begin{adjustbox}{width=\textwidth}
\renewcommand{\arraystretch}{1.2}
\begin{tabular}{clll}
\hline
\rowcolor{Gray}
\textbf{Notation}&\textbf{Description}&\textbf{Type}&\textbf{Equations}\\
\hline
\rowcolor{Gray}
\multicolumn{2}{l}{\textbf{Allocation}} \\
\hdashline
$x_{j k}$ & \textbf{1}, if vaccination center (VC) $k$ is assigned to distribution center (DC) $j$. & B & \ref{eq:DC-VC-linking},\ref{eq: equity based objective 2} \\
$x_{s k}$ & Assignment variable (or percentage, Continuous). Assigning population site (PS) $s$ to vaccination center (VC) $k$. & B/C & \ref{eq: p-median and max coverage mixture} \\
$\mathfrak{z}_s$ & \textbf{1}, if population site (PS) $s$ is covered by a facility for vaccinating its population. & B & \ref{eq: p-median2} \\
$\mathfrak{z}_{sq}$ & \textbf{1}, if population site (PS) $s$ is covered by a facility at its distance based level $q \in \mathcal{Q}_s$. & B & \ref{eq: step-wise coverage},\ref{eq: step-wise detail coverage1} \\
$y^{dv}_{j k t}$ & \textbf{1}, if supply happens to vaccination center (VC) $k$ by distribution center (DC) $j$ at time period $t$. & B & \textit{In Text} \\
 $x^{v}_{g k p t}$ & Number of people of group $g$ vaccinated by vaccine type $p$ at vaccination center (VC) $k$ during time period $t$. & I/C & \ref{eq: Inventory_balance_VC},\ref{eq: healthcare personnel availability},\ref{eq: vial vs dosage balancing},\ref{eq: equity based objective} \\
 $x^v_{k p t}$ & Number of people vaccinated by vaccine type $p$, at vaccination center (VC) $k$ during time $t$. & I/C & \ref{eq:shortage-demand} \\
 $x^{v}_{g k p t t'}$ & Number of people of group $g$ vaccinated by vaccine type $p$ at vaccination center (VC) $k$ during time period $t'$ from a vial opened at time period $t$. & I/C & \ref{eq:OPen Vial Wastage in forward time},\ref{eq:OPen Vial Wastage in forward time2},\ref{eq:OPen Vial Wastage in backward time},\ref{eq:OPen Vial Wastage in backward time2} \\
$y^v_{g k p t}$ & Number of vaccine type $p$ available to fulfill demand or exactly meeting demand of priority group $g$ in vaccination center (VC) $k$ at time period $t$. & I/C & \ref{eq:priority group focus},\ref{eq:priority group focus2},\ref{eq:priority group focus3},\ref{eq:priority group focus4} \\
$s^v_{g k p t}$ & Amount of shortages of vaccine type $p$ in vaccination center (VC) $k$ at time period $t$ by group $g$. & I/C & \ref{eq:OPen Vial Wastage in backward time2},\ref{eq:env obj}\\
$s^v_{k p t}$ & Amount of shortages of vaccine type $p$ in vaccination center (VC) $k$ at time period $t$. & I/C & \ref{eq:shortage-demand} \\
$L^{v}_{k p t t'}$ & Number of type $p$ vaccine vials allocated to vaccination center (VC) $k$ at time $t$ and opened in time $t'$. & I/C & \ref{eq: Shelf life of vaccine vial 1},\ref{eq: Shelf life of vaccine vial 2},\ref{eq: Shelf life of vaccine vial 4},\ref{eq: Shelf life of vaccine vial 3} \\
$N^v_{k p t'}$ & Number of type $p$ vaccine vials opened at vaccination center (VC) $k$ in time $t'$. & I/C & \ref{eq: Shelf life of vaccine vial 4},\ref{eq: Shelf life of vaccine vial 3}\\
$N^v_{\nu k p t'}$ & Number of type $p$ vaccine vial of size $\nu$ opened at vaccination center (VC) $k$ in time $t'$. & I/C & \ref{eq: vial vs dosage balancing},\ref{eq:OPen Vial Wastage in forward time},\ref{eq:OPen Vial Wastage in forward time2} \\
$\mathfrak{x}_{\ell k}$ & \textbf{1}, if vaccination center (VC) $k$ served by an outreach center/team $\ell \in \mathcal{O}$. & B & \ref{eq: step-wise detail coverage1} \\
$A^{v}_{k p 0 t'}$ & Number of type $p$ vaccine vials available in vaccination center (VC) $k$ at initial time period and opened in time $t'$. & I/C & \ref{eq: Shelf life of vaccine vial 3},\ref{eq: Shelf life of vaccine vial 5}\\
$\mathtt{\Delta}^+_{r g}, \mathtt{\Delta}^-_{r g}$ & Positive and negative deviation from a quantity particular to a region $r$ for group $g$. & C & \ref{eq: equity based objective3.1},\ref{eq: equity based objective3.2} \\
$\mathtt{\Pi}_{r g}$ & Positive difference between minimum percentage of demand to be met and allocation in region $r$ for group $g$. & C & \ref{eq: equity based objective3.1} \\
$f_r (f_g)$ & Fraction of total population in region $r$ (of risk group $g$). & C & \ref{eq: herd-immunity}\\
\hdashline
\rowcolor{Gray}
\multicolumn{2}{l}{\textbf{Inventory}} \\
\hdashline
$I^t_j$ & Inventory available at distribution center (DC) $j$ at the end of time period $t$. & I/C & \ref{eq: Inventory_balance},\ref{eq: Inventory_balancee},\ref{eq: inventory-capacity},\ref{eq: Inventory_balance with lead time1}\\
$I^t_k$ & Inventory available at vaccination center (VC) $k$ at the end of time period $t$. & I/C & \ref{eq: Inventory_balance_VC},\ref{eq:priority group focus},\ref{eq:priority group focus2},\ref{eq:priority group focus3},\ref{eq:priority group focus4} \\ 
$I^0_{k p}$ & Number of type $p$ vaccine vials available at initial time period at vaccination center (VC) $k$. & I/C & \ref{eq: Shelf life of vaccine vial 5}\\
$Q_j$ & Order quantity of distribution center (DC) $j$. & I/C & \ref{eq:uncertain-DC} \\
\hdashline
\rowcolor{Gray}
\multicolumn{2}{l}{\textbf{Location}} \\
\hdashline
$Y^{D(c)}_j$ & \textbf{1}, if distribution center (DC) at location $j$ with usual cold storage facility $c$ is launched. & B & \ref{eq: Capacity with Forcing},\ref{eq: Forcing individually second echelon},\ref{eq: Forcing individually first echelon2},\ref{eq:env obj} \\
 $Y^{D(vc)}_j$ & \textbf{1}, if distribution center (DC) at location $j$ with very cold storage facility $vc$ is launched. & B & \ref{eq: Forcing individually second echelon},\ref{eq: Forcing individually first echelon2},\ref{eq:env obj} \\
 $Y^{D(uc)}_j$ & \textbf{1}, if ultra cold storage facility $uc$ is added to an \textit{existing} distribution center (DC) $j$. & B & \ref{eq: Capacity with Forcing2},\ref{eq: Capacity with Forcing3},\ref{eq: Forcing individually second echelon},\ref{eq:env obj} \\
 $Y^{D}_j$ & \textbf{1}, if distribution center (DC) $j$ having only usual cold storage facility $j$ is launched. & B & \ref{eq: outflow-location},\ref{eq: inflow-location},\ref{eq: inflow-locationn},\ref{eq: inventory-capacity},\ref{eq:env obj} \\
$Z^V_k$ & \textbf{1}, if vaccination center (VC) is established. & B & \ref{eq: p-median and max coverage mixture},\ref{eq: p-median2},\ref{eq: step-wise coverage},\ref{eq:env obj} \\
\hdashline
\rowcolor{Gray}
\multicolumn{2}{l}{\textbf{Routing}} \\
\hdashline
$y_{r r'}$ & \textbf{1}, if region $r$ and $r'$ are served in order in a vehicle routing problem. & B & \ref{eq: vehicle-routing} \\
$z_{r h}$ & \textbf{1}, if region $r$ is served by vehicle $h$. & B & \ref{eq: vehicle-routing} \\
$T_r$ & Service start time of region $r$. & C & \ref{eq: vehicle-routing} \\
$Z^R_r$ & \textbf{1}, if region $r$ is visited by a vehicle. & B & \textit{In Text} \\
\hdashline
\rowcolor{Gray}
\multicolumn{2}{l}{\textbf{Transshipment}} \\
\hdashline
$x^{m d}_{i j p t}$ & Number of vaccine-type $p$ arriving from manufacturer $i$ to distribution center (DC) $j$ at the beginning of time period $t$. & I/C & \ref{eq:capacity_man0},\ref{eq: Inventory_balance with lead time1},\ref{eq:env obj} \\
$x^{m d}_{i j p t-l}$ & Number of vaccine-type $p$ arriving from manufacturer $i$ to distribution center (DC) $j$ at the beginning of time period $t$, ordered at $t-l$. & I/C & \ref{eq: Inventory_balance},\ref{eq: inflow-location} \\
$x^{md}_{i j t' t}$ & Amount of vaccines shipped from manufacturer $i$ to distribution center (DC) $j$, ordered at time $t'$, received at time $t$. & I/C & \ref{eq: Capacity with Forcing},\ref{eq: Capacity with Forcing2},\ref{eq: Capacity with Forcing3},\ref{eq: Forcing individually first echelon2} \\
 $x^{d v}_{j k p t} (x^{d v}_{j k p t h})$ & Number of vaccine type $p$ transferred from distribution center (DC) $j$ to vaccination center (VC) $k$ at time $t$ (using vehicle $h$). & I/C & \ref{eq: Inventory_balance},\ref{eq: Inventory_balancee},\ref{eq: outflow-location},\ref{eq: Inventory_balance with lead time1},\ref{eq: Capacity with Forcing}, \\
 & & & \ref{eq: Capacity with Forcing2},\ref{eq: Capacity with Forcing3},\ref{eq: Forcing individually second echelon},\ref{eq: fleet size capacity},\ref{eq: Inventory_balance_VC}, \\
 & & & \ref{eq:DC-VC-linking},\ref{eq: Shelf life of vaccine vial 1},\ref{eq: Shelf life of vaccine vial 2},\ref{eq:priority group focus},\ref{eq:priority group focus2}, \\
 & & & \ref{eq:priority group focus3},\ref{eq:priority group focus4},\ref{eq: CC normal},\ref{eq:env obj} \\
 $x^{d d'}_{j j' p t}$ & Number of vaccine type $p$ transshipped from distribution center (DC) $j$ to another distribution center (DC) $j'$ at time period $t$. & I/C & \ref{eq: Inventory_balance},\ref{eq: Inventory_balancee},\ref{eq: outflow-location},\ref{eq:env obj}\\
 $x^{d' d}_{j' j p t}$ & Number of vaccine type $p$ arriving to distribution center (DC) $j$ from another distribution center (DC) $j'$ via transshipment at time period $t$. & I/C & \ref{eq: Inventory_balance},\ref{eq: Inventory_balancee},\ref{eq: inflow-location},\ref{eq: inflow-locationn}\\
  $x^{v v'}_{k k' p t}$ & Number of vaccine type $p$ transshipped from vaccination center (VC) $k$ to another vaccination center (VC) $j'$ at time period $t$. & I/C & \ref{eq:priority group focus},\ref{eq:priority group focus2},\ref{eq:env obj} \\
 $x^{v' v}_{k' k p t}$ & Number of vaccine type $p$ arriving to vaccination center (VC) $k$ from another vaccination center (VC) $k'$ via transshipment at time period $t$. & I/C & \ref{eq:priority group focus},\ref{eq:priority group focus2} \\
$X^m_{i t}$ & \textbf{1}, if selection of manufacturer $i$ at time $t$ for vaccine production occurs. & B & \ref{eq:capacity_man0} \\
 $X^m_{i p t' t}$ & \textbf{1}, if vaccine type $p$ ordered from manufacturer $i$ at time $t'$ to be delivered to distribution center (DC) at the beginning of time $t$. & B & \ref{eq: manufacturer restriction} \\
\hdashline
\rowcolor{Gray}
\multicolumn{2}{l}{\textbf{Vaccine Wastage}} \\
\hdashline
 $w_{j p t}$ & Wasted amount of type $p$ vaccine vial at distribution center (DC) $j$ at time $t$. & I/C & \ref{eq: Inventory_balance},\ref{eq: Inventory_balancee},\ref{eq: outflow-location},\ref{eq: Inventory_balance with lead time1}\\
$w_{k p t}$ & Wasted amount of type $p$ vaccine vial at vaccination center (VC) $k$ at time $t$. & I/C & \ref{eq: Inventory_balance_VC},\ref{eq: Shelf life of vaccine vial 1},\ref{eq: Shelf life of vaccine vial 2},\ref{eq:priority group focus},\ref{eq:priority group focus2}, \\
& & & \ref{eq:priority group focus3},\ref{eq:priority group focus4} \\
$W^{v}_{k p 0 t'}$ & Amount of type $p$ vaccine vials available in vaccination center (VC) $k$ at initial time period and wasted in time $t'$ due to time expiration. & I/C & \ref{eq: Shelf life of vaccine vial 5}\\
$w^{o v}_{k p t}$ & Wasted dosage of vaccine type $p$ in vaccination center (VC) $k$ at time period $t$ from opened vial $o\nu$. & I/C & \ref{eq: vial vs dosage balancing},\ref{eq:OPen Vial Wastage in backward time}\\
\hdashline
\rowcolor{Gray}
\multicolumn{2}{l}{\textbf{Workforce}} \\
\hdashline
$F^{R}_{k t}$ & Required number of health care worker in vaccination center (VC) $k$ at time period $t$. & I & \ref{eq: healthcare personnel availability} \\
$F^{A}_{k t}$ & Additional healthcare workers required in vaccination center (VC) $k$ at time period $t$. & I & \ref{eq: healthcare personnel availability} \\
\hdashline
\rowcolor{Gray}
\multicolumn{2}{l}{\textbf{Uncertainty Modeling}} \\
\hdashline
$\mathsf{p}_{\omega}$ & Probability that a scenario $\omega \in \Omega$ would happen. & C & \ref{eq:DROvsSP}\\
\hline
\end{tabular}
\end{adjustbox}
\end{center}
\label{Tab:Decision Variable Description}
\end{table}

\FloatBarrier

\bgroup

\begin{table}[!hbt] \setlength{\tabcolsep}{8pt}
\begin{center}
\caption{Decisions considered in each optimization article.}
\begin{adjustbox}{width=0.75\textwidth}
\renewcommand{\arraystretch}{0.95}
\begin{tabular}{l|ccccc|c}
\hline
\rowcolor{Gray}
& \multicolumn{5}{c}{\textbf{Type of Decision}} & \\
\rowcolor{Gray}
\textbf{Article} & \textbf{Location} & \textbf{Inventory} & \textbf{Allocation} & \textbf{Distribution} & \textbf{Routing} & \textbf{SIR} \\
\hline

\rowcolor{Gray}
\multicolumn{7}{c}{\textbf{\textcolor{black}{Supply Chain Models (See Section \S \ref{sec:model})}}} \\

\hline
\rowcolor{Gray}
\textbf{Resource Allocation} \\

\hline

\citet{abbasi2020modeling}  &     &     &  $\checkmark$  &  $\checkmark$  &     &    \\

\citet{anahideh2022fair}  &     &     &  $\checkmark$  &     &     &    \\

\citet{balcik2022mathematical}  &     &     &  $\checkmark$  &  $\checkmark$  &     &    \\

\citet{bandi2021optimal} &  &  &  $\checkmark$  &  $\checkmark$  &  &    \\

\citet{barthspatiotemporal} &  &  &  $\checkmark$  &  &  &  $\checkmark$ \\

\citet{bennouna2022covid}  &     &     &  $\checkmark$  &  $\checkmark$  &     &  $\checkmark$ \\

\citet{bertsimas2020optimizing}  &     &     &  $\checkmark$  &     &     &  $\checkmark$ \\

\citet{enayati2020optimal}  &     &     &  $\checkmark$  &  $\checkmark$  &     &  $\checkmark$ \\

\citet{hu2023first} &  &  &  $\checkmark$  &  $\checkmark$  &  &  $\checkmark$ \\

\citet{huang2017equalizing}  &     &     &  $\checkmark$  &     &     &    \\

\citet{jadidi2021two}  &     &     &  $\checkmark$  &     &     &  $\checkmark$ \\

\citet{jarumaneeroj2022epidemiology}  &     &     &  $\checkmark$  &     &     &  $\checkmark$ \\

\citet{liquantifying} &  &  &  $\checkmark$  &  $\checkmark$  &  &  $\checkmark$ \\

\citet{minoza2021covid}  &     &     &  $\checkmark$  &  $\checkmark$  &     &  $\checkmark$ \\

\citet{munguia2021fair}  &     &     &  $\checkmark$  &     &     &    \\

\citet{orgut2023equitable} &  &  &  $\checkmark$  &  &  &    \\

\citet{rao2021optimal}  &     &     &  $\checkmark$  &     &     &  $\checkmark$ \\

\citet{roy2021optimal}  &     &     &  $\checkmark$  &     &     &  $\checkmark$ \\

\citet{shukla2022optimizing}  &     &     &  $\checkmark$  &  $\checkmark$  &     &    \\

\citet{thul2023stochastic}  &     &     &  $\checkmark$  &     &     &  $\checkmark$ \\

\citet{yang2022comparison}  &     &     &  $\checkmark$  &     &     &  $\checkmark$ \\

\citet{yarmand2014optimal}  &     &     &  $\checkmark$  &     &     &    \\

\hline
\rowcolor{Gray}
\textbf{Inventory Management} \\

\hline

\citet{azadi2020optimization} &   & $\checkmark$ &   & $\checkmark$ &   &   \\

\citet{azadi2020developing} &   & $\checkmark$ &   &   &   &   \\

\citet{bonney2011environmentally} &   & $\checkmark$ &   &   &   &   \\

\citet{fadaki2022multi} &   & $\checkmark$ & $\checkmark$ & $\checkmark$ &   &   \\

\citet{georgiadis2021optimal} &   & $\checkmark$ & $\checkmark$ & $\checkmark$ &   &   \\

\citet{hovav2015network} &   & $\checkmark$ & $\checkmark$ & $\checkmark$ &   &   \\

\citet{icsik2023optimizing} &  &  $\checkmark$  &  $\checkmark$  &  $\checkmark$  &  &    \\

\citet{jahani2022covid} &   & $\checkmark$ & $\checkmark$ & $\checkmark$ &   &   \\

\citet{karakaya2023developing} &   & $\checkmark$ & $\checkmark$ &   &   &   \\

\citet{mak2022managing} &   & $\checkmark$ & $\checkmark$ &   &   & $\checkmark$ \\

\citet{mohammadi2022bi} &   & $\checkmark$ & $\checkmark$ & $\checkmark$ &   &   \\

\citet{shah2022investigation} &   & $\checkmark$ &   &   &   &   \\

\citet{sinha2021strategies} &   & $\checkmark$ & $\checkmark$ &   &   &   \\

\hline
\rowcolor{Gray}
\textbf{Supply Chain Management} \\

\hline

\citet{abbasi2023designing} &  $\checkmark$  &  &  $\checkmark$  &  $\checkmark$  &  &    \\

\citet{basciftci2023resource} & $\checkmark$ & $\checkmark$ & $\checkmark$ &   &   &   \\

\citet{chowdhury2022modeling} & $\checkmark$ & $\checkmark$ & $\checkmark$ & $\checkmark$ & $\checkmark$ &   \\

\citet{gilani2022data} & $\checkmark$ & $\checkmark$ & $\checkmark$ & $\checkmark$ &   &   \\

\citet{goodarzian2022sustainable} & $\checkmark$ & $\checkmark$ & $\checkmark$ & $\checkmark$ &   &   \\

\citet{goodarzian2022designing} & $\checkmark$ &   & $\checkmark$ & $\checkmark$ &   &   \\

\citet{habibi2023designing} & $\checkmark$ & $\checkmark$ & $\checkmark$ & $\checkmark$ & $\checkmark$ &   \\

\citet{kohneh2023optimization} & $\checkmark$ & $\checkmark$ & $\checkmark$ & $\checkmark$ &   &   \\

\citet{lai2021multi} & $\checkmark$ & $\checkmark$ & $\checkmark$ &   &   &   \\

\citet{li2021locate} & $\checkmark$ & $\checkmark$ & $\checkmark$ &   &   &   \\

\citet{lim2022redesign} & $\checkmark$ & $\checkmark$ & $\checkmark$ & $\checkmark$ &   &   \\

\citet{manupati2021multi} & $\checkmark$ & $\checkmark$ & $\checkmark$ & $\checkmark$ &   &   \\

\citet{rahman2023optimising} & $\checkmark$ & $\checkmark$ & $\checkmark$ & $\checkmark$ &   &   \\

\citet{rastegar2021inventory} & $\checkmark$ & $\checkmark$ & $\checkmark$ & $\checkmark$ &   &   \\

\citet{sazvar2021capacity} & $\checkmark$ & $\checkmark$ & $\checkmark$ & $\checkmark$ &   &   \\

\citet{shiri2022equitable} & $\checkmark$ & $\checkmark$ & $\checkmark$ & $\checkmark$ &   &   \\

\citet{tang2022bi} & $\checkmark$ & $\checkmark$ & $\checkmark$ &   &   &   \\

\citet{tavana2021mathematical} & $\checkmark$ & $\checkmark$ & $\checkmark$ & $\checkmark$ &   &   \\

\citet{wang2023robust} & $\checkmark$ & $\checkmark$ & $\checkmark$ & $\checkmark$ &   &   \\

\citet{xu2021hub} & $\checkmark$ &   & $\checkmark$ & $\checkmark$ &   &   \\

\citet{yang2021optimizing} & $\checkmark$ & $\checkmark$ & $\checkmark$ & $\checkmark$ &   &   \\

\hline
\rowcolor{Gray}
\multicolumn{7}{c}{\textbf{\textcolor{black}{Location and Routing Models (See Section \S \ref{sec: pure facility location model} and Appendix \S \ref{sec: vehicle routing})}}} \\

\hline
\rowcolor{Gray}
\textbf{Location-Allocation} \\

\hline

\citet{bertsimas2022locate} & $\checkmark$ &   & $\checkmark$ &   &   & $\checkmark$ \\

\citet{bravo2022optimal} & $\checkmark$ &   & $\checkmark$ &   &   &   \\

\citet{cabezas2021two} & $\checkmark$ &   & $\checkmark$ &   &   &   \\

\citet{cao2023digital} & $\checkmark$ &   & $\checkmark$ &   &   &   \\

\citet{dastgoshade2022social} & $\checkmark$ &   & $\checkmark$ &   &   &   \\

\citet{emu2021validating} & $\checkmark$ &   &   & $\checkmark$ &   &   \\

\citet{enayati2023multimodal} &  $\checkmark$  &  &   &  $\checkmark$  &  &    \\

\citet{enayati2023vaccine} &  $\checkmark$  &  &   &  $\checkmark$  &  &    \\

\citet{kumar2022optimal} & $\checkmark$ &   & $\checkmark$ &   &   &   \\

\citet{leithauser2021quantifying} & $\checkmark$ &   & $\checkmark$ &   &   &   \\

\citet{lim2016coverage} & $\checkmark$ &   & $\checkmark$ &   &   &   \\

\citet{luo2023service} & $\checkmark$ &   & $\checkmark$ &   &   &   \\

\citet{lusiantoro2022locational} & $\checkmark$ &   & $\checkmark$ &   &   &   \\

\citet{polo2015location} & $\checkmark$ &   & $\checkmark$ &   &   &   \\

\citet{soria2021proposal} & $\checkmark$ &   &   & $\checkmark$ &   &   \\

\citet{srivastava2021strengthening} & $\checkmark$ &   & $\checkmark$ &   &   &   \\

\citet{zhang2022mass} & $\checkmark$ &   & $\checkmark$ &   &   &   \\

\hline
\rowcolor{Gray}
\textbf{Routing} \\

\hline

\citet{shamsi2021novel} &   &   &   & $\checkmark$ & $\checkmark$ & $\checkmark$ \\

\citet{yang2021outreach} & $\checkmark$ &   &   & $\checkmark$ & $\checkmark$ &   \\

\citet{yucesoymobile} & $\checkmark$ &   &   & $\checkmark$ & $\checkmark$ &   \\

\hline
\end{tabular}
\end{adjustbox}
\label{table:optimization-based-classification}
\end{center}
\end{table}

\FloatBarrier

\newpage

\section{Supply Chain Models} \label{sec:model}

In this section, we present a comprehensive perspective on the diverse factors considered in the existing literature concerning the COVID-19 vaccine supply chain across multiple echelons. We achieve this by examining various objective functions and constraints that have been utilized for different components of the supply chain. Specifically, we emphasize key constraints derived from mathematical models documented in relevant literature, while also highlighting the potential challenges associated with modeling the intricate network of vaccine allocation and distribution for manufacturers (Section \S \ref{subsection:manufacturers}), distribution centers (Section \S \ref{subsection:DCs}), and  vaccination centers (Section \S \ref{subsection:VCs}).

\subsection{Manufacturers} \label{subsection:manufacturers}

In this section, we focus on the production capacity constraint applicable to manufacturers in the vaccine supply chain. Several articles have addressed inventory management \citep{georgiadis2021optimal, hovav2015network, jahani2022covid, shah2022investigation} and supply chain management \citep{chowdhury2022modeling, gilani2022data, goodarzian2022sustainable, manupati2021multi, sazvar2021capacity, tavana2021mathematical} aspects related to manufacturers' requirements. In the vaccine supply chain, capacity constraints can exist between the manufacturer and the distribution center (DC), as well as between the distribution center (DC) and the vaccination center (VC). These capacity constraints can pertain to various resources, such as the number of vaccines, the number of vaccine lots, fleet size, workforce, time, budget, and other restricted resources necessary for production. Additionally, the capacity constraints can be specific to commodities, time, or facilities. Moreover, lead time considerations may arise due to the time required for production.
In its simplest form, the capacity constraint can be expressed as in \eqref{eq:capacity_man0}. 

\noindent \begin{align}
\label{eq:capacity_man0}
& \underbrace{\sum_{j \in \mathcal{N}} x^{md}_{i j p t} \leq C^M_{i p t} X^m_{it}~, \hspace{0.25cm}  \forall \;i,p,t}_{\substack{\text{Product-and time-specific production capacity} \\ \text{constraint of a manufacturer.}}}
\end{align}

The capacity constraint can take various forms depending on the specific problem context. To achieve a more aggregated representation of the constraints in terms of vaccine type, time, and facility, the decision variables, parameters, summations, and set indexes in \eqref{eq:capacity_man0} can be redefined. This allows for a more generalized formulation of the capacity constraint, accommodating different variations based on the problem.

In practice, it is common for multiple distribution centers (DCs) to request vaccines from the same manufacturer. However, due to the manufacturer's capacity limitations, it can only fulfill the requests of a distribution center (DC) once its current order is processed. To address this consideration in the optimization model, two time indices ($t'$ representing the ordering time and $t$ representing the receiving time $t' \leq t$) are used (e.g., see \citet{tavana2021mathematical}).

\noindent \begin{align}
& \underbrace{\sum_{t' = 1}^{|\mathcal{T}|} X^m_{i p t' t} \leq 1, \hspace{0.25cm} \forall i, p, t}_{\substack{\text{Manufacturer can receive at most one order} \\ \text{for the same delivery time $t$ over all DCs.}}}; \quad \underbrace{\sum_{t' <\hat{t'} < t} \sum_{\hat{t'} \leq \hat{t}} X^m_{i p \hat{t'} \hat{t}} \leq \mathbf{M}\cdot (1 - X^m_{i p t' t} ), \hspace{0.25cm} \forall i, p,t, t' \leq t}_{\substack{\text{Manufacturer can not accept a new order} \\ \text{from a DC unless its previous order is served.}}} \label{eq: manufacturer restriction}
\end{align}

While the majority of research on COVID-19 focuses on shipments from manufacturers to distribution centers (DCs), it is worth noting that manufacturers can also directly transship vaccines to vaccination centers (VCs). Additionally, in the context of COVID-19, all manufacturers typically produce a single type of vaccine. For example, Johnson \& Johnson plants exclusively manufacture the Janssen Vaccine, while Pfizer plants exclusively produce the BioNTech Vaccine. However, in cases where a manufacturing plant produces multiple types of vaccines, the set $p \in \mathcal{P}_i$ can be incorporated into \eqref{eq:capacity_man0} and \eqref{eq: manufacturer restriction} to account for this variation. Furthermore, some articles consider the flow of raw materials from suppliers to manufacturers \citep{chowdhury2022modeling, gilani2022data, sazvar2021capacity}, the manufacturing process itself \citep{shah2022investigation}, and the location decisions of manufacturers with different technology levels \citep{sazvar2021capacity}. These aspects contribute to a more comprehensive understanding of the COVID-19 vaccine supply chain.

\subsection{Distribution Center} \label{subsection:DCs}

In this section, we present the considerations related to distribution centers (DCs) in terms of flow conservation (Section \S \ref{subsubsec:flow-conservation}), cold chain requirements (Section \S \ref{subsubsection:cold-chain-requirement}), and fleet size (Section \S \ref{subsubsection:fleet-size}). These considerations are based on optimization models proposed in recent COVID-19 literature addressing resource allocation \citep{abbasi2020modeling, anahideh2022fair, balcik2022mathematical, bennouna2022covid, bertsimas2020optimizing, enayati2020optimal, huang2017equalizing, jadidi2021two, jarumaneeroj2022epidemiology, minoza2021covid, munguia2021fair, rao2021optimal, roy2021optimal, shukla2022optimizing, thul2023stochastic, yang2022comparison, yarmand2014optimal}, inventory management \citep{azadi2020optimization, bonney2011environmentally,azadi2020developing,fadaki2022multi, georgiadis2021optimal, hovav2015network, karakaya2023developing, mak2022managing, mohammadi2022bi, shah2022investigation, sinha2021strategies}, and supply chain management \citep{basciftci2023resource, chowdhury2022modeling, gilani2022data, goodarzian2022sustainable, habibi2023designing,goodarzian2022designing, kohneh2023optimization, lim2022redesign, manupati2021multi, rahman2023optimising, rastegar2021inventory, sazvar2021capacity, shiri2022equitable, tavana2021mathematical, wang2023robust, xu2021hub, yang2021optimizing}. These models contribute to a comprehensive understanding of the distribution center (DC) aspects within the COVID-19 vaccine supply chain.

\subsubsection{Flow Conservation} \label{subsubsec:flow-conservation}

Distribution centers (DCs) play a critical role in the vaccine supply chain, facilitating both inflows and outflows of vaccines. The inflows to distribution centers (DCs) include procurement from manufacturers, their own inventory carried over from previous time periods, and shipments received from other distribution centers (DCs). Conversely, the outflows involve transshipments to vaccination centers (VCs), maintaining inventory for future time periods, shipments to other distribution centers (DCs), and addressing vaccine waste resulting from expiration, damage (e.g., physical or temperature-related), or mishandling. If distribution centers (DCs) can also function as vaccination centers (VCs), an additional vaccine outflow can be the administration of vaccines to individuals.

To maintain flow conservation in distribution centers (DCs) \eqref{eq: Inventory_balance} to \eqref{eq: inventory-capacity} are utilized. Equations \eqref{eq: Inventory_balance} and \eqref{eq: Inventory_balancee} represent the inventory flow conservation constraints for distribution centers (DCs). Equation \eqref{eq: outflow-location} ensures that if there is an outflow from a distribution center (DC), it must correspond to an existing outflow location. Similarly, \eqref{eq: inflow-location} and \eqref{eq: inflow-locationn} ensure that if there is an inflow to a distribution center (DC), it must originate from an existing inflow location. It should be noted that not all distribution centers (DCs) may serve all vaccination centers (VCs). In such cases, the flow balance constraints can be modified by incorporating the index set $\mathcal{V}_j$, which represents the set of vaccination centers (VCs) exclusively served by distribution center (DC) $j$.

Furthermore, in practical scenarios, there may exist a lead time ($l$) between the ordering and receiving of vaccines from one location to another. This lead time can be either deterministic or uncertain. In the current formulation, let us consider the case where there is only a deterministic lead time between the manufacturer and the distribution center (DC).
For $l = 0$, the formulation has no lead time. In such a case, \eqref{eq: Inventory_balancee} and \eqref{eq: inflow-locationn} become redundant and the decision variable $x^{md}_{i j p t-l}$ becomes $x^{md}_{i j p t}$ in \eqref{eq: Inventory_balance} and \eqref{eq: inflow-location}.

\noindent \begin{align}
& I^{t}_j = \underbrace{I^{t-1}_j + \sum_{i = 1}^{|\mathcal{M}|} \sum_{p = 1}^{|\mathcal{P}|} x^{m d}_{i j p t-l} + \sum_{\substack{j' \in \mathcal{N} \\ j' \neq j}} \sum_{p = 1}^{|\mathcal{P}|} x^{d' d}_{j' j p t}}_\text{Vaccine inflow.} \underbrace{- \sum_{k = 1}^{|\mathcal{V}|} \sum_{p = 1}^{|\mathcal{P}|} x^{d v}_{j k p t} - \sum_{\substack{j' \in \mathcal{N} \\ j' \neq j}} \sum_{p = 1}^{|\mathcal{P}|} x^{d d'}_{j j' p t} - \sum_{p=1}^{|\mathcal{P}|} w_{jpt}}_\text{Vaccine outflow.}, & \forall j,t > l \label{eq: Inventory_balance} \\
& I^{t}_j = \underbrace{I^{t-1}_j + \sum_{\substack{j' \in \mathcal{N} \\ j' \neq j}} \sum_{p = 1}^{|\mathcal{P}|} x^{d' d}_{j' j p t}}_\text{Vaccine inflow.} \underbrace{- \sum_{k = 1}^{|\mathcal{V}|} \sum_{p = 1}^{|\mathcal{P}|} x^{d v}_{j k p t} - \sum_{\substack{j' \in \mathcal{N} \\ j' \neq j}} \sum_{p = 1}^{|\mathcal{P}|} x^{d d'}_{j j' p t} - \sum_{p=1}^{|\mathcal{P}|} w_{jpt}}_\text{Vaccine outflow.}, & \forall j,t \leq l \label{eq: Inventory_balancee} \\
& \underbrace{\sum_{p=1}^{|\mathcal{P}|} \sum_{k = 1}^{|\mathcal{V}|} x^{dv}_{j k p t}}_\text{Outflow to VCs.} + \underbrace{\sum_{\substack{j' \in \mathcal{N} \\ j' \neq j}} \sum_{p = 1}^{|\mathcal{P}|} x^{d d'}_{j j' p t}}_\text{Outflow to other DCs.} + \underbrace{\sum_{p=1}^{|\mathcal{P}|} w_{jpt}}_\text{Wasted vaccines.} \leq \mathbf{M}\cdot Y^{D}_{j}, & \forall j,t \label{eq: outflow-location} \\
& \underbrace{\sum_{i = 1}^{|\mathcal{M}|} \sum_{p = 1}^{|\mathcal{P}|} x^{m d}_{i j p t-l}}_\text{Inflow from manufacturers.} + \underbrace{\sum_{\substack{j' \in \mathcal{N} \\ j' \neq j}} \sum_{p = 1}^{|\mathcal{P}|} x^{d' d}_{j' j p t}}_\text{Inflow from other DCs.} \leq \mathbf{M}\cdot Y^{D}_{j}, & \forall j,t > l \label{eq: inflow-location} \\
& \underbrace{\sum_{\substack{j' \in \mathcal{N} \\ j' \neq j}} \sum_{p = 1}^{|\mathcal{P}|} x^{d' d}_{j' j p t}}_\text{Inflow from other DCs.} \leq \mathbf{M}\cdot Y^{D}_{j}, & \forall j,t \leq l \label{eq: inflow-locationn} \\
& \underbrace{I^{t}_j \leq C^D_j Y^d_{j}, \hspace{0.25cm} \forall j,t}_\text{Inventory capacity constraint.}; \quad \underbrace{I^{t}_j \geq S^D_j Y^d_{j}, \hspace{0.25cm} \forall j,t}_\text{Inventory safety stock constraint.} \label{eq: inventory-capacity} 
\end{align}

\noindent \citet{georgiadis2021optimal} formulate an aggregated version of the inventory balancing constraint with lead time ($l$) using \eqref{eq: Inventory_balance with lead time1}: 

\noindent \begin{align} \label{eq: Inventory_balance with lead time1}
& \underbrace{\sum_{t' \leq t} \sum_{k = 1}^{|\mathcal{V}|} x^{dv}_{j k p t'}}_{\substack{\text{Outflow from a DC} \\ \text{until time $t$.}}} \leq \underbrace{I^0_j}_\text{Initial inventory.} + \underbrace{\sum_{t'' \leq (t - l)^+} \sum_{i = 1}^{|\mathcal{M}|} x^{md}_{i j p t''}}_{\substack{\text{Quantity of vaccines ordered at or before} \\ \text{time $(t-l)^{+}$, and arrived at or before time $t$.}}} - \underbrace{\sum_{t'' \leq (t - l)^+} w_{jt''p}}_\text{Vaccine wastage.}, & \forall j, t, p
\end{align}

In \citet{fadaki2022multi}, flow units of different decision variables or parameters differ. For example, if inflow is in terms of vaccine lots and outflow is in terms of vaccine vials, corresponding variables must be multiplied by the appropriate conversion coefficient. Moreover, \citet{hinojosa2000multiperiod} consider fractional flow (i.e., using continuous decision variables instead of discrete decision variables).

\subsubsection{Cold Storage Refrigeration Requirement} \label{subsubsection:cold-chain-requirement}

Vaccines may have varying storage requirements in terms of temperature, e.g., cold, very cold, or ultra cold. Ultra cold storage satisfies very cold storage requirements, but not the other way around. Additionally, distribution centers (DCs) can have different capacities for each level of refrigeration, which can vary over time. For example, \citet{tavana2021mathematical} employ \eqref{eq: Capacity with Forcing} to maintain the capacity restriction for very cold storage, \eqref{eq: Capacity with Forcing2} to sustain the capacity restriction for  very cold storage, and \eqref{eq: Capacity with Forcing3} to ensure the capacity restriction for ultra cold storage at distribution center (DC) $j$ during time period $t$. These equations assume that everything received by a distribution center (DC) will be distributed to vaccination centers (VCs) during the same time duration.

\noindent \begin{align}
& \sum_{i=1}^{|\mathcal{M}|} \sum_{t'} x^{md}_{i j t' t} (1 - \delta^{(vc)}_i) = \sum_{p=1}^{|\mathcal{P}|} \sum_{k = 1}^{|\mathcal{V}|} x^{d v}_{j k p t} (1 - \delta^{(vc)}_p)  \leq C^{D(c)}_j Y^{D(c)}_{j}, & \forall j,t \label{eq: Capacity with Forcing} \\
& \sum_{i=1}^{|\mathcal{M}|} \sum_{t'} x^{md}_{i j t' t} \delta^{(vc)}_i (1 - \delta^{(uc)}_i) = \sum_{p=1}^{|\mathcal{P}|} \sum_{k = 1}^{|\mathcal{V}|} x^{d v}_{j k p t} \delta^{(vc)}_p (1 - \delta^{(uc)}_p)  \leq C^{D(vc)}_j -  C^{D(uc)}_j Y^{D(uc)}_{j}, & \forall j,t \label{eq: Capacity with Forcing2} \\
& \sum_{i=1}^{|\mathcal{M}|} \sum_{t'} x^{md}_{i j t' t} \delta^{(vc)}_i \delta^{(uc)}_i = \sum_{p=1}^{|\mathcal{P}|} \sum_{k = 1}^{|\mathcal{V}|} x^{d v}_{j k p t} \delta^{(vc)}_p \delta^{(uc)}_p  \leq  C^{D(uc)}_j Y^{D(uc)}_{j}, & \forall j,t \label{eq: Capacity with Forcing3} 
\end{align}

\noindent The binary location decisions also need to be linked to the outflow and inflow from a distribution center (DC). Equation \eqref{eq: Forcing individually second echelon} links the outflow to the location decision variable, based on the required refrigeration type. 

\noindent \begin{align}
& x^{d v}_{j k p t} (1 - \delta^{(vc)}_p) \leq \mathbf{M}\cdot Y^{D(c)}_{j}, \hspace{0.25cm} \forall j, k, p, t; \quad x^{d v}_{j k p t} \delta^{(vc)}_p (1 - \delta^{(uc)}_p)  \leq  \mathbf{M}\cdot Y^{D(vc)}_{j}, \hspace{0.25cm} \forall j, k, p, t \label{eq: Forcing individually second echelon} \\ 
& x^{d v}_{j k p t} \delta^{(vc)}_p \delta^{(uc)}_p  \leq  \mathbf{M}\cdot Y^{D(uc)}_{j}, \hspace{0.25cm} \forall j, k, p, t \nonumber
\end{align} 

\noindent Similarly, \eqref{eq: Forcing individually first echelon2} ensures the same relation for inflow respecting the assumption that unlike cold and very cold storage type, no new distribution centers (DCs) will be set up for ultra cold storage; ultra cold storage facility can only be added to an existing distribution center (DC).

\noindent \begin{align} \label{eq: Forcing individually first echelon2}
& x^{md}_{i j t' t} \leq \mathbf{M}\cdot Y^{D(c)}_{j}, \hspace{0.25cm} \forall i, j, t', t; \quad x^{md}_{i j t' t} \leq \mathbf{M}\cdot Y^{D(vc)}_{j}, \hspace{0.25cm} \forall i, j, t', t
\end{align}

\subsubsection{Fleet and Drones} \label{subsubsection:fleet-size}
Vaccine vials can be transported between different echelons using various refrigerated vehicles, each with its own capacity for different levels of refrigeration \citep{goodarzian2022sustainable,dastgoshade2022social,chowdhury2022modeling,habibi2023designing,hovav2015network,lim2022redesign}. However, it is worth noting that the majority of existing studies in the literature do not consider modeling of vehicles and exclude them from their scope. Among the works that do incorporate vehicles, refrigerated trucks are commonly utilized as the primary transportation mode. Similarly, \citet{lim2022redesign} explore the use of motorbikes.

In practice, each distribution center (DC) may have its own dedicated fleet of vehicles, which can be employed to serve different vaccination centers (VCs). Detailed considerations regarding vehicle routing are presented in Section \S \ref{sec: vehicle routing}. Equation \eqref{eq: fleet size capacity} provides an aggregated representation of fleet availability across all distribution centers (DCs), capturing the overall capacity of the available vehicles.
 
\noindent \begin{align} \label{eq: fleet size capacity}
& \sum_{j = 1}^{|\mathcal{V}|} \sum_{p = 1}^{|\mathcal{P}|} \sum_{k = 1}^{|K|} x^{dv}_{j k p t} \leq C^F, \hspace{0.25cm} \forall t 
\end{align}

\noindent \citet{lim2016coverage} use vehicle trips required to transport vaccines in a replenishment as a decision variable for each specific inventory replenishment frequency. Moreover, \citet{yang2021optimizing} consider transportation capacity-based required trips for two different annual inventory replenishment frequencies. Travel time to complete a trip is an important parameter in this setting.

\srevision{Notable studies by \citet{wang2023robust}, \citet{enayati2023multimodal}, and \citet{enayati2023vaccine} explore the use of drones for vaccine distribution. While \citet{wang2023robust} optimizes drone capacity with vaccination scheduling in facility location problems, models in the latter two works involve location decision of drone base, drone relay station, and distribution centers (DCs) layering each possible origin-destination path for drone delivery with controlled drone stops, drone range,  and cold-time.}

\subsection{Vaccination Center} \label{subsection:VCs}

In this section, we discuss the vaccination center (VC) considerations in terms of flow conservation (Section \S \ref{VCsubsubsection:flow}), DC-VC assignment (Section \S \ref{sec:assigning a VC to a DC}), healthcare personnel capacity/scheduling (Section \S \ref{sec: healthcare personnel capacity}), vaccine shelf life (Section \S \ref{sec: vaccine_shelf_life}), vaccine vials (Section \S \ref{VCsubsubsection:vaccine-vials}), and priority groups (Section \S \ref{VCsubsubsection:priority-group}), based on the optimization models presented in the recent COVID-19 literature.

\subsubsection{Inventory and Flow Balancing} \label{VCsubsubsection:flow}
Vaccination centers (VCs), similar to distribution centers (DCs), experience both inflows and outflows of vaccines. The inflows of vaccines at vaccination centers (VCs) can originate from several sources, including procurement from manufacturers, their own inventory carried over from previous time periods, and shipment from other vaccination centers (VCs) (lateral transshipment) or distribution centers (DCs). Conversely, the outflows from vaccination centers (VCs) involve shipment to other vaccination centers (VCs), maintaining inventory for future periods, returning leftover vaccines to distribution centers (DCs), managing vaccine waste (such as open vial wastage or improper administration), and administering vaccines to individuals. In the context of this article, we assume that a vaccination center (VC) can only have inventory, vaccination, and waste as outflows. Additionally, the center can only receive inflows from its own inventory from the previous periods and distribution centers (DCs). It is important to note that vaccine demand can be either deterministic or stochastic. Furthermore, when managing lateral shipment to minimize open vial wastage, the flow should be considered at both the dosage and vial levels. In this formulation, the demand is expressed in terms of vaccine vials. If the demand were stated in terms of vaccine doses, appropriate conversion coefficients would need to be included for accurate calculations.

Equation \eqref{eq: Inventory_balance_VC} is the inventory flow conservation constraint for vaccination centers (VCs), assuming there is no supply shortage. Similar to flow conservation constraints at distribution center (DC) given in constraints \eqref{eq: outflow-location} to \eqref{eq: inventory-capacity}, constraints for location-flow binding, inventory capacity and safety stock level can be formulated for vaccination centers (VCs). Section \S \ref{sec: pure facility location model} further elaborates on facility location decision for vaccination centers (VCs) and the routing dynamics when a vaccination center (VC) is mobile.

\noindent \begin{align} 
& I^{t}_k = I^{t-1}_k + \underbrace{\sum_{j = 1}^{|\mathcal{N}|} \sum_{p = 1}^{|\mathcal{P}|} x^{d v}_{j k p t}}_\text{Inflow from DCs.} - \underbrace{\sum_{p = 1}^{|\mathcal{P}|} d_{k p t}}_\text{Vaccine demand.} - \underbrace{\sum_{p = 1}^{|\mathcal{P}|} w_{kpt}}_\text{Wasted vaccines.}, \hspace{0.25cm} \forall k,t \label{eq: Inventory_balance_VC}
\end{align}

\noindent Under supply scarcity, \eqref{eq: Inventory_balance_VC} would create infeasibility, therefore $d_{k p t}$ must be replaced by decision variable $x^v_{k p t}$ that denotes allocation decisions. In addition, \eqref{eq:shortage-demand} must be added to ensure that shortages are accounted for. 
\noindent \begin{align} \label{eq:shortage-demand}
& \underbrace{x^v_{k p t}}_\text{Vaccination.} + \underbrace{s^v_{k p t}}_\text{Shortage.} = \underbrace{d_{k p t}}_\text{Demand.}, \hspace{0.25cm} \forall k,p,t
\end{align}
Alternatively, the constraint $x^v_{k p t} \leq d_{k p t} \; \forall k,p,t$ can be used if tracking shortages is not necessary. \citet{hovav2015network} penalize the vaccine shortages using a shortage cost $c^{SH}_{k p}$ in the objective function. This practice is very common in works that consider uncertain demand. 

\subsubsection{Assigning a Vaccination Center to a Distribution Center} \label{sec:assigning a VC to a DC}

The linkage between distribution centers (DCs) and vaccination centers (VCs) can be established in multiple ways. Several articles in the literature explore different assignment methods, including one-to-one assignment, set covering assignment, and set packing assignment. In the case of the flow of vaccines from a distribution center (DC) to its corresponding vaccination center (VC), this can be modeled using \eqref{eq:DC-VC-linking}. It is worth noting that similar types of assignments can also be considered for the linkage between manufacturers and distribution centers (DCs) \citep{manupati2021multi} or between vaccination centers (VCs) themselves \citep{fadaki2022multi}.

\noindent \begin{align} \label{eq:DC-VC-linking}
\underbrace{\sum_{p, t} x^{d v}_{j k p t} \leq \mathbf{M}\cdot x_{j k}, \hspace{0.25cm} \forall j, k}_\text{Direct pairwise assignment.}; \quad \underbrace{\sum_{j} a_{k j} x_{j k} \geq 1, \hspace{0.25cm} \forall k}_\text{Set cover type constraint.}; \quad \underbrace{\sum_{j} a_{k j} x_{j k} \leq 1, \hspace{0.25cm} \forall k}_\text{Set packing type constraint.}
\end{align}

Once the distribution center (DC) assignment is specified, the next step is to allocate vaccines from the distribution center (DC) to various destinations, such as regional depots \citep{balcik2022mathematical,huang2017equalizing,munguia2021fair}, vaccination centers (VCs) \citep{abbasi2020modeling,anahideh2022fair,bertsimas2020optimizing,roy2021optimal,shukla2022optimizing,thul2023stochastic,yang2022comparison,yarmand2014optimal}, or specific sub-population groups \citep{bennouna2022covid,enayati2020optimal,jadidi2021two,jarumaneeroj2022epidemiology,minoza2021covid,rao2021optimal}. Different allocation rules and prioritization strategies are used in these models to determine the amount of vaccine to be allocated, with the objective of minimizing the spread of infection while considering the limited availability of vaccines. In addition to efficiency concerns, equity among the served sub-populations is often a key objective or constraint in these models. For more detailed information on \textit{resource allocation} articles, refer to Table~\ref{table:resource-allocation-articles} in Appendix.

\subsubsection{Healthcare Personnel Capacity and Vaccination Scheduling} \label{sec: healthcare personnel capacity} 
The healthcare workforce plays a crucial role in pandemic management, as skilled personnel are needed to administer vaccines to individuals. The allocation of healthcare workforce is addressed in several COVID-19 articles \citep{leithauser2021quantifying,georgiadis2021optimal,kumar2022optimal,lai2021multi,rahman2023optimising}. The availability and sustainability of healthcare personnel in vaccination centers (VCs) are captured by \eqref{eq: healthcare personnel availability}.

\noindent \begin{align}
& \underbrace{\sum_{g, p} x^v_{g k p t} \leq F^{R}_{k t}C^{HW}, \hspace{0.25cm} \forall k, t}_\text{Healthcare workforce availability constraint.}; \quad \underbrace{F^{R}_{k t} \leq F^{A}_{k t} + C^{EH}_{k t}, \hspace{0.25cm} \forall k, t}_\text{Additional healthcare workforce constraint.} \label{eq: healthcare personnel availability}
\end{align}

\citet{wang2023robust,zhang2022mass} focus on vaccination scheduling. \citet{zhang2022mass} propose a mass vaccination appointment scheduling model that has the potential to be extended to healthcare personnel scheduling, including overtime requirements. Their linear ordering formulation considers location, population size, and time window inputs to determine optimal vaccination center (VC) location, appointment acceptance, assignment of appointments to selected vaccination centers (VCs), and the scheduling of the vaccination timetable while accounting for possible tardiness. The formulation is further enhanced using dominance rules and valid inequalities, and is solved using a logic-based Benders approach. \srevision{\citet{bandi2021optimal} emphasize the importance of efficient appointment slot management, second-dose reserves, and flexible vaccine preferences for optimizing vaccination rates during a two-dose roll-out.}

\subsubsection{Vaccine Vial Shelf Life} \label{sec: vaccine_shelf_life}

Vaccine vials are considered perishable and have a limited shelf life (denoted by $\lambda^{sl}$). Once the shelf life of a vaccine vial expires, the vaccine doses within the vial become unsafe for use, resulting in waste. Many articles do not explicitly account for the shelf life of vaccine vials and instead define the problem time horizon to be shorter than the shelf life. However, \citet{georgiadis2021optimal} address this issue by introducing \eqref{eq: Shelf life of vaccine vial 1} to \eqref{eq: Shelf life of vaccine vial 5} to ensure the proper flow of vaccine vial shelf life. These equations aim to balance the forward flow of vaccine vial shelf life by fixing the receiving time ($t$) (see \eqref{eq: Shelf life of vaccine vial 1} and \eqref{eq: Shelf life of vaccine vial 2}) and the backward flow by fixing the opening time ($t'$) (see \eqref{eq: Shelf life of vaccine vial 3} and \eqref{eq: Shelf life of vaccine vial 4}). It is assumed that there is no lead time in transhipped vaccines between distribution centers (DCs) and vaccination centers (VCs). Other key assumptions include the shelf life applying to a vial immediately after its receiving time and a vial being opened one day after its receiving time. Thus, the relationship $t + 1 \leq t' \leq t + \lambda^{sl}$ holds, where $(t' - \lambda^{sl})^{+} \leq t \leq t' - 1$, and all times are bounded by the time horizon $|\mathcal{T}|$. Additionally, in \eqref{eq: Shelf life of vaccine vial 1}, it is assumed that all vial openings occur within the planning time horizon.

\noindent \begin{align}
& \underbrace{\sum_{t' = t + 1}^{t + \lambda^{sl}} L^v_{k p t t'}}_{\substack{\text{Vaccines can be opened at $t'$} \\ \text{where $t + 1 \leq t' \leq t + \lambda^{sl}$.}}} = \underbrace{\sum_{j \in \mathcal{N}} x^{d v}_{j k p t}}_\text{Vaccines received at $t$.} - \underbrace{w_{k p t}}_\text{Wasted vaccines.}, & \forall k, p , t \in [|\mathcal{T}| - \lambda^{sl}] \label{eq: Shelf life of vaccine vial 1} \\
& \underbrace{\sum_{t' = t+1}^{|\mathcal{T}|} L^v_{k p t t'} \leq \sum_{j \in \mathcal{N}} x^{d v}_{j k p t} - w_{k p t}}_{\substack{\text{Similar to Constraint \eqref{eq: Shelf life of vaccine vial 1}, by not considering} \\ \text{the vaccines opened after time horizon $|\mathcal{T}|$.}}}, & \forall k, p , |\mathcal{T}| - \lambda^{sl} < t \leq |\mathcal{T}| \label{eq: Shelf life of vaccine vial 2} \\
& \underbrace{N^v_{k p t'}}_{\substack{\text{Number of vaccines} \\ \text{opened at $t'$.}}} = \underbrace{\sum_{t' - \lambda^{sl} \leq t}^{t' - 1} L^v_{k p t t'}}_{\substack{\text{Number of vaccines received} \\ \text{in $t' - \lambda^{sl} \leq t \leq t' - 1$.}}}, & \forall k, p, t' > \lambda^{sl} \label{eq: Shelf life of vaccine vial 4} \\
& \underbrace{N^v_{k p t'} = {A^{v}_{k p 0 t'}} +  \sum_{(t' - \lambda^{sl})^+ \leq t}^{t' - 1} L^v_{k p t t'}}_{\substack{\text{Similar to Constraint \eqref{eq: Shelf life of vaccine vial 4}, by considering} \\ \text{initial net vaccine allocation.}}}, & \forall k, p, t' \leq \lambda^{sl} \label{eq: Shelf life of vaccine vial 3} \\
& \sum_{1 \leq t' \leq \lambda^{sl}} (\underbrace{A^{v}_{k p 0 t'}}_{\substack{\text{Initial net vaccine} \\ \text{availability.}}} + \underbrace{W^v_{k p 0 t'}}_{\substack{\text{Initial vaccine} \\ \text{waste.}}}) = \underbrace{I^0_{k p}}_{\substack{\text{Initial gross vaccine} \\ \text{inventory.}}}, & \forall k, p \label{eq: Shelf life of vaccine vial 5}
\end{align}

\subsubsection{Relationship Between Vaccine Vials, Dosages and Open Vial Wastage} \label{VCsubsubsection:vaccine-vials}

Several articles employ inventory replenishment policies in multi-echelon systems \citep{fadaki2022multi,georgiadis2021optimal,hovav2015network} to optimize vaccination delivery and minimize waste, considering the unique requirements of considered vaccines. These requirements encompass aspects such as cold chain management, storage conditions, vaccine-specific vial sizes at distribution (DCs) and vaccination centers (VCs), and expiration timelines. It is important to note that a single vaccine vial can contain multiple vaccine doses. In general, when a vaccine vial is opened, all doses within the vial should be administered to individuals before any doses expire and go to waste. It is worth mentioning that the open vial shelf life is typically shorter than the shelf life of the unopened vial. For instance, the Pfizer COVID-19 vaccine is packaged in vials containing five doses. While the vial shelf life is 30 days, the open vial shelf life is only 24 hours. This consideration underscores the need to efficiently utilize the doses within an opened vial to prevent wastage. Equation \eqref{eq: vial vs dosage balancing} ensures vial-dosage balancing, by equating the total number of opened doses to the summation of vaccinations among all priority groups and open vial wastage.

\noindent \begin{align}\label{eq: vial vs dosage balancing}
& \underbrace{\sum_{\nu \in \mathcal{V}} \iota_{\nu} N^v_{\nu k p t}}_\text{Opened vaccine dosages.} = \underbrace{\sum_{g \in \mathcal{G}} x^v_{g k p t}}_\text{Vaccinated people.} + \underbrace{w^{ov}_{k p t}}_\text{Open vial wastage.}, \hspace{0.25cm} \forall k,p,t
\end{align}

In their works, \citet{azadi2020optimization} and \citet{azadi2020developing} address the impact of vials on inventory management. They propose an integrated vial opening-inventory replenishment policy and develop a framework that takes into account the specific vaccine type and vial size. To mitigate the issue of open vial wastage, \citet{azadi2020developing} focus on childhood vaccines with different vial sizes within a two-stage stochastic model that considers demand uncertainty. The first stage of their model involves determining the timing, quantity, and size of vials to be ordered. The model aims to achieve a balance between the number of vaccines in terms of vials at the first stage and the number of doses at the second stage. Equations \eqref{eq:OPen Vial Wastage in forward time} to \eqref{eq:OPen Vial Wastage in backward time2} illustrate the process of dose balancing through the opening time ($t$) and vaccination time ($t'$). It is important to note that vaccinations can be safely administered before $t + \tau$, but not at $t + \tau$.

\noindent \begin{align}
& \underbrace{\sum_{t' = t}^{t' + \tau - 1} \sum_{g \in \mathcal{G}} x^v_{g k p t t'}}_{\substack{\text{Vaccinated people in terms of doses,} \\ \text{vaccination time in $t \leq t' \leq t + \tau - 1$.}}} + \underbrace{x^v_{g k p t t+\tau}}_\text{Expired vaccine.} = \underbrace{\sum_{\nu \in \mathcal{V}}  \iota_{\nu p} N^v_{\nu k p t}}_{\substack{\text{Amount of doses in} \\ \text{the opened vials.}}}, & \forall k, p, t \leq |\mathcal{T}| - \tau + 1 \label{eq:OPen Vial Wastage in forward time} \\
& \underbrace{\sum_{t' = t}^{|\mathcal{T}|} \sum_{g \in \mathcal{G}} x^v_{g k p t t'} \leq  \sum_{\nu \in \mathcal{V}}  \iota_{\nu p} N^v_{\nu k p t}}_{\substack{\text{Similar to Constraint \eqref{eq:OPen Vial Wastage in forward time}, by not considering} \\ \text{the vaccination after time horizon $|\mathcal{T}|$.}}}, & \forall k, p, t > |\mathcal{T}| - \tau + 1 \label{eq:OPen Vial Wastage in forward time2} \\
& \sum_{t \geq (t' - \tau + 1)^+}^{t'} \sum_{g \in \mathcal{G}} x^v_{g k p t t'}  + \underbrace{w^{ov}_{k p t'}}_\text{Open vial wastage.}  \leq \sum_{g \in \mathcal{G}} d_{g k p t'}, & \forall p, k, t'  \leq \tau \label{eq:OPen Vial Wastage in backward time} \\
& \underbrace{\sum_{t \geq (t' - \tau + 1)}^{t'} \sum_{g \in \mathcal{G}} x^v_{g k p t t'}}_{\substack{\text{Vaccinated people in terms of doses,} \\ \text{opening time in $t' - \tau + 1 \leq t \leq t'$.}}} + \underbrace{\sum_{g \in \mathcal{G}} s^v_{g p k t'}}_\text{Vaccine dose shortage.}  = \underbrace{\sum_{g \in \mathcal{G}} d_{g k p t'}}_\text{Vaccine demand.}, & \forall p, k, t'  > \tau \label{eq:OPen Vial Wastage in backward time2}
\end{align}

Finally, \citet{mak2022managing} assess the inventory dynamics of the rollout process for three strategies: (i) holding back second doses, (ii) releasing second doses, and (iii) stretching the lead time between doses.

\subsubsection{Priority Group Focus} \label{VCsubsubsection:priority-group}

In order to effectively address the pandemic, vaccine prioritization among different subgroups is crucial, particularly when faced with limited supply. Various factors can be taken into account for prioritization, such as age, profession, social contact, health condition, income, and race. Optimization modeling can be employed to incorporate vaccination strategies.
For instance, \citet{fadaki2022multi} allocate COVID-19 vaccines based on a \textcolor{blue}{``community transmission"} parameter, which reflects the number of unvaccinated individuals at a specific location at a given time. On the other hand, \citet{shiri2022equitable} introduce a formulation that divides the population into a high-priority group ($g = 1$) and a low-priority group ($g = 2$). Equations \eqref{eq:priority group focus} to \eqref{eq:priority group focus4} are used to define this grouping. The formulation includes lateral vaccine shipment between vaccination centers (VCs) exclusively during the first time period. With these constraints in place, the low-priority group receives vaccinations only after the high-priority group's needs have been met. Any remaining vaccines at the end of the time period are carried over to the subsequent period through inventory management.

\noindent \begin{align}
& y^v_{1 k p 1} = \min \{d_{1 k p 1}, I^{0}_k + \sum_{j \in \mathcal{N}} x^{d v}_{j k p 1}  + \sum_{\substack{k' \in \mathcal{V} \\ k' \neq k}} x^{v' v}_{k' k p 1} - \sum_{\substack{k' \in \mathcal{V} \\ k' \neq k}} x^{v v'}_{k k' p 1} - w_{k p 1} \}, & \forall k, p \label{eq:priority group focus} \\
& y^v_{2 k p 1} = \min \{d_{2 k p 1}, I^{0}_k + \sum_{j \in \mathcal{N}} x^{d v}_{j k p 1}  + \sum_{\substack{k' \in \mathcal{V} \\ k' \neq k}} x^{v' v}_{k' k p 1} - \sum_{\substack{k' \in \mathcal{V} \\ k' \neq k}} x^{v v'}_{k k' p 1} - y^v_{1 k p 1} - w_{k p 1} \}, & \forall k, p \label{eq:priority group focus2} \\
& y^v_{1 k p t} = \min \{d_{1 k p t}, I^{t-1}_k + \sum_{j \in \mathcal{N}} x^{d v}_{j k p t}   - w_{k p t} \}, & \forall k, p, t \geq 2 \label{eq:priority group focus3} \\
& y^v_{2 k p t} = \min \{d_{2 k p t}, I^{t-1}_k + \sum_{j \in \mathcal{N}} x^{d v}_{j k p t}  - y^v_{1 k p t} - w_{k p t} \}, & \forall k, p, t \geq 2 \label{eq:priority group focus4}
\end{align}

\section{Location and Routing Models}\label{sec: pure facility location model}

Decisions on distribution center (DC)  \citep{dastgoshade2022social,emu2021validating,soria2021proposal,srivastava2021strengthening} and vaccination center (VC) locations \citep{bertsimas2022locate,bravo2022optimal,cabezas2021two,leithauser2021quantifying,lim2016coverage,lusiantoro2022locational,polo2015location,kumar2022optimal,luo2023service,zhang2022mass} are important for vaccine supply-demand matching in the context of COVID-19. Models are based on demand, distance, coverage, priority, required service levels, and fairness. A recent article considers locating waste disposal center \citep{cao2023digital}. In this section, we present studies that focus on a variant of the assignment type (Section \S \ref{subsection:FLP}) or coverage type (Section \S \ref{subsection:coverage}) facility location models. In the formulations used in this section, vaccination center (VC) is the facility and population site (PS) is the demand source. Similar formulations can be used for the case where distribution center (DC) is the facility and vaccination center (VC) is the demand source. We present vehicle routing in Appendix \S \ref{sec: vehicle routing}.

\subsection{Facility Location Models} \label{subsection:FLP}

\citet{lusiantoro2022locational} consider a multi-objective facility location model (see \eqref{eq: p-median and max coverage mixture}), where the first objective function ensures effectiveness by maximizing demand coverage and the second objective function ensures efficiency by minimizing demand weighted distance, simultaneously \citep{gralla2014assessing}.

\noindent \begin{align}\label{eq: p-median and max coverage mixture}
& \text{(Obj. 1) } \max \underbrace{\sum_{s \in \mathcal{S}} \sum_{k \in \mathcal{V}} d^{P}_{s} x_{s k}}_\text{Demand coverage.}; \quad \text{(Obj. 2) } \min \underbrace{\sum_{s \in \mathcal{S}} \sum_{k \in \mathcal{V}} d^{P}_{s} \eta_{s k} x_{s k}}_\text{Demand weighted distance.} \nonumber \\
& \text{s.t. }  \underbrace{\sum_{k \in \mathcal{V}} x_{s k} \leq 1, \hspace{0.25cm} \forall s}_\text{Assigning a PS to at most one VC.}; \quad \underbrace{\sum_{k \in \mathcal{V}} Z^V_k \leq |\Bar{\mathcal{V}}|}_\text{There can be at most $|\Bar{\mathcal{V}}|$ VCs.}; \quad \underbrace{\sum_{s \in \mathcal{S}} d^{P}_{s} x_{s k} \leq C^V_k, \hspace{0.25cm} \forall k}_\text{VC capacity constraint.}\\
& \underbrace{0 \leq x_{s k} \leq Z^V_k, \hspace{0.25cm} \forall s,k}_\text{Can not assign to a non-existing VC.}; \quad  \underbrace{x_{s k} \leq a_{s k}, \hspace{0.25cm} \forall s,k}_\text{Viability constraint.}; \quad Z^V_{k} \in \{0, 1\}, \hspace{0.25cm} \forall k \nonumber
\end{align}

\citet{bravo2022optimal} consider a formulation similar to the one given in \eqref{eq: p-median and max coverage mixture}, where the second objective function is not considered and the decision variable $x_{s k}$ is a binary assignment variable, which captures percent coverage in the current formulation. In addition to the classical facility location factors, \citet{polo2015location} consider healthcare workforce planning and spatial accessibility using GIS data \citep{radke2000spatial}. \citet{cabezas2021two} use stochastic optimization to incorporate patient preference of vaccination centers (VCs) in terms of distance. \citet{bertsimas2022locate} is the only paper that considers an epidemiological model within a location model. \citet{dastgoshade2022social} develop an allocation optimization framework using Rawls' theory, Sadr's theory, and utilitarianism fairness theories. Unlike other studies, \citet{soria2021proposal} locate intermediate distribution centers (DCs) based on production plants and existing distribution centers (DCs). \citet{kumar2022optimal} compute number of optimal days a vaccination center (VC) should operate. \citet{luo2023service} integrate linear utility function with COVID-19 service center location decision that depends on travel distance, waiting time and location based service features. \citet{zhang2022mass} merge appointment scheduling with vaccination center (VC) location problem. The location decisions are dependent on the accepted appointments and vaccination sequencing in each site. \citet{cao2023digital} find the best locations for COVID-19 waste disposal centers by considering stochastic infection risk, while carrying the waste and the population around candidate disposal centers.  

\subsection{Coverage Models} \label{subsection:coverage}

\citet{lim2016coverage}, \citet{polo2015location} and \citet{srivastava2021strengthening} consider the maximum coverage problem (see \eqref{eq: p-median2}).

\noindent \begin{align}\label{eq: p-median2}
& \max \underbrace{\sum_{s \in \mathcal{S}} \bar{P}_s \mathfrak{z}_s}_\text{Total coverage.} \nonumber\\
& \text{s.t. } \underbrace{\mathfrak{z}_s  \leq \sum_{k \in \mathcal{V}|\eta_{sk} \leq \eta_\text{max}} Z^V_k, \hspace{0.25cm} \forall s}_{\substack{\text{A PS cannot be covered by} \\ \text{a VC not within $\eta_\text{max}$.}}}; \quad \underbrace{\sum_{k \in \mathcal{V}} Z^V_k \leq |\Bar{\mathcal{V}}|}_\text{There can be at most $|\Bar{\mathcal{V}}|$ VCs.}; \quad \underbrace{\sum_{k \in \mathcal{V}} {c}_k Z^V_k \leq \mathsf{b}}_\text{Budget constraint.} \\ 
& \qquad Z^V_k, \mathfrak{z}_s \in \{0,1\}, \hspace{0.25cm} \forall k,s \nonumber
\end{align}

\noindent There are multiple ways to calculate the coverage coefficient ($\Bar{P}_s$). The most common is to equate it to the population size in a population site (PS), assuming that in a covered site all population is covered. In practice, coverage can be affected by factors like accessibility and willingness to be vaccinated. \citet{xu2021hub} develop two different $\Bar{P}_s$ measures based on the number of public transportation spots and the social vulnerability index. \citet{lim2016coverage} provide a step-wise distance based coverage formulation, where coverage percentage and distance are inversely proportional (see \eqref{eq: step-wise coverage}).

\noindent \begin{align}
& \max \underbrace{\sum_{s \in \mathcal{S}} \bar{P}_s \sum_{q \in \mathcal{Q}_s}\theta_q \mathfrak{z}_{s q}}_\text{Step-wise coverage.} \nonumber \\
& \text{s.t. } \underbrace{\mathfrak{z}_{s q} \leq \sum_{k \in \mathcal{V}_s} Z^V_k, \mathcal{V}_s = \{k \in \mathcal{V}| \mathfrak{y}_{q-1} \leq \eta_{s k} \leq \mathfrak{y}_q \}, \hspace{0.25cm} \forall q, s}_\text{A PS cannot be covered by a VC outside allowed level.}; \quad \underbrace{\sum_{q \in \mathcal{Q}_s} \mathfrak{z}_{s q} \leq 1, \hspace{0.25cm} \forall s}_\text{A PS can be covered in at most one level $q$.} \label{eq: step-wise coverage} \\
& \qquad \sum_{k \in \mathcal{V}} Z^V_k \leq |\Bar{\mathcal{V}}|; \quad \sum_{k \in \mathcal{V}} {c}_k Z^V_k \leq \mathsf{b};\quad  Z^V_k, \mathfrak{z}_{s q} \in \{0, 1\}, \hspace{0.25cm} \forall k, q, s \nonumber
\end{align}

\noindent \citet{lim2016coverage} consider outreach centers (OC) for covering vaccination center (VC). Outreach happens when healthcare workers take vaccines to remote, temporary locations according to accessibility needs. Generally, these outreach centers (OCs) cover multiple rural vaccination centers (VCs) (see \eqref{eq: step-wise detail coverage1}).

\noindent \begin{align} \label{eq: step-wise detail coverage1}
& \underbrace{\mathfrak{z}_{s q} \leq \sum_{k \in \mathcal{V}_s} \sum_{\ell \in \mathcal{O}} \mathfrak{x}_{\ell k}, \mathcal{V}_s = \{k \in \mathcal{V}| \mathfrak{y}_{q-1} \leq \eta_{s k} \leq \mathfrak{y}_q \}, \hspace{0.25cm} \forall q,s}_\text{OC coverage constraint.}; \quad \sum_{q \in \mathcal{Q}_s} \mathfrak{z}_{s q} \leq 1, \hspace{0.25cm} \forall s  \\
& \underbrace{\sum_{k \in \mathcal{V}} \mathfrak{x}_{\ell k} \leq C^O_{\ell}, \hspace{0.25cm} \forall \ell}_\text{OC capacity constraint.}; \quad \underbrace{\eta_{\ell k}\mathfrak{x}_{\ell k} \leq \eta_\text{max}, \hspace{0.25cm} \forall \ell, k}_\text{OC distance constraint.}; \quad \underbrace{\sum_{\ell \in \mathcal{O}} \mathfrak{x}_{\ell k} \leq 1, \hspace{0.25cm} \forall k}_\text{A VC can be covered by one OC.} \nonumber
\end{align}

\noindent These types of location models support regional decision making by considering the specific requirements that may arise due to factors such as population density, existing pharmacy stores, and other infrastructure related factors. Further information on \textit{location-allocation} articles is given in Table~\ref{table:location-allocation-articles} in Appendix \S \ref{sec: vehicle routing}. 

\section{Uncertainty Modeling in COVID-19 Vaccine Supply Chain Models} \label{sec: uncertainty-scm}

In COVID-19 vaccine supply chain literature, demand uncertainty is the predominant form of uncertainty \citep{bennouna2022covid,bertsimas2020optimizing,enayati2020optimal,jadidi2021two,jarumaneeroj2022epidemiology,mehrotra2020model,minoza2021covid,rao2021optimal,thul2023stochastic,yang2022comparison,bertsimas2022locate,dastgoshade2022social,azadi2020optimization,mak2022managing,azadi2020developing,sinha2021strategies,basciftci2023resource,chowdhury2022modeling,goodarzian2022sustainable,habibi2023designing,kohneh2023optimization,lai2021multi,manupati2021multi,sazvar2021capacity,shiri2022equitable,shamsi2021novel,yang2021outreach}. However, research has also considered uncertainty in supply \citep{dastgoshade2022social,jahani2022covid,karakaya2023developing,mohammadi2022bi,chowdhury2022modeling,wang2023robust}, lead time \citep{jahani2022covid,karakaya2023developing,habibi2023designing,manupati2021multi,yang2021outreach}, infection \citep{roy2021optimal,cao2023digital,basciftci2023resource}, vaccination location choice \citep{cabezas2021two,luo2023service}, vaccine accessibility \citep{yarmand2014optimal,gilani2022data} and vaccine wastage \citep{shah2022investigation,rahman2023optimising}. While demand from a distibution center (DC) to manufacturer \citep{gilani2022data,goodarzian2022sustainable,manupati2021multi,shah2022investigation} and manufacturing lead time \citep{manupati2021multi} impact the DC-manufacturer level uncertainty, vaccination center (VC) level demand uncertainty and their variation across regions is influenced by numerous factors including gender, education, employment, exposure to misinformation, and hesitancy \citep{seboka2021factors}. These uncertainties affect both DC-VC and VC-PS related echelons. 

In terms of mathematical modeling of such uncertainties, stochastic inventory model \citep{manupati2021multi, shah2022investigation, sinha2021strategies, mak2022managing, yarmand2014optimal, chick2008supply}, two-stage stochastic programming \citep{lai2021multi, basciftci2023resource, mohammadi2022bi, karakaya2023developing, shiri2022equitable, azadi2020developing}, robust and distributionally robust optimization approaches \citep{basciftci2023resource, luo2023service, mohammadi2022bi, gilani2022data, wang2023robust}, chance-constraint modeling \citep{goodarzian2022sustainable, rahman2023optimising, azadi2020optimization} are noteworthy. Forecasting, predictive analytics, machine learning methods, simulation tools as well as their ensembles \citep{bennouna2022covid} are also used to obtain estimates of the uncertain parameters and aid optimization under various types of uncertainties. Other practices of modeling uncertainty found in COVID-19 vaccine literature include queuing \citep{jahani2022covid}, fuzzy programming \citep{kohneh2023optimization}, game theory \citep{sinha2021strategies}, decision trees \citep{manupati2021multi}, epidemiological modeling \citep{bertsimas2020optimizing,enayati2020optimal,jadidi2021two,jarumaneeroj2022epidemiology,minoza2021covid,rao2021optimal,thul2023stochastic,yang2022comparison}, etc. The latter is separately explored in Section \S\ref{sec:SEIRD models}. In this section, we present stochastic inventory models (Section \S \ref{ss:SIM}), two stage stochastic programming approaches (Section \S \ref{sec: TSSP}), robust/distributional robust approaches (Section \S \ref{ss:DRO}), and chance constrained modeling (Section \S \ref{ss:CCM}).

\subsection{Stochastic Inventory Models} \label{ss:SIM}

In response to supply chain uncertainties, vaccine literature, besides leveraging the inventory balancing equations discussed earlier, also uses to continuous or periodic inventory review policies. \citet{mak2022managing} focus on the scenario of very limited supply during the early period of COVID-19, examine varied two-dose vaccine rollout policies by analyzing the related inventory dynamics.

Under a continuous review policy, the reorder point and reorder quantity are the key decisions to make such that the remaining inventory ($d^{@l}_{\max}$) satisfies the lead time demand ($d(l)$) at the desired service level ($1 - \alpha$) (i.e., $\mathbb{P} (d(l) \leq d^{@l}_{\max} ) \geq 1 - \alpha$). \citet{manupati2021multi} use a continuous review policy and EOQ structure to model this uncertainty assuming Normally distributed demand.
In their work, resulting total cost for a distribution center (DC) and reorder point is expressed in terms of the order quantity as shown in \eqref{eq:uncertain-DC}. \citet{shah2022investigation} also use similar type of ensembles of continuous review and EOQ policy and investigate carbon emission from a vaccine inventory system.

\noindent \begin{align}
& \underbrace{(\sum_{i}c^{OC}_{i j})(\frac{d^D_j}{Q_j}) + c^{HC}_j(\frac{Q_j}{2}+ \mathsf{Z}_{1 - \alpha} \sqrt{(l_j)(\sigma^2_{{d}^D_j})})}_\text{Total cost for DC $j$.}; \quad \underbrace{({d}^D_j) (l_j) + \mathsf{Z}_{1 - \alpha} \sqrt{({d}^D_j)^2 (\sigma^2_{l_j}) + (l_j)^2 (\sigma^2_{{d}^D_j})}}_\text{Reorder point of DC $j$.} \label{eq:uncertain-DC}
\end{align}

Periodic review policy-based inventory decisions that review stocks periodically are considered using a newsvendor modeling framework. \citet{yarmand2014optimal} propose a closed-form solution to Phase-1 allocation of a two-phase vaccine allocation problem using the newsvendor approach. This work uses the ratio of underage cost to summation of underage and overage cost and finds a closed-form solution for vaccine allocation.  The newsvendor structure is also utilized to model contracts for supply chain coordination \citep{chick2008supply}. 

\citet{sinha2021strategies} use inventory balancing equations similar to those described in the earlier sections to model uncertainties. For instance, using forecasts, \citet{sinha2021strategies} replace the $\sum_{p} d_{k p t}$ term by $\gamma \sum_{p} d_{k p t}$ in \eqref{eq: Inventory_balance_VC} for \textit{critical nodes} of a supply chain network. Due to budget and infrastructure inefficiencies (especially in developing countries), vaccination centers (VCs) may not get sufficient number of vaccines from distribution centers (DCs) on a regular basis. A model in such a context adds binary variables $y^{dv}_{j k t}$ to decide at which time period a distribution center (DC) should supply vaccine to a vaccination center (VC). For a similar reason, lead time can be uncertain as well. \citet{sinha2021strategies} analyze major disruption scenarios and store additional vaccine supplies in strategic higher-echelon nodes to prevent transshipment lead time uncertainty.

\subsection{Two-stage Stochastic Programming Approach}\label{sec: TSSP}

In a two-stage stochastic programming (TSSP) approach, the first stage makes here-and-now strategic decisions while, given the first stage decisions, the second stage makes wait-and-see decisions. Below, we outline a COVID-19 motivated TSS mixed-integer linear programming model from a number of papers \citep{mehrotra2020model, lai2021multi, basciftci2023resource, mohammadi2022bi, karakaya2023developing}. 
\begin{align}\label{eq: TSSP}
 &\min \underbrace{\mathtt{c}^{\top} \mathtt{x}}_{\text{First-stage cost.}} + \underbrace{\mathbb{E}_{\xi \sim \mathbb{P}} [\mathfrak{P}(\mathtt{x}, \xi)]}_{\text{Expected optimal cost of second-stage problem.}} \,\, \text{s.t.} \,\, \underbrace{\mathtt{A x} \geq \mathtt{a}, \,\,  \text{$\mathtt{x}$ is Integer}}_{\text{First-stage constraints.}} \\
&\underbrace{\mathfrak{P}(\mathtt{x}, \xi^{\omega})}_{\substack{\text{Recourse problem parameterized by $\mathtt{x}$ and} \\ \text{random parameters $\xi^{\omega}$ under scenario $\omega \in \Omega$.}}} := \hspace{1.5cm} \min \, \underbrace{\mathtt{e}^{\omega^{\top}} \mathtt{y}^{\omega}}_{\text{Second-stage cost for scenario $\omega$.}} \\
 &\hspace{5.75cm} \text{s.t.} \, \underbrace{\mathtt{B}^{\omega} \mathtt{y}^{\omega} \geq \mathtt{b}^{\omega} - \mathtt{T}^{\omega} \mathtt{x}^{\omega} }_{\substack{\text{Second-stage coupling constraint} \\ \text{with technology matrix $\mathtt{T}^{\omega}$.}}}; \underbrace{\mathtt{D}^{\omega} \mathtt{y}^{\omega} \geq \mathtt{d}^{\omega}}_{\substack{\text{Second stage constraint} \\ \text{free of first-stage variable $\mathtt{x}$.}}} \nonumber\\
 & \hspace{6.8cm} \mathtt{y}^{\omega} \, \text{is continuous or mixed-integer} \nonumber
\end{align}

\citet{mehrotra2020model} were the earliest to use a TSSP framework in the context of ventilator allocation. The TSSP model proposed by \citet{lai2021multi} minimizes setup costs by considering the opening of vaccination centers (VCs), vaccination center (VC) assignment to population sites, and required healthcare personnel as first-stage variables. Each discrete scenario of the second stage model incorporates multi-period planning and considers vaccine transportation, inventory, shortages, and satisfied demand at each location within the planning horizon as decision variables. The number of personnel required at each vaccination center (VC) in each time period is also considered, ensuring that the minimum average realized demand threshold is met based on their service rate. Additionally, a coupling constraint prevents the number of healthcare personnel needed to meet demand from exceeding the total workers assigned in the strategic stage at any given time period. Similarly, demand of a population site is satisfied by a vaccination center (VC) if only if the first stage establishes a one-to-one assignment between them. 

In the work of \citet{basciftci2023resource}, the TSSP model minimizes setup and operational costs of different distribution centers (DCs), while considering inventory balancing and lead time. However, it does not account for the connection between facility requirements and healthcare workers. \citet{karakaya2023developing}, on the other hand, propose a TSSP model that focuses on supply level uncertainty and aims to minimize the expected weighted sum of deviations from the vaccination timeline for each priority group made in the first stage. In each scenario, their second stage model calculates earliness and lateness deviation, taking into account factors such as scenario-specific vaccination starting and completion time, dosage scheme, inventory availability, and other factors with necessary coupling of the first stage decisions. \citet{mohammadi2022bi} present a TSSP model that incorporates uncertainty in both supply and demand, including parameters related to pandemic progression, vaccination effectiveness, social contacts, and uncertain death rates. Their objective is to minimize cost and the number of deaths, with the added consideration of vaccination center (VC) capacity determined through a queuing system.

\subsection{Robustness and Distributionally Robust Approaches} \label{ss:DRO}

In addition to their TSSP approach, \citet{mohammadi2022bi} test a robust approach that minimizes deaths using a scenario-dependent second objective value ($\vartheta^{\omega}_2$) of a bi-objective model. The resulting objective given in \eqref{eq: robustify by std} is subject to the same constraints as their TSSP model. To solve the problem, the non-linearities due to the death minimization objective and robustification are linearized.
\begin{align}\label{eq: robustify by std}
\min \underbrace{\sum_{\omega \in \Omega} \mathsf{p}^{\omega} \vartheta^{\omega}_2}_{\text{Expected number of deaths.}} + \underbrace{\check{\mathbf{M}} \sum_{\omega \in \Omega} \mathsf{p}^{\omega} |\vartheta^{\omega}_2 - \sum_{\omega' \in \Omega} \mathsf{p}^{\omega'} \vartheta^{\omega'}_2|}_{\text{Weighted mean deviation from mean objective value.}}
\end{align}

Motivated by \citet{bertsimas2016analytics}, \citet{gilani2022data} propose a dynamic robust optimization approach based on uncertainty sets and cutting hyperplanes, to reduce conservatism compared to its static counterpart while maintaining robustness. \citet{wang2023robust} develop a two-stage robust optimization model for facility location and scheduling under supply uncertainty, with the goal of  distributing vaccines to remote areas using drones.

Besides a TSSP model, \citet{basciftci2023resource} also utilize a distributinally robust optimization (DRO) technique for robustification \citep{rahimian2019distributionally}. The DRO approach selects the worst-case distribution from an ambiguity set ($\mathcal{B}$) defined on some support set for random parameters, without assuming any specific probability distribution for the expected cost calculation in the TSSP model. \citet{basciftci2023resource} assume finite support set and moment-based ambiguity set as defined in \eqref{eq:DROvsSP}. \citet{luo2023service} propose a utility-robust service center location model with a pandemic case study. The ambiguity set used for uncertain linear utility function parameter are a generalized version of a moment-based ambiguity set \eqref{eq:DROvsSP} because of its additional dependence on the location decision.

\noindent \begin{align}
& \underbrace{\mathcal{B} = \big\{\mathsf{p} \in \mathbb{R}^{|\Omega|}_+ | \mu_{\Tilde{d}_{k t}} \hspace{-1mm} - \epsilon^{\mu}_{k t} \leq \sum_{\omega \in \Omega} \mathsf{p}_{\omega} \Tilde{d}_{\omega k t} \leq  \mu_{\Tilde{d}_{k t}} \hspace{-1mm} + \epsilon^{\mu}_{k t}, \sigma_{\Tilde{d}_{k t}} \underline{\epsilon}^{\sigma}_{k t} \leq \sum_{\omega \in \Omega} \mathsf{p}_{\omega} \Tilde{d}^2_{\omega k t} \leq  \sigma_{\Tilde{d}_{k t}} \overline{\epsilon}^{\sigma}_{k t}, \hspace{0.25cm} \forall k, t, \mathsf{p}^{\top} \mathbf{1}= 1\big\}}_{\text{Finitely supported moment based ambiguity set under imperfect moment information.}} \label{eq:DROvsSP}
\end{align}

\subsection{Chance-constrained Modeling} \label{ss:CCM}

In stochastic optimization, handling chance-constrained (CC) problems can be challenging. However, in the context of the vaccine literature, the normality assumption of random parameters simplifies the problem and is commonly used. For example, \citet{azadi2020optimization} use chance constraints to ensure that supplied amounts of pediatric vaccine at every vaccination center (VC) in every time period ($x^{dv}_{j k p t}$), meet the random demand ($\Tilde{d}_{k p t}$) at least $1 - \alpha$ percent of time (i.e., $\mathbb{P} \bigg( \sum_{j} x^{dv}_{j k p t} \geq \Tilde{d}_{k p t} \bigg) \geq 1 - \alpha, \hspace{0.25cm} \forall k, p, t$).  Similarly, \citet{rahman2023optimising} impose chance constraints for total cost and make span of the vaccination program to remain within their random budget with at least a probability of $1-\alpha$. If normality assumption remains in place, such presence of a random parameter only in the right hand side of the expression (not as a coefficient of variables) is very easy to deal with since resultant deterministic equivalent of CC is linear. Other than \citet{azadi2020optimization}, \citet{goodarzian2022designing} also leverage this technique to impose chance constraints for every vehicle responsible to distribute vaccine between each of the possible supply chain players: manufacturer (main and local plant) to distribution center (DC), or warehouse, warehouse to vaccination center (VC) (hospital, pharmacy), according to the respective uncertain demands. If $\Tilde{d}_{ t h}$ is the random demand at time $t$ to be transported by vehicle $h$, a linear reformulation of ${\mathbb{P} \big(\sum_{j, k, p} x^{dv}_{j k p t h} \geq \Tilde{d}_{t h}\big) \; \forall t,  h}$ is
\begin{align}\label{eq: CC normal}
\underbrace{\sum_{j, k, p} x^{dv}_{j k p t h} \geq \mu_{\Tilde{d}_{t h}} + \mathsf{Z}_{1 -\alpha} \sigma_{\Tilde{d}_{t h}}, \hspace{0.25cm} \forall t, h}_{\text{Each vehicle $h$ should at least transport $1-\alpha$ confidence level (upper) of the demand at time $t$.}} 
\end{align}

\section{Equity, Sustainability and Multi-objective Modeling} \label{sec:eq-and-sus}

\subsection{Equity Considerations} 

In public health decision-making, equity plays a significant role. A significant number of papers has studied this concern in the context of COVID-19 \citep{balcik2022mathematical,bertsimas2020optimizing,enayati2020optimal,munguia2021fair,dastgoshade2022social,lim2016coverage,fadaki2022multi,rastegar2021inventory,shiri2022equitable,tavana2021mathematical,wang2023robust,bennouna2022covid,basciftci2023resource,mohammadi2022bi}. Equity can be integrated into decision models by incorporating it into the objective function and/or the constraints. Within the COVID-19 literature, two main dimensions of equity are examined: priority group-based \citep{fadaki2022multi,enayati2020optimal,shiri2022equitable,basciftci2023resource,mohammadi2022bi}, geography-based \citep{lim2016coverage,wang2023robust,bertsimas2020optimizing,munguia2021fair} or both \citep{rastegar2021inventory,tavana2021mathematical,balcik2022mathematical,dastgoshade2022social,bennouna2022covid}. Priority group equity involves distributing vaccines equitably among diverse sub-populations, taking into account factors like risk, age, occupation, and other demographic characteristics. This approach ensures that individuals with higher vulnerability or essential roles have equitable access to vaccines. Geography-based equity aims to ensure that vaccines are distributed equitably across different geographical locations. Decision makers consider factors such as population density, healthcare infrastructure, and disease prevalence to achieve an equitable allocation of vaccines.

The commonly used equity measures are: equitable demand satisfaction \citep{tavana2021mathematical, rastegar2021inventory,shiri2022equitable}, deviation from fair coverage level \citep{balcik2022mathematical}, prohibiting daily variations in the number of vaccinations distributed to each location \citep{bertsimas2020optimizing} and using equity frameworks (e.g., Social Welfare, Nash, Rawlsian Justice) \citep{munguia2021fair}. For example, \citet{tavana2021mathematical}, \citet{rastegar2021inventory} and \citet{shiri2022equitable} use \eqref{eq: equity based objective} and \eqref{eq: equity based objective 2} to model equity in the form of demand satisfaction:

\noindent \begin{align}
& \underbrace{\max \min_{k, g} \bigg\{ \frac{\sum_{p, t} x^v_{g k p t}}{\sum_{p, t} d_{g k p t}}\bigg\}}_{\substack{\text{Maximizing the minimum} \\ \text{vaccination percentage.}}}; \quad \underbrace{\sum_{p, t} x^v_{g k p t} \geq \gamma \sum_{p, t} d_{g k p t} \; \forall g,  k}_{\substack{\text{Minimum demand satisfaction} \\ \text{rate constraint.}}} \label{eq: equity based objective}\\ 
& \quad \underbrace{\frac{\sum_{i, t} x^{md}_{i j p t}}{\sum_{g, k , t} d_{g k p t} x_{j k }} \geq \gamma, \hspace{0.25cm} \forall j, p}_{\substack{\text{Minimum demand satisfaction rate} \\ \text{constraint with assignment.}}}  \label{eq: equity based objective 2}
\end{align}

\citet{balcik2022mathematical} incorporate equity across different regions and priority groups using \eqref{eq: equity based objective3.1} and \eqref{eq: equity based objective3.2} based on the deviation of actual allocation from equitable quantity estimated. Negative deviation with respect to demand (first term of \eqref{eq: equity based objective3.1}) is  discouraged  by using a penalty term $\check{\mathbf{M}}$.
\noindent \begin{align}
\min \,\underbrace{\check{\mathbf{M}} \sum_{r = 1}^{|\mathcal{R}|} \sum_{g = 1}^{|\mathcal{G}|}  \frac{\varsigma \mathtt{\Delta}^{-}_{r g}}{d_{r g}} + \sum_{r = 1}^{|\mathcal{R}|} \sum_{g = 1}^{|\mathcal{G}|}  \frac{(1 - \varsigma) \mathtt{\Delta}^{+}_{r g}}{d_{r g}}}_{\substack{\text{Weighted sum of deviation from} \\ \text{fair amount to demand ratio.}}} + \underbrace{\sum_{r = 1}^{|\mathcal{V}|} \sum_{g = 1}^{|\mathcal{G}|}  \mathbf{M}\cdot  \mathtt{\Pi}_{r g}}_{\substack{\text{Penalizing gap between minimum} \\ \text{percentage demand and allocation.}}} \label{eq: equity based objective3.1}\\
\underbrace{\sum_{p = 1}^{|\mathcal{P}|} x^v_{g r p} = e^v_{r g} + \mathtt{\Delta}^+_{r g} - \mathtt{\Delta}^-_{r g} }_{\substack{\text{Deviation of allocation} \\ \text{from fair amount.}}}; \quad
\underbrace{e^v_{r g} + \mathtt{\Delta}^+_{r g} - \mathtt{\Delta}^-_{r g} \leq d_{r g}}_{\substack{\text{Allocation should not be} \\ \text{more than the demand.}}}, \hspace{0.25cm} \forall r, g \label{eq: equity based objective3.2}
\end{align}

The approaches discussed in \citet{munguia2021fair} include Rawlsian justice and the Social Welfare II approaches. Rawlsian justice principle maximizes the minimum allocated vaccine over all regions. Social Welfare II ensures vaccine allocation per region by maximizing the sum of their squared deviations from the allocation quantity possible in each region without equity consideration to achieve the overall maximum allocation.
\citet{bertsimas2022locate} consider equity in deciding the location of vaccination centers (VCs) at a strategic level. The number of sites opened in each state is determined based on the total target facility setup weighted by the state's population share, with efficiency-fairness trade-offs controlled through model parameters. \citet{enayati2020optimal} incorporate equity using the well-known Gini coefficient (or index) expressed as $GI(f) = \frac{\sum_{g' \in \mathcal{G}} \sum_{g \in \mathcal{G}} |f_g - f_{g'}|}{2 |\mathcal{G}| \sum_{g' \in \mathcal{G}} f_{g'}}$, where $f_{g'}$ is a decision variable representing the fraction of population covered from population class $g'$.

\subsection{Sustainability Considerations}

Sustainability focused models consider the environmental cost and social benefit, with the aim of minimizing carbon emissions and adhering to carbon emission budgets. The quantification of environmental cost typically involves measuring carbon emissions from activities such as truck usage, inventory holding, and waste generation or disposal \citep{gilani2022data,goodarzian2022designing,chowdhury2022modeling,goodarzian2022sustainable,sazvar2021capacity}. On the other hand, the social cost is assessed by considering factors such as waiting times for vaccination, transmission risk at workplaces, minimizing unmet demand, delivery time reduction, and staff balancing. Meanwhile, social benefits are in the form of job creation \citep{gilani2022data,goodarzian2022designing,shamsi2021novel,goodarzian2022sustainable,chowdhury2022modeling,sazvar2021capacity}. An illustrative example can be found in \eqref{eq:env obj} where the two equations are two objective functions, based on the work of \citet{goodarzian2022designing}.

\noindent \begin{align}
& \min \; \underbrace{c^{EI} \big(\sum_{j} (Y^{D(c)}_j + Y^{D(vc)}_j + Y^{D(uc)}_j) + \sum_{k} Z^V_k\big)}_\text{Carbon emission cost of establishing DCs and VCs.} + \underbrace{c^{CR} \hspace{-5mm} \sum_{i, j, j', k, k', p, t} \hspace{-5mm} \big(\eta_{i j} x^{m d}_{i j p t}  + \eta_{j k} x^{d v}_{j k p t}  + \eta_{j j'} x^{d d'}_{j j' p t}  + \eta_{k k'} x^{v v'}_{k k' p t} \big)}_\text{Carbon emission cost of vaccine transshipments.} \nonumber\\
& \min \underbrace{\sum_{g,k,p,t} s^v_{g k p t}}_\text{Total vaccine shortages.} \label{eq:env obj} 
\end{align}

\citet{bonney2011environmentally} use carbon emissions, waste disposal cost and social cost, in addition to ordering cost per replenishment cycle, to derive the optimal order quantity of an environmentally-friendly EOQ model. Moreover, \citet{shah2021prevalence} consider carbon tax and amount of carbon emissions during manufacturing, holding, preparation, purchase and inspection processes to build an inventory model that minimizes total cost of an inventory system per unit time.

\srevision{
\subsubsection{Waste Management}
There are two optimization papers that consider waste management in the context of COVID-19 vaccine supply chain. \citet{abbasi2023designing} emphasize the importance of building a green SCM that not only addresses the environmental impact of biomedical and plastic waste but also considers the hygiene costs associated with waste management. In parallel, \citet{icsik2023optimizing} delves into provincial-level waste management in Turkey, with a specific focus on medical waste generated during the vaccine distribution process.}

\subsection{Multi-objective Considerations}
Equity, sustainability, and economic objectives often conflict. Several papers  consider conflicting multiple objective functions in resource allocation \citep{anahideh2022fair,minoza2021covid}, location-allocation \citep{bertsimas2022locate,cao2023digital,dastgoshade2022social,leithauser2021quantifying,lusiantoro2022locational,polo2015location}, inventory management \citep{fadaki2022multi,jahani2022covid,mohammadi2022bi}, supply chain management \citep{chowdhury2022modeling,goodarzian2022designing,gilani2022data,goodarzian2022sustainable,kohneh2023optimization,li2021locate,sazvar2021capacity,tang2022bi,wang2023robust, mohammadi2022bi} and routing \citep{shamsi2021novel}. For example, expressions in \eqref{eq:env obj} are conflicting, since, in order to prevent vaccine shortage, more vaccine should be shipped and more facilities should be located. In the bi-objective TSSP and robust optimization model of \citet{mohammadi2022bi}, death and cost minimization are also conflicting to each other. To handle the conflicting nature of the models, different multi-objective approaches have been employed. Regardless of the number of objectives in these multi-objective models, most literature eventually works with a single objective after using other objective(s) as $\varepsilon-$constraints.

\section{Models Considering Epidemiological Dynamics}\label{sec:SEIRD models}

Demand estimation in vaccine management is accomplished using epidemiological models. Several articles (especially in resource allocation category)  unify mathematical optimization with compartmental epidemiological modeling \citep{bertsimas2020optimizing,enayati2020optimal,jadidi2021two,jarumaneeroj2022epidemiology,minoza2021covid,rao2021optimal,roy2021optimal,thul2023stochastic,yang2022comparison,bertsimas2022locate,mak2022managing,shamsi2021novel,bennouna2022covid}, for decision making under demand uncertainty. SIR (Susceptible-Infectious-Recovered) model is the most commonly used compartmental epidemiological model, which is represented by a set of ordinary differential equations (ODE). The SIR model has many extensions (i.e., SEIR, SEIQR, SIRD, DELPHI). \citet{bertsimas2020optimizing} and \citet{bertsimas2022locate} consider the DELPHI-V model, representing the rate of susceptible ($S$), exposed ($E$), infectious ($I$), undetected ($U$), quarantined ($Q$) and dead ($D$) population in a certain region, using compartment-wise flow rates. The compartmental flows for a region are illustrated in Figure~\ref{fig:SEIRD_bertsimas}. Table~\ref{Tab:epidemiological parameter description} outlines the nomenclature used in the DELPHI-V model and the differential equations given in \eqref{eq:susceptible_risk_class_l} through \eqref{eq:immunity_risk_class_l} provide the ODE representation of the system dynamics model.

\begin{figure}[htp]
\centering
\includegraphics[width=0.6\textwidth]{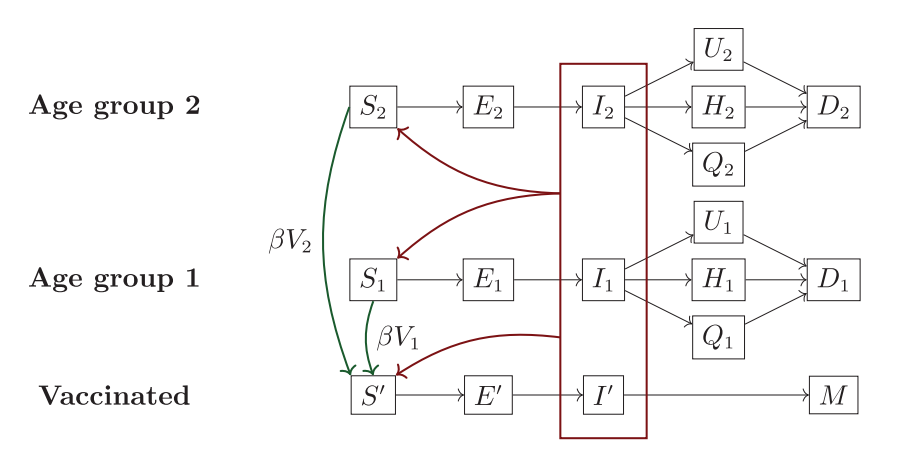}
\caption{DELPHI-V Compartmental Model (Adapted from \citep{bertsimas2022locate})}
\label{fig:SEIRD_bertsimas}
\end{figure}

\begin{table} \setlength{\tabcolsep}{4pt}
\vspace{-0.3cm}
\caption{Parameter Description of Compartmental Epidemiological Model}
\vspace{-0.5cm}
\begin{center}
\begin{adjustbox}{width=\textwidth}
\renewcommand{\arraystretch}{1}
\begin{tabular}{cl}
\hline
\rowcolor{Gray}
\textbf{Notation}&\textbf{Description}\\
\hline
$\beta$ & Vaccine effectiveness percentage. \\
$\alpha$ & Nominal infection rate. \\
$\gamma(t)$ & Governmental and societal response function modulating infection rate over different time. \\
 $r^I$ &  Infection progression rate. \\
$r^d$ &  Detection rate. \\
$r^D$ & Death rate. \\
$r^U_g(t)$ & Transition rate to undetected compartment for risk group $g$. \\
$r^H_g(t)$ & Hospitalization rate  for risk-based population risk group $g$. \\
$r^Q_g(t)$ & Quarantine rate for risk-based population risk group $g$. \\
\hdashline
$S_{rgt} (S'_{r t})$  & Number of susceptible people in region $r$, risk group $g$ (vaccinated class, see Fig. \ref{fig:SEIRD_bertsimas}) and time $t$. \\
$\bar{S}_{rgt}$ & Number of eligible people in region $r$, risk group $g$ and time $t$. \\
$V_{rgt}$ & Number of vaccinated people in region $r$, risk group $g$ and time $t$. \\
$E_{rgt} (E'_{r t})$ & Number of exposed people in region $r$, risk group $g$ (vaccinated class) and time $t$. \\
$I_{rgt} (I'_{r t})$ & Number of infectious people in region $r$, risk group $g$ (vaccinated class) and time $t$. \\
$U_{rgt}$ & Number of undetected people in region $r$, risk group $g$ and time $t$. \\
$D_{rgt}$ & Number of dead people in region $r$, risk group $g$ and time $t$. \\
$M_{rgt}$ & Number of immunized people in region $r$, risk group $g$ and time $t$. \\
$Q_{rgt}$ & Number of quarantined people in region $r$, risk group $g$ and time $t$. \\
$H_{rgt}$ & Number of hospitalized people in region $r$, risk group $g$ and time $t$. \\
\hline
\end{tabular}
\end{adjustbox}
\end{center}
\label{Tab:epidemiological parameter description}
\end{table}

\FloatBarrier

\noindent \begin{align}
& \underbrace{\frac{dS_g}{dt}}_{\substack{\text{Rate at which susceptible population under risk} \\ \text{class $g$ flows out to other compartments.}}} = - \beta V_g(t) - \alpha \gamma(t) ( S_g(t) - \beta V_g(t) ) \underbrace{\bigg ( \sum_{g \in \mathcal{G} } I_g(t) + I'(t) \bigg)}_{\substack{\text{Currently infected} \\ \text{population.}}} \label{eq:susceptible_risk_class_l}\\
& \frac{dS'}{dt} = + \beta \sum_{g \in \mathcal{G}} V_g(t) - \alpha \gamma (t) \bigg(S'(t) + \beta \sum_{g \in \mathcal{G}} V_g(t)\bigg) \bigg(\sum_{g \in \mathcal{G}} I_g(t) + I'(t)\bigg)\label{eq:susceptible_risk_class_vaccinated}\\
& \frac{dE_g}{dt} = \alpha \gamma(t) (S_g(t) - \beta V_g(t)) \bigg(\sum_{g \in \mathcal{G}} I_g(t) + I'(t) \bigg) - r^I E_g(t)\label{eq:exposed_risk_class_l}\\
& \frac{dE'}{dt} = \alpha \gamma(t) \bigg(S'(t) + \beta \sum_{g \in \mathcal{G}} V_g(t)\bigg) \bigg(\sum_{g \in \mathcal{G}} I_g(t) + I'(t) \bigg) - r^I E'(t)\label{eq:exposed_risk_class_vaccinated}\\
& \frac{dI_g}{dt} = r^I E_g(t) - r^d I_g(t); \;\; \frac{dI'}{dt} =  r^I E'(t) - r^d I'(t)\label{eq:infectious_risk_class_vaccinated}\\
& \frac{dU_g}{dt} = r^U_g(t) I_g(t) - r^D U_g(t); \;\ \frac{dH_g}{dt} = r^H_g(t) I_g(t) - r^D H_g(t); \; \frac{dQ_g}{dt} = r^Q_g(t) I_g(t) - r^D Q_g(t)\label{eq:quarantined_risk_class_l}\\
& \frac{dD_g}{dt} = r^D(U_g(t) + H_g(t) + Q_g(t))\label{eq:death_risk_class_l}\\
& \frac{dM}{dt} = r^d I'(t) \label{eq:immunity_risk_class_l}
\end{align}

\subsection{Embedding SIR Dynamics into Optimization Modeling}

In order to embed a compartmental model into an optimization model, ODEs are discretized within a time horizon $\mathcal{T}$ with discretization unit $\Delta t$ (generally articles segment $\mathcal{T}$ into $|\mathcal{T}|$ many days). Equations \eqref{seir-opt-obj} to \eqref{eq: SEIRD variables} represent a unified epidemiological-optimization framework for region and compartment-based, multi-period vaccine allocation. The key decision variable $x^{\kappa}$ is related to vaccine supply chain management while $y^{\kappa}$ is related to the number of people in each compartment at each time period for each region  (see \eqref{eq: targeted vaccination SEIRD} to~\eqref{eq: SEIRD variables}).

\noindent \begin{align}
&  \min \underbrace{c^{\kappa} \; x^{\kappa}}_\text{Cost.}  + \underbrace{d^{\kappa} \; y^{\kappa}}_{\substack{\text{Epidemiological spread and} \\ \text{deaths over time period.}}} \label{seir-opt-obj} \\
& \text{s.t.} \;\; \underbrace{A^{\kappa} x^{\kappa} \leq a^{\kappa}}_{\substack{\text{Manufacturer, DC and VC} \\ \text{spesific constrains.}}} \label{eq: supply constraint SEIRD}\\
& \qquad  \underbrace{B^{\kappa} y^{\kappa} \geq b^{\kappa}; \;\; f(y^{\kappa}) \geq 0}_{\substack{\text{Discretized version of} \\ \text{ODEs \eqref{eq:susceptible_risk_class_l} to \eqref{eq:immunity_risk_class_l}.} \label{eq: targeted vaccination SEIRD}}} \\
& \qquad \underbrace{C^{\kappa} x^{\kappa} \leq D^{\kappa} y^{\kappa} + e^{\kappa}}_{\substack{\text{Linking constraint sets \eqref{eq: supply constraint SEIRD} and \eqref{eq: targeted vaccination SEIRD} for meeting} \\ \text{vaccination target, schedule, etc.}}} \label{eq: SEIRD linking supply chain}\\
& \qquad x^{\kappa} \;\; \text{is mixed binary or integer}, \;\; y^{\kappa} \; \text{is continuous or pure integer}
\label{eq: SEIRD variables} \\
& \qquad \{\mathbf{\overline{S}}, \mathbf{S}, \mathbf{V}, \mathbf{E}, \mathbf{I}, \mathbf{U}, \mathbf{D}, \mathbf{M}, \mathbf{Q}\} \subseteq y^{\kappa} 
\nonumber
\end{align}

\noindent Equation \eqref{eq: targeted vaccination SEIRD} builds the population upper bounds in each compartment. For example, \eqref{eq:susceptible_risk_class_l} can be used in the form {$S_{rgt+1} \geq S_{rgt} - \beta V_{rgt} - \alpha_r \gamma_{rt} ( S_{rgt} - \beta V_{rgt} ) \bigg ( \sum_{g \in \mathcal{G} } I_{rgt} + I'_{rt} \bigg)\Delta t \text{ where, } \Delta t = 1$.} Embedding $|\mathcal{T}|$ many forward difference based representative equations for the system of ODEs given in \eqref{eq:susceptible_risk_class_l} to \eqref{eq:immunity_risk_class_l} into this optimization model for each region of interest gives an estimation of the proportion of population in each compartment in each time period, alongside of decisions on vaccine distribution. Note that, non-linearity arises from the multiplication of $S_{rgt}$, $I_{rgt}$ and $V_{rgt}$, $I_{rgt}$ pairs. Finally, in \eqref{eq: SEIRD linking supply chain}, not all the elements are required to be linked.

It is common to write constraints in terms of epidemiological compartments. For instance, \citet{bertsimas2022locate} use \eqref{eq: supply and vaccine eligibility SEIRD} to limit vaccination ($\mathbf{V}$) by eligible population ($\mathbf{\overline{S}}$), since number of susceptible individuals ($S_{rgt}$) is highly dependent on vaccine availability. Note that eligible population ($\mathbf{\overline{S}}$) is determined by susceptible population ($\mathbf{S}$) and vaccine effectiveness ($\beta$).

\noindent \begin{align}
& \overline{S}_{rgt+1} \leq \overline{S}_{rgt} - \beta V_{rgt} - (S_{rgt} - S_{rgt+1}) \; \forall r, g, t; \quad V_{rgt} \leq \overline{S}_{rgt} \; \forall r, g, t \label{eq: supply and vaccine eligibility SEIRD}
\end{align}

\citet{bertsimas2022locate} also decide where to locate mass vaccination centers (VCs) and how to assign each state to a vaccination center (VC) by embedding an epidemiological model into a facility location optimization model. The model aims to minimize total mortality and hospitalizations. Finally, their analytical approach determines the quantity of vaccines allocated to each vaccination center (VC) and the eligibility of individuals within a particular risk category in a state for receiving the vaccination. Several SIR related literature can be considered in the context of the general structure sketched above. The optimization model in \citet{jarumaneeroj2022epidemiology}, for instance, is similar but considers an infectious case minimization objective. \srevision{\citet{barthspatiotemporal} incorporated epidemiological dynamics with time-varying contact rate. The corresponding model is thus able to reflect multiple waves during pandemic. Note that other than an SIR type structure, a vaccine allocation model can reflect infection spread dynamics by incorporating subject-specific risk and contact information accounting social network, such as the model presented by \citet{liquantifying} to underscore the benefits of customized vaccine distribution strategies.}

\subsection{Graph Laplacian-based Dynamics Representation} 
Epidemiological dynamics are also represented by graph Laplacian-based dynamic models. \citet{enayati2020optimal} compute the minimum number of required vaccines to prevent the virus spread, while ensuring sufficient vaccine coverage. They model the infection growth in terms of generation instead of time \citep{hill2003critical}. Using the current population in each compartment, their model gives the number of new infections due to contact with exposed or infectious individuals. 

\subsection{Model Parameterization and Herd Immunity}
Embedding a discretized SIR model results in a large optimization model that can be computationally expensive to solve. Consequently, several papers incorporate SIR dynamics indirectly into their model. For example, \citet{rao2021optimal} use a SIR model-generated patient population over time in each compartment as parameters in their optimization model. Their models examine different objectives separately including minimization of infections, deaths, life years lost, and quality adjusted life years. The  model requires vaccination of a minimum fraction of population in each region. \citet{minoza2021covid} also use the output from a SIR model as parameters in an optimization model that focuses on risk-class based prioritization for front-line workers and the elderly. As an alternative approach to incorporate an SIR model into optimization framework, \citet{shamsi2021novel} leverage optimal control theory to determine the proportion of susceptible population targeted for vaccination in a region.

On the other hand, the model proposed in \citet{jadidi2021two, westerink2017mathematical} aim to maximize the vaccinated population size and population achieving herd immunity using a herd immunity function (see \eqref{eq: herd-immunity}). The herd immunity function, Herd$(f_r)$ (e.g., Lambert herd function \citep{westerink2017mathematical}) is described as a variable of the fraction of total vaccinated population in region $r$.
\begin{align}\label{eq: herd-immunity}
&  \max \; \underbrace{\sum_{r} \hat{P}_r f_r}_{\substack{\text{Population size vaccinated} \\  \text{over all regions.}}} + \; \underbrace{\sum_{r} \hat{P}_r \,\text{Herd}\; (f_r)}_{\substack{\text{Population size benefiting from} \\ \text{herd immunity over all regions.}}} \\
&  \,\, \text{s.t.} \underbrace{\sum_{r} \hat{P}_r f_r \leq C^V}_{\substack{\text{Vaccinated Population cannot} \\ \text{outnumber vaccine availability.}}} \quad \underbrace{0 \leq  f_r \leq S_r/\hat{P}_r, \hspace{0.25cm} \forall r}_{\substack{\text{At most susceptible} \\ \text{population fraction} \\ \text{can be vaccinated.}}}  \nonumber
\end{align}

Instead of using a herd immunity function, \citet{macintyre2022modelling, rahman2023optimising} use vaccine coverage coefficient,  $\frac{1}{\beta} (1-\frac{1}{\alpha})$, as a measure for herd immunity. When vaccine effectiveness $\beta$ is small and nominal infection rate $\alpha$ is large, vaccine coverage coefficient becomes larger than 1, implying that we need to vaccinate more people to achieve the desired herd-immunity level. Vaccine effectiveness is also addressed in \cite{mohammadi2022bi} while minimizing expected number of deaths resulting from unvaccinated, single-dose and double-dose vaccinated groups. 

\srevision{
\subsection{Vaccine Roll-out Policy Evaluation With SIR Dynamics}
\citet{hu2023first} propose a vaccine allocation model motivated by the trade-offs between partially vaccinating a segment of population versus fully vaccinating half of that segment by splitting each full vaccine dose into two under supply scarcity. Their work considers an SIR model with two vaccinated compartments, containing individuals with fractional and full dosage. \citet{mak2022managing} evaluated vaccine roll-out policies such as hold-back, release, and single-dose-first strategies using SEIR dynamics, focusing on symptomatic cases, hospitalizations, and deaths.}

\srevision{
\section{Methodological Considerations and Case Studies} \label{sec:Methods}

\subsection{Methodological Considerations} 

Initial studies for COVID-19 vaccine supply chain management literature predominantly used commercial solvers to solve the analytical models developed. As the pandemic evolved, researchers integrated more realistic features into existing models or introduced new representations of previously studied features while adapting existing solution approaches. The following sections outline decomposition-based global optimization approaches (Section \S \ref{sec:811}), decomposition-based heuristic approaches (Section \S \ref{sec:812}), other heuristic and meta-heuristic approaches (Section \S \ref{sec:813}), non-linear optimization approaches (Section \S \ref{sec:814}) and multi-objective optimization algorithms (Section \S \ref{sec:815}). 

\subsubsection{Decomposition and Global Optimization Approaches} \label{sec:811}

These methods guarantee solution quality using a defined value for the optimality gap. \citet{azadi2020developing} use the L-shaped method to solve their two-stage problem. Motivated by their problem structure, \citet{zhang2022mass} utilize logic-based Benders to decompose location-assignment and appointment scheduling decisions. \citet{wang2023robust} use a column-and-cut generation (C\&CG) algorithm to solve their two-stage robust optimization problem, working with both primal and dual forms of C\&CG. The primal approach requires exploring vertices of a specially structured uncertainty set, while a vertex-traversing algorithm helps obtain the worst-case scenario among all vertices. The dual approach dualizes the innermost problem, and the resultant bilinear problem is reformulated as a mixed-binary linear program by exploiting a vertex-traversal algorithm. This approach avoids the traditional Big-M approach. \citet{enayati2020optimal} use the multi-parametric-disaggregation technique, a global optimization approach, to solve their proposed bi-linear problem. This method introduces two new sets of variables (continuous and binary) for each bilinear term, in addition to a variable representing each bi-linear term. These representative variables are then expressed as a summation of a continuous set of variables with varying parameters as coefficients. 

\subsubsection{Decomposition-based Heuristic Approaches} \label{sec:812}
These methods, while not assuring solution quality like global optimization, can still efficiently achieve good solutions depending on the problem and decomposition type.  \citet{lai2021multi} utilize a Benders decomposition-based heuristic to solve a two-stage stochastic model. Their method solves the LP relaxation of a mixed-integer linear subproblem to generate an optimality cut. Once a targeted optimality gap is reached, the fractional solution of the second stage is rounded up for integer components providing the heuristic solution. \citet{georgiadis2021optimal} treat hub/DC-VC level allocation as a subproblem in a vaccine supply chain network from manufacturer to VC. Their method first allocates VCs to nearby hubs and then condenses the model by clustering VCs within existing political boundaries, subsequently addressing the smaller hub-VC assignment problems with aggregated parameters. By treating binary solutions from this assignment as incumbent, original large-scale MILP subproblems become LPs, consequently reducing computational burden. 

The heuristic decomposition method by \citet{yang2021design} first breaks down the main network using a hierarchical clustering technique, ensuring the resulting sub-network can handle a MILP problem. After solving the sub-problem corresponding to a sub-network, the method sequentially aggregates other sub-networks, keeping certain optimal hub-opening decisions fixed, while allowing adjustments to others. The nodes that are reconsidered as decision variables in the merged sub-networks depend on their positions. The aggregation-based algorithm by \citet{karakaya2023developing}, on the other hand, performs aggregation on the priority groups. The algorithm for priority group based vaccination timeline decision simultaneously deals with two groups: the highest priority group and the group of all others post-aggregation. The timeline established for the first group serves as a parameter in the subsequent iteration, which focuses on the second priority class and the rest in an aggregated form. To prevent potential delays in vaccinating the aggregated group due to its high demand (due to aggregation), the minimum coverage threshold is also suitably adjusted, guaranteeing an optimal timeline for each class at the end. \citet{karakaya2023developing} also offer comprehensive guidance on modifying the initial vaccination calendar in response to unforeseen events, such as severe vaccine side effects during the planned period.

\subsubsection{Other Heuristics and Meta-Heuristics} \label{sec:813}

Heuristics are problem-specific strategies or rules of thumb that guide the search for scalable solutions but offer no performance guarantee and sometimes result in very poor solution quality. In comparison to heuristic approaches, meta-heuristics often yield better solutions, but they are generally computationally more intensive for exploring a large search space. 

A commonly-practiced algorithm by local and state governments during the COVID-19 pandemic was the pro-rata heuristic, which divides the available vaccines equally based on the population size. The vaccine allocation heuristic by \citet{orgut2023equitable} distributes vaccines to population sites following capacity-to-demand ratio in ascending order. If a site's demand is unmet, the heuristic checks neighboring sites in the order of proximity, adjusting available capacity with each allocation. After fulfilling a region's entire demand, the process repeats for the next region, using the updated capacity-to-demand sequence. The greedy approach outlined in \citet{liquantifying} ranks vaccination candidates based on a value derived from multiplying the associated risk measure, the likelihood of post-infection fatality, and the reduction in fatality rate after vaccination. Vaccine needs of a selected candidate are then fully met unless constrained by vaccine availability or budget. Although such an approach is optimal when restricted to a single period design, due to its myopic nature, it can be always far from an optimal allocation in a multi-period framework.          

\citet{enayati2023vaccine},  \citet{enayati2023multimodal} in their layered feasible path network  formulation leverages depth-first-search algorithm to construct all viable paths between each origin-destination pair. Subsequent heuristic reduces layered network size through parameters that control the number of transshipments, and drone launch-landings in a vaccine delivery path. Thus, the computational cost to solve the resultant model is reduced. In \citet{zhang2022mass}, a meta-heuristic approach is proposed for a combined location-assignment and scheduling problem, which mitigates NP-hardness of the scheduling subproblem by using a heuristic sequencing policy. This policy allows for the elimination of binary variables that determine task sequences. The heuristic prioritizes tasks based on the smallest values of the earliest start and latest completion times. Consequently creating an approximation of the original MILP is created, which, in addition to other feasibility constraints, includes constraints for the start times of each vaccination task utilizing the predetermined sequences.

\subsubsection{Nonlinear Optimization Approaches} \label{sec:814}
A two-stage stochastic formulation of the facility location problem with location preference order by \citet{cabezas2021two} uses the accelerated dual ascent (ADA) algorithm. ADA in the first step solves the Lagrangian dual problem by a subgradient algorithm while the semi-Lagrangian-dual problem in the second step is solved by a dual ascent algorithm using the solution from the first step, followed by a heuristic for accelerating the entire procedure. Iterative coordinate descent algorithm, which takes a gradient descent step with respect to only one variable at a time, is used in \citet{bertsimas2020optimizing, bertsimas2022locate} to solve the bilinear model, resulting from epidemiological dynamics considerations. \citet{liquantifying} also employed the coordinate descent method, alternating between assignment variable and the product term associated to risk over multiple time periods, exploiting their multiplicative relationship. Among other methods, \citet{barthspatiotemporal} used sequential quadratic programming (SQP) to solve a discretized version of continuous time model that, alongside bilinearity also involves non-linearity resulting from time-dependent contact-rate considerations. We note that these non-linear optimization approaches for non-convex problems typically converges to a local solution, but the quality of the solution can be enhanced by exploring different initial solutions.

\subsubsection{Algorithms for Multi-Objective Optimization} \label{sec:815}

Among the types of algorithms used in multi-objective vaccine supply chain literature, especially including the sustainability criteria, weighted-sum approach, and epsilon-constraint approach are the two simplest ones used (e.g., \citet{tang2022bi}). Nature-inspired and evolutionary-type algorithms, also known as meta-heuristics, appear next in the list of most common algorithm types. Under this category, grey-wolf optimization and its modifications (e.g., \citet{goodarzian2022designing}), particle swarm optimization (PSO) (e.g., \citet{chowdhury2022modeling}), teaching-learning-based optimization (TLBO) and its hybrid with PSO and genetic algorithm (GA) (e.g., \citet{goodarzian2022sustainable}), non-dominated sorting genetic algorithm II (NSGA-II) (e.g., \citet{tang2022bi}), evolutionary strategies (ES) (e.g., \citet{lim2022redesign}), can be found in the literature. Other methods include, but not limited to, variable neighborhood search and LP-metric method (e.g., \citet{goodarzian2022designing}). We note that constructing a Pareto frontier is an integral part of these works. We also note that evolutionary type algorithms are better regarded for building a Pareto frontier with a diverse set of non-dominated solutions compared to other multi-objective optimization approaches used.

\subsection{Insights from Numerical Studies}
In this section, we outline key highlights from representative numerical or case studies in the literature. While some studies use synthetic data to emphasize their main points, others integrate real data (e.g., \citet{hajibabai2022using}), supplemented with rough approximations. These studies either provide an artificial view of national vaccine supply chains, like those in the U.S. or India, or focus on specific regional details. In deterministic modeling, some papers validate models using a range of artificially generated parameter data and assessing their interpretational consistency, whereas under uncertainty, model outcomes are assessed under various uncertainty scenarios. Another aspect of model validation, observed in a few studies, involves stating modeling assumptions and evaluating their reasonableness in reflecting real-world conditions without oversimplification or misrepresentation. Regardless, sensitivity analysis is often performed with respect to central parameters such as budget in facility number, cost budget, vaccine availability, supplier/DC/VCs capacity, healthcare personnel service rate, minimum coverage threshold, minimum service ratio, cost-coefficient, parameters defining uncertainty/ambiguity set, etc.

Case studies also highlight the importance of using an optimal allocation policy compared to the pro-rata policy by using metrics such as region or risk class-wise number of active cases, deaths, (e.g., \citet{bertsimas2020optimizing}), infected days (infection rate multiplied by length of time horizon in days) (e.g., \citet{barthspatiotemporal}, percentage of unused vaccines, percentage of unmet demand (e.g., \citet{orgut2023equitable}), etc. For example, \citet{bertsimas2020optimizing} showed that optimized vaccine allocation following their model can reduce 10-25\% of the fatalities over the 90-day planning horizon. \citet{balcik2022mathematical} compared several coverage-based key performance indicators (KPIs) across allocation policies systematically, factoring in priority groups and/or socio-economic differences. \citet{barthspatiotemporal} additionally assessed the efficiency of different heuristic vaccine allocation strategies, remarking on the consistency of the pro-rata policy compared to others. For more on such policies, see Section 4.1 of \citet{barthspatiotemporal}. 

\subsubsection{Insights from Parameter Variation Effect}
\citet{mohammadi2022bi} demonstrated how decisions such as the number of vaccine centers, vaccine orders, transfers, trans-shipments, and first-dose and second-dose vaccinations per risk class vary across different uncertainty scenarios. Their study revealed that as vaccine campaign parameters (like warehouse capacity, dose time-lag, and manufacturer trust rate) enhance, deaths decrease and first-dose vaccinations increase. Longer vaccine delivery lead times, however, have the opposite effect. Mortality also increases with increased contact and death rates but reduces with vaccine effectiveness. The authors suggest prioritizing high-risk individuals for vaccination when supplies are limited but high-contact individuals otherwise. Vaccination strategy and dosing intervals depend on the situation: in the face of a lethal and highly contagious virus variant, it is important to fully vaccinate the high-risk group rather than just giving one dose to many. Boosting vaccination capacity can combat such variants, and resource sharing offers flexible capacity planning. 

Sensitivity analyses in the literature explicitly show that equity can be improved via a change in facility capacity (e.g., \citet{balcik2022mathematical}, \citet{orgut2023equitable}). In addition to evaluating fairness and efficiency through case studies, many studies assess other performance metrics as indicated earlier in this section. When optimization models introduce or emphasize specific features, assessing their outcomes numerically is standard to highlight their importance. For instance, \citet{balcik2022mathematical} evaluated the fairness measure they use against other known equity measures. \cite{barthspatiotemporal} analyzed how the ratio of available vaccine percentage to population size affects infected days, Gini-coefficients of vaccinations as well as Gini-coefficients of infected days.}

\srevision{\section{Discussion} \label{sec:conclusion}
This section provides a comprehensive discussion of the above discussed models and their usefulness in a future pandemic. The discussion includes unique features of a vaccine supply chain during a pandemic (\S \ref{sec:91}), trends in the literature (\S \ref{sec:92}), and authors' perspectives for model selection (\S\ref{sec:94}).

\subsection{Unique Challenges in Pandemic Vaccine Supply Chain} \label{sec:91}
The vaccine supply chain presents unique challenges that distinguish it from other supply chains. First, it requires rigorous cold chain management, especially for vaccines necessitating ultra-cold storage. Secondly, it necessitates the development of adaptable distribution strategies that take into account limited manufacturing sites and unpredictable epidemiological patterns. Third, complex inventory management is needed due to diverse factors such as vaccine dose policies, vial sizes, vial shelf life, open vial shelf life, wastage, and the need for specialized healthcare staff. Additionally, the supply chain should maintain a delicate balance between equity and efficiency while addressing vaccine hesitancy and priority groups. Extending vaccine access to remote areas poses another hurdle, wherein drone delivery does not fully resolve the issue of healthcare personnel availability. Moreover, the integration of demand into allocation models is critical, considering vaccine preferences and epidemiological dynamics in quantifying demand. 

\subsection{Temporal Trends in the Literature} \label{sec:92}
In consideration of the above mentioned distinctions, several temporal trends in the literature can be noted. In the initial phases of the pandemic, models considered are similar to those of conventional supply chains, albeit with some added emphasis on cold chain management. Additionally, real-life case studies were rare, most probably due to limitations in the availability of data. However, as the pandemic progressed, the number of real-life case studies increased, providing a rich repository of insights and experiences.  In this stage, researchers also produced more realistic models that consider different sources of uncertainty and epidemiological dynamics, enhancing the relevance of these tools in public health decision-making. 

A later wave of studies has aimed to advance equity and efficiency in vaccine deployment beyond that enabled by prior work. Recent work also features models that address vaccine efficacy in preventing transmission and infections, as well as vaccine hesitancy to take the vaccine among different segments of the population. Hence, we see a strong interest in modeling characteristics pertaining to human behavior and social dynamics as well as clinical performance of the developed vaccines against mutations of the vaccine. 

Despite these advances, however, there are still important issues to study. For example, the more general question of how to prioritize the inoculation of various segments of the population (e.g., fractional dosing, prioritizing individuals of high-risk occupations/social interaction) while at the same time considering realistic models of vaccine characteristics (e.g., type, dosing, storage requirements) population mobility, social determinants of health, and access to healthcare are lacking. 

Regardless of the specific circumstances, it is clear that vaccine allocation models that incorporate epidemiological dynamics and realistic mobility behavior can greatly enhance the effectiveness of public health decision-making. However, even without the explicit consideration of these aspects of the problem, optimization models that facilitate resource sharing among different geographical areas and possess the flexibility to adapt to changing environments can still prove to be highly valuable.

\subsection{Perspectives for Model Selection} \label{sec:94}

Models that primarily address location selection and opening decisions, as these are strategic choices that are not easily altered once made, should incorporate different practical factors. In parallel, there should be solution approaches that provide solutions of high quality. On the other hand, operational-level decision-making models, which prioritize enhancing coverage in hard-to-reach regions or fortifying vulnerable supply chains to mitigate the aftermath of a major impact, warrant attention for enhancing global vaccination efforts. Furthermore, it is advisable to develop models that considers epidemiological dynamics, either via compartment or network structure, to account for long-term effects. These models should focus on the time-varying nature of a pandemic and should be capable of optimizing resource allocations, while also considering factors such as vaccine availability, quarantine, lockdown measures, and other preventive strategies. Although such models may pose computational challenges, they are deemed to be ideal, even when epidemiological dynamics is represented in approximate or simulated forms.

Effective utilization of Operations Research often demands a high level of expertise, and without it, there exists a risk of improperly applying models and techniques. Hence, models should be carefully developed, providing clear and simple guidelines for public health officers to facilitate implementation and execution. Additionally, caution is required as recommendations from models based on unreliable or insufficient data can deviate significantly from reality, particularly in scenarios with rapidly changing dynamics such as the COVID-19 pandemic. It is also worth noting that vaccine supply chains can be influenced by cultural, social, or behavioral factors, which may prove challenging to quantify within a model's framework. Moreover, if the chosen solution approach to address the modeled problem lacks quality or fails to align with the specific purpose, the effort and expenses invested in the model may outweigh the benefits it provides.}

\srevision{
\section{Concluding Remarks on Future Research Directions} \label{sec:Critic}

A pandemic vaccine supply chain can be viewed in three phases: Phase-I represents the initial period following vaccine development, characterized by a pronounced gap between vaccine supply and demand for inoculation. During Phase-II, the supply of vaccines increases, meeting a substantial portion of the demand but still falling short of full coverage. Phase-III defines the state in which the vaccine supply matches or even surpasses the demand. Note that the pandemic vaccine deployment literature, including the planning of distribution center (DC) and vaccination center (VC) locations, is predominantly centered on problems pertaining to Phase-II. Below, we include some concluding remarks on fruitful research directions in modeling (\S \ref{sec:101} and \S \ref{sec:102}) and solution methods (\S \ref{sec:103}).

\subsection{Models with Explicit Consideration of Epidemiological Dynamics} \label{sec:101}

Traditionally, an estimated size of the susceptible population is regarded as the demand in vaccine allocation models. However, several factors, such as the administration of multiple vaccine doses, variations in vaccine efficacy, and fundamental epidemiological dynamics strongly influence the size of the population susceptible to infection over time. Hence, it is imperative that such decision making models incorporate effective epidemiological dynamics models to estimate the demand for vaccines over time in different geographical regions that may be going through different phases of the pandemic. 

Although several papers incorporate compartmentalized epidemiological models (e.g., SIR-type models) into vaccine allocation models, present studies exhibit several limitations. As discussed in Section~\ref{sec:SEIRD models}, directly integrating SIR dynamics into optimization models or employing separate models with diverse SIR parameters results in computational complexity or challenges in achieving interpretability. This leads to a critical question: is it possible to develop a computationally efficient SIR-surrogate that preserves essential epidemiological dynamics?

Moreover, the ever-evolving nature of virus challenges the accuracy of models incorporating epidemiological dynamics, even using advanced machine-learning models. Addressing this uncertainty underscores the need for using multiple prediction models. However, current studies lack a unified approach for integrating these models into optimization algorithms. Additionally, concerning the dynamic environment of a pandemic, it is worth exploring whether prediction and optimization should be performed sequentially or concurrently in a decision-making framework. Furthermore, in a sequential decision-making framework, the potential utilization of an integrated approach that combines offline and online learning for a fast-yet-more informed decision-making warrants new studies. The question of how contextual optimization model outcomes, modeled either with real or simulated data, respond to contextual changes should also be examined. The OR/MS community needs to address these critical issues to advance the field beyond superficial adjustments to carefully-studied vaccine supply chain models with potential for real-world impact.

\subsection{Additional Modeling Considerations} \label{sec:102}

Several additional gaps exist in the modeling of vaccine supply chain models. The uncertainty-driven dynamism of these problems poses the question of whether distribution networks aided by novel comprehensive topological metrics and graph representation learning can yield decisions that are more adaptable. Research examining the interplay between quarantine/isolation policies and vaccine allocation is notably scarce, as is the development of models that account for time-varying contact and mortality rates. Recognizing that vaccination does not uniformly reduce the risks of infection and mortality across various subpopulations in different risk categories, equity metrics based on demand or coverage may prove inadequate, prompting the exploration of alternative methods. Exploring vaccine rollout strategies that consider time-varying contact rates and the occurrence of multiple infection spikes, and utility-driven location models that factor in diverse determinants such as distance, regional income, and education, present promising avenues for further research. Modeling regulated by risk-neutral or risk-aversion type utility can further serve as a viable alternative to conventional vaccine hesitancy models.

Current discrete scenario-based approaches to two-stage stochastic programming fall short of addressing multiple spikes of infections during a pandemic. Multi-stage stochastic programming, while computationally intensive, could more accurately capture epidemiological dynamics over time. The challenge lies in devising effective surrogates for this approach and rigorously assessing their performance.

Finally, further developments on model validation approaches are required for ascertaining the adaptability of a vaccine supply chain network to rapid shifts in infection dynamics. Defining criteria for allocating limited vaccine supplies in a prioritized order during the initial vaccination phase (Phase-I) necessitates systematic analysis, streamlining long-term risks and benefits, while considering access to healthcare, social determinants of health, human behaviour, mobility and adherence to non-pharmaceutical interventions such as social distancing and vaccine hesitancy.

\subsection{Solution Methodologies} \label{sec:103}

The majority of reviewed work directly adopts solution methods from similar literature, which, regrettably, hampers the advancement of comprehensive models equipped with practical features. To address this gap, the OR/MS community needs to innovate effective solution strategies that leverage the inherent problem structure, decomposition techniques (e.g., graph-based methodologies), and valid inequalities synchronized with model condensation techniques, and harness combinatorial heuristic algorithms to solve realistic problems with large network structures. Such methods will be particularly crucial to efficiently solve multi-objective models with conflicting objectives as well as bilinear/biconvex models resulting from robust optimization formulations that are needed to represent the plethora of sources of uncertainties in this problem domain.
}

\ACKNOWLEDGMENT{This research is supported by the NIH Grant R01AI168144. We would like to thank Jiaqi Lei for helping proofread this article. \textcolor{black}{We also sincerely thank the Associate Editor and the two anonymous reviewers for their insightful comments and recommendations, which allowed us to improve the manuscript significantly.}}
 
\bibliographystyle{plainnat}
\bibliography{references}

\newpage

\begin{landscape}

\begin{figure}[htp]
\centering
\includegraphics[width=1.4\textwidth]{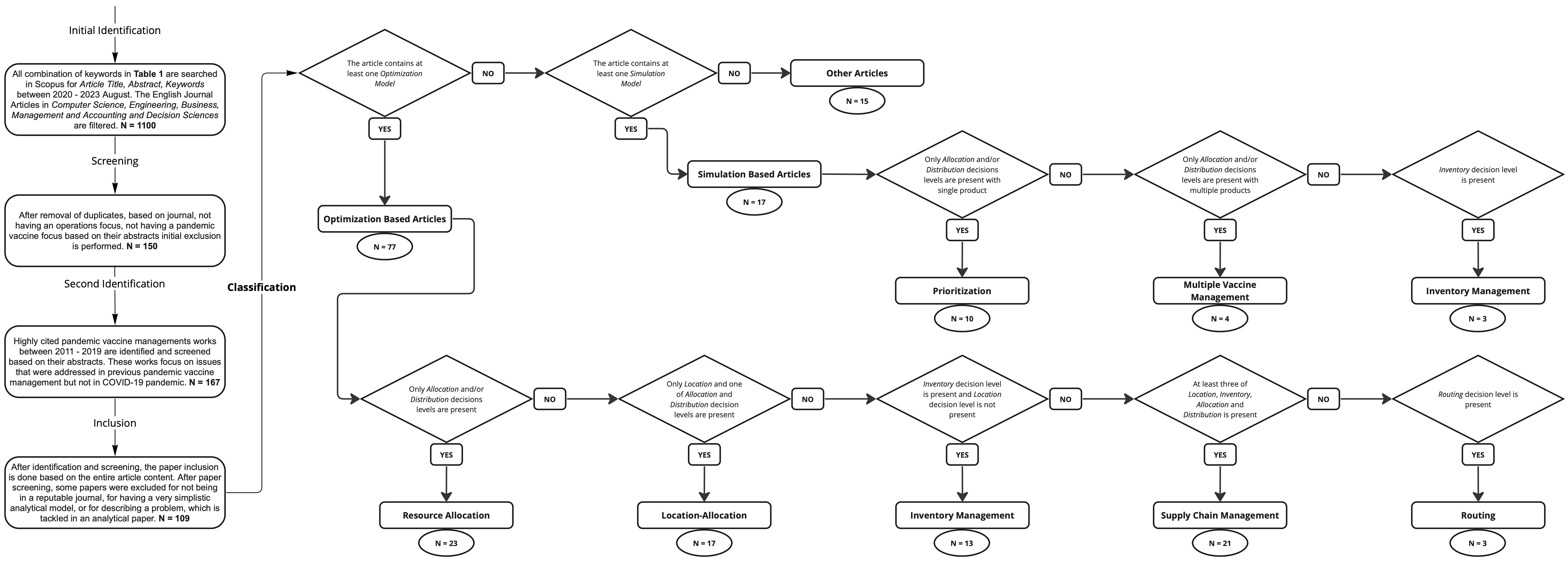}
\caption{Flowchart of the process used to select and classify papers on pandemic vaccine supply chain management.}
\label{fig:PRISMA}
\end{figure}

\end{landscape}

\newpage

\begin{APPENDIX}{Detailed Features of Optimization and Simulation Articles in Pandemic Vaccine Management}

\begin{table}[htp] \setlength{\tabcolsep}{8pt}
\begin{center}
\caption{List of abbreviations used in Appendix.}
\begin{adjustbox}{width=\textwidth}
\begin{normalsize}
\renewcommand{\arraystretch}{1}
\begin{tabular}{llll}
\hline
\rowcolor{Gray}
\textbf{Model} & & & \\
\hline
CCP: Chance Constraint Programming  & CoP: Copositive Programming  & CSP: Constraint Satisfaction Programming  & DT: Decision Trees  \\
FP: Fuzzy Programming & GT: Game Theory  & I: Inventory Modeling  & LP: Linear Programming  \\
LRP: Location Routing Problem & MC: Markov Chains  & MILP: Mixed Integer Linear Programming  & MINLP: Mixed Integer Non-linear Programming  \\
NLP: Non-linear Programming & OC: Optimal Control & Q: Queuing  & RL: Reinforcement Learning  \\
RO: Robust Optimization  & S: Simulation & SCP: Set Covering Problem  & SN: Social Networks \\
SP: Stochastic Programming  & U: Utility Function  & VRP: Vehicle Routing Problem &  \\
\hline
\rowcolor{Gray}
\textbf{Uncertainty} & & & \\
\hline
D: Demand & DR: Deterioration Rate & I: Infection & VLC: Vaccination Location Choice \\
LT: Lead Time & S: Supply & VA: Vaccine Accessibility & VW: Vaccine Wastage \\
\hline
\rowcolor{Gray}
\textbf{Commodity} & & & \\
\hline
\textit{Quantity:} & M: Multiple & S: Single & \\
\hdashline
\textit{Type:} & ChV: Childhood Vaccine & CP: COVID-19 Pills & CT: COVID-19 Test \\
& CV: COVID-19 Vaccine & CW: COVID-19 Waste & GHS: General Health Service \\
& GV: General Vaccine & IV: Influenza Vaccine & \\
\hdashline
\textit{Characteristics:} & MD: Multiple Doses & P: Perishable & SD: Single Dose \\
\hline
\rowcolor{Gray}
\textbf{Period} & & & \\
\hline
M: Multiple & S: Single & & \\
\hline
\rowcolor{Gray}
\textbf{Echelon} & & & \\
\hline
DC: Distribution Center & GHSC: General Health Service Center & M: Manufacturer & MVC: Mobile Vaccination Clinic \\
OC: Outreach Center & PS: Population Sub-group & RT: Refrigerated Truck & VC: Vaccination Center \\
WDC: Waste Disposal Center & & & \\
\hline
\rowcolor{Gray}
\textbf{Objective Function} & & & \\
\hline
\textit{Quantity:} & [Max]: Maximize & [Min]: Minimize & [$n$]: $n$ Different Separate Objective Functions \\
\hdashline
\textit{Operator:} & M: Multiple Objectives & S: Single Objective & \\
\hline
\rowcolor{Gray}
\textbf{Vehicle} & & & \\
\hline
\textit{Type:} & MVC: Mobile Vaccination Clinic & RT: Refrigerated Truck & \\
\hdashline
\textit{Characteristics:} & He: Heterogeneous & Ho: Homogeneous & \\
\hline
\end{tabular}
\end{normalsize}
\end{adjustbox}
\label{table:abbreviations}
\end{center}
\end{table}

\begin{figure}[htp]
\centering
\includegraphics[width=0.5\textwidth]{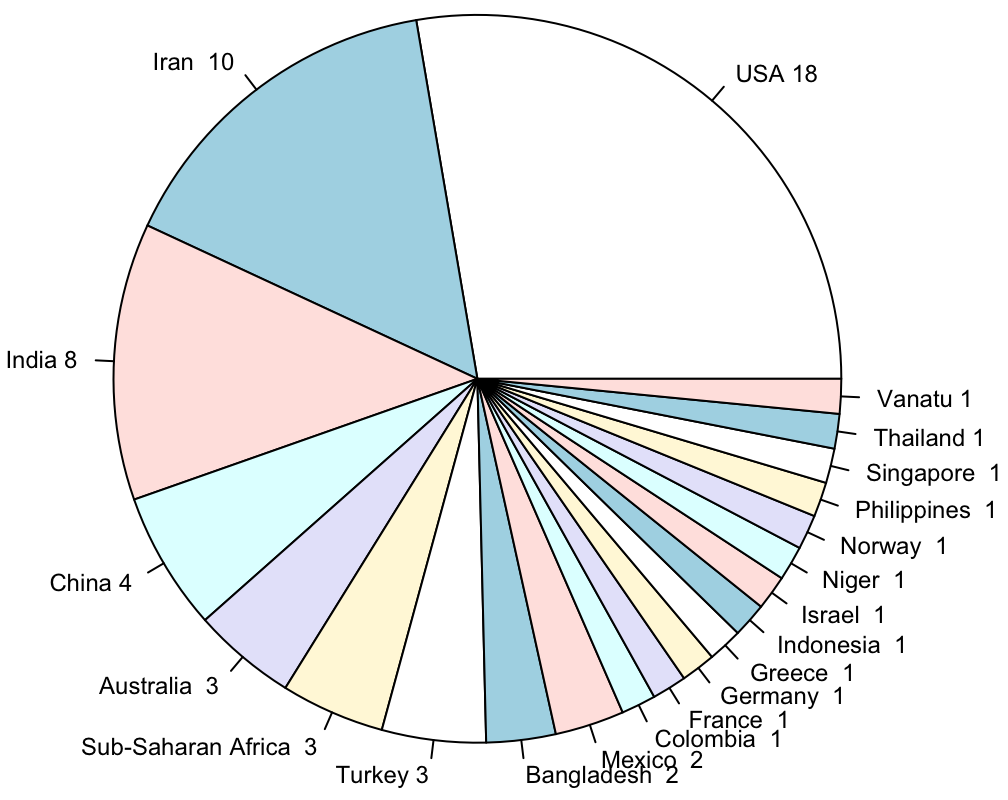}
\caption{Number of real-life case studies used by optimization-based articles according to countries.}
\label{fig:CaseStudies}
\end{figure}

\begin{table}[!hbt] \setlength{\tabcolsep}{8pt}
\begin{center}
\caption{Detailed features of the optimization articles studying resource allocation problems.}
\begin{adjustbox}{width=\textwidth}
\renewcommand{\arraystretch}{1.25}
\begin{tabular}{p{1.5in}|p{0.6in}|p{0.25in}|p{0.25in}|p{0.80in}|p{0.25in}|p{0.25in}|p{0.65in}|p{2.5in}|p{1in}|p{1.5in}}
\hline
\rowcolor{Gray}
\rotatebox{270}{\textbf{Article}} & \rotatebox{270}{\textbf{Model}} & \rotatebox{270}{\textbf{Prioritization}} & \rotatebox{270}{\textbf{Uncertainty}} & \rotatebox{270}{\textbf{Commodity}} & \rotatebox{270}{\textbf{Period}} & \rotatebox{270}{\textbf{Echelon: From}} & \rotatebox{270}{\textbf{Echelon: To}} & \rotatebox{270}{\textbf{Objective}} & \rotatebox{270}{\textbf{Case Study}} & \rotatebox{270}{\textbf{Solution Method   }} \\
\hline

\citet{abbasi2020modeling} & MILP & $\checkmark$ & $\times$ & S,CV,SD & S & DC & PS & S, [Min] Vaccine oversupply, transshipment time. & Victoria, Australia & Risk-based heuristic algorithm. \\
\hline
\citet{anahideh2022fair} & MILP & $\times$ & $\times$ & S,GV,SD & S & DC & VC & M, [Min] Geographical diversity, [Min] social fairness gap. & Baltimore, Chicago, New York, USA & LP Relaxation, Penalty Method, Binary Search \\
\hline
\citet{balcik2022mathematical} & MILP & $\checkmark$ & $\times$ & M,CV,SD & S & DC & DC & S, [Min] Deviation from the fair coverage levels. & Turkey & Commercial Solver \\
\hline
\citet{bandi2021optimal} & CoP & $\checkmark$ & D & M,CV,MD,P & M & VC & PS & M, [Max] Vaccination rates, [Min] appointment waiting times. & Singapore & Copositive Program Reformulation \\
\hline
\citet{barthspatiotemporal} & OC & $\checkmark$ & D & S,CV,SD & M & VC & PS & S, [Min] Total number of infected people. & Synthetic & NLP reformulation by discretizing the time dimension. \\
\hline
\citet{bennouna2022covid} & SP & $\checkmark$ & D & M,CV,MD & M & DC & PS & S, [Min] Expected number of deaths. & Massachusetts, USA & ML-based aggregation method. \\
\hline
\citet{bertsimas2020optimizing} & NLP & $\checkmark$ & D & S,CV,SD & M & DC & PS & S, [Min] Number of deaths. & New York, USA & Iterative Coordinate Descent \\
\hline
\citet{enayati2020optimal} & NLP & $\checkmark$ & D & S,IV,SD & M & DC & PS & S, [Min] Number of vaccine doses distributed. & Synthetic (Asian Flu Pandemic) & Discretization with Multi-parametric Disaggregation Algorithm \\
\hline
\citet{hu2023first} & OC & $\checkmark$ & D & S,CV,SD & M & VC & PS & S, [Min] Total number of infected people. & World & Optimal Control with Binary Search \\
\hline
\citet{huang2017equalizing} & MILP & $\checkmark$ & $\times$ & M,IF,SD & S & DC & DC & S, [Max] Equal access. & Texas, USA & Commercial Solver \\
\hline
\citet{jadidi2021two} & NLP,S & $\checkmark$ & D & S,CV,SD & S & DC & PS & S, [Max] Total immunity among populations. & Synthetic (Another Paper) & Graph Theory Algorithm \\
\hline
\citet{jarumaneeroj2022epidemiology} & NLP & $\checkmark$ & D & M,CV,SD & M & DC & PS & S, [Min] Stress on healthcare system. & Thailand & Iterative Algorithm \\
\hline
\citet{liquantifying} & NLP,SN,S & $\checkmark$ & $\times$ & S,CV,SD & M & VC & PS & S, [Min] Number of total fatalities. & Texas & Greedy solution scheme, solution algorithm using coordinate descent. \\
\hline
\citet{minoza2021covid} & LP,S & $\checkmark$ & D & S,CV,SD & S & DC & PS & M, [Max] Vaccine to be allocated for each location, [Max] priority factor. & Quezon City, Philippines & Simulating Scenarios \\
\hline
\citet{munguia2021fair} & LP,NLP & $\checkmark$ & $\times$ & S,CV,SD & S & DC & DC & [\textit{3}], [Max] Total vaccines allocated, [Max] smallest allocated vaccines, [Max] social welfare. & Mexico & Commercial Solver \\
\hline
\citet{orgut2023equitable} & MILP & $\checkmark$ & $\times$ & S,CV,SD & S & DC & VC & S, [Min] Maximum travel duration of patients. & Alabama & Pro Rata Vaccine Allocation Heuristic \\
\hline
\citet{rao2021optimal} & MILP & $\checkmark$ & D & S,CV,SD & S & DC & PS & [\textit{4}], [Min] New infections, [Min] deaths, [Min] life years lost, [Min] quality-adjusted life years. & New York, USA & Exact Analytical Conditions \\
\hline
\citet{roy2021optimal} & LP & $\times$ & I & S,CV,SD & M & DC & PS & S, [Min] Transportation cost. & New York, USA & Exact Solution Algorithm \\
\hline
\citet{shukla2022optimizing} & MILP & $\checkmark$ & $\times$ & M,CV,SD & M & DC & VC,MVC & S, [Max] Vaccine distribution efficiency. & USA & Python NetworkX, Minimum Cost Flow Problem \\
\hline
\citet{thul2023stochastic} & SP,RL & $\checkmark$ & D & S,CV,SD,CT & M & DC & PS & S, [Max] Expected sum of cumulative rewards. & USA & Partially-Observable Markov Decision Process \\
\hline
\citet{yang2022comparison} & SP & $\checkmark$ & D & S,CV,SD & S & DC & PS & S, [Min] Total number of confirmed cases. & Synthetic & Commercial Solver \\
\hline
\citet{yarmand2014optimal} & SP & $\times$ & VA & S,GV,SD & M & DC & PS & S, [Min] Phase-I vaccination cost, expected Phase-II vaccination cost. & North Carolina & Newsvendor Closed Form Solution, Heuristic Algorithm \\

\hline
\end{tabular}
\end{adjustbox}
\label{table:resource-allocation-articles}
\end{center}
\end{table}

\begin{table}[!hbt] \setlength{\tabcolsep}{8pt}
\begin{center}
\caption{Detailed features of the optimization articles studying location-allocation problems.}
\begin{adjustbox}{width=\textwidth}
\renewcommand{\arraystretch}{1.25}
\begin{tabular}{p{1.5in}|p{0.55in}|p{0.35in}|p{0.60in}|p{0.15in}|p{0.40in}|p{0.30in}|p{3in}|p{1in}|p{1.5in}}
\hline
\rowcolor{Gray}
\rotatebox{270}{\textbf{Article}} & \rotatebox{270}{\textbf{Model}} & \rotatebox{270}{\textbf{Uncertainty}} & \rotatebox{270}{\textbf{Commodity}} & \rotatebox{270}{\textbf{Period}} & \rotatebox{270}{\textbf{Facility Echelon}} & \rotatebox{270}{\textbf{Demand Echelon}} & \rotatebox{270}{\textbf{Objective}} & \rotatebox{270}{\textbf{Case Study}} & \rotatebox{270}{\textbf{Solution Method  
 }} \\
\hline
\citet{bertsimas2022locate} & NLP & D & S,CV,SD & M & VC & PS & M, [Min] Deaths, [Min] exposed, [Min] distance. & New York, USA & Iterative Coordinate Descent \\
\hline
\citet{bravo2022optimal} & MILP & $\times$ & S,CV,SD & S & VC & PS & S, [Max] Aggregate vaccinations. & California, USA & Commercial Solver \\
\hline
\citet{cabezas2021two} & SP & VLC & M,CV,SD & S & VC & PS & S, [Min] Demand coverage and facility cost. & Synthetic & Lagrangian Relaxation Heuristic \\
\hline
\citet{cao2023digital} & RO & I & S,CW & M & WDC & VC & M, [Min] Transportation and waste disposal cost, [Min] infectious risk. & Maharashtra, India & KKT, Branch and Bound, Polyhedral Uncertainty Set \\
\hline
\citet{dastgoshade2022social} & SP & D,S & S,CV,SD & M & VC & PS & M, [Min] Network cost, [Max] social equity. & Iran & Sample Average Approximation, Lexicographic Goal Programming \\
\hline
\citet{emu2021validating} & CSP & $\times$ & S,CV,SD & S & DC & PS & [\textit{4}], [Max] Overall vaccine distribution, [Max] vaccine distribution among the higher priority groups, [Max] vaccine distribution among people located closer to DCs, [Max] vaccine distribution among the higher priority groups and people located closer to DCs. & Chennai, India & K-medoids Algorithm \\
\hline
\citet{enayati2023multimodal} & MILP & $\times$ & S,CV,SD,P & S & DC & DC & S, [Min] Total operating costs. & Synthetic & Commercial Solver \\
\hline
\citet{enayati2023vaccine} & MILP & $\times$ & S,CV,SD & S & DC & VC & S, [Min] Total operating costs. & Vanatu Islands & Commercial Solver \\
\hline
\citet{kumar2022optimal} & MINLP,I & $\times$ & S,CV,SD & S & VC & PS & S, [Min] Operational (transportation, location, salary) costs and travel distance. & Chandigarh, India & Two-stage MILP \\
\hline
\citet{leithauser2021quantifying} & MILP & $\times$ & S,CV,SD & S & VC & PS & M, [Min] Number of vaccination centers, [Min] number of physicans, [Min] sum of patient travel distance. & Germany & Different Assignment Heuristics (Greedy, Same State) \\
\hline
\citet{lim2016coverage} & MILP & $\times$ & S,ChV,SD & S & OC & PS & S, [Max] Demand coverage. & Bihar, India & Commercial Solver \\
\hline
\citet{luo2023service} & RO,U & VLC & S,CT & S & VC & PS & S, [Max] Total utility of the service received by the customers. & San Diego, USA & Mixed 0-1 second-order cone (MISOCP) formulation and cutting-plane algorithm. \\
\hline
\citet{lusiantoro2022locational} & MILP & $\times$ & S,CV,SD & S & VC & PS & M, [Max] Coverage, [Min] weighted distance. & Yogyakarta, Indonesia & Commercial Solver \\
\hline
\citet{polo2015location} & MILP & $\times$ & S,GHS & S & GHSC & PS & M, [Max] Coverage, [Min] total distance. & Bogota, Colombia & Commercial Solver \\
\hline
\citet{soria2021proposal} & MILP & $\times$ & M,CV,SD & S & M,DC & DC & S, [Min] Total operational cost. & Mexico & Commercial Solver \\
\hline
\citet{srivastava2021strengthening} & MILP & $\times$ & S,CV,SD & S & DC & VC & S, [Max] Coverage of demand points, within a maximum time. & Madhya Prades, India & Greedy Adding Algorithm-Based Optimization \\
\hline
\citet{zhang2022mass} & MILP & $\times$ & S,CV,SD & S & VC & PS & S, [Min] Operational cost (location, traveling distance, appointment rejection, tardiness). & Yiwu, China & Logic Based Benders Decomposition and Metaheuristic \\
\hline
\end{tabular}
\end{adjustbox}
\label{table:location-allocation-articles}
\end{center}
\end{table}

\begin{table}[!hbt] \setlength{\tabcolsep}{8pt}
\begin{center}
\caption{Detailed features of the optimization articles studying inventory management problems.}
\begin{adjustbox}{width=\textwidth}
\renewcommand{\arraystretch}{1.25}
\begin{tabular}{p{1.5in}|p{0.45in}|p{0.4in}|p{0.7in}|p{0.2in}|p{0.8in}|p{2.5in}|p{1in}|p{2in}}
\hline
\rowcolor{Gray}
\rotatebox{270}{\textbf{Article}} & \rotatebox{270}{\textbf{Model}} & \rotatebox{270}{\textbf{Uncertainty}} & \rotatebox{270}{\textbf{Commodity}} & \rotatebox{270}{\textbf{Period}} & \rotatebox{270}{\textbf{Echelons}} & \rotatebox{270}{\textbf{Objective}} & \rotatebox{270}{\textbf{Case Study}} & \rotatebox{270}{\textbf{Solution Method  
 }} \\
\hline
\citet{azadi2020optimization} & CCP & D & M,ChV,SD & M & DC,VC & S, [Max] Number of fully immunized children, vaccine availability. & Niger & Sample Average Approximation \\
\hline
\citet{azadi2020developing} & SP & D & M,ChV,SD & M & DC,VC & S, [Min] Cost and open vial waste. & Bangladesh & Stochastic Benders Decomposition, L-Shaped Method \\
\hline
\citet{bonney2011environmentally} & I & $\times$ & S,GHS & S & DC & S, [Min] Environmental cost. & $\times$ & Closed Form \\
\hline
\citet{fadaki2022multi} & MINLP & $\times$ & S,CV,SD & M & DC,VC & M, [Min] Weighted risk of unvaccinated population, transshipment mechanism and return mechanism. & Victoria, Australia & Augmented Multi-period Algorithm \\
\hline
\citet{georgiadis2021optimal} & MILP & $\times$ & M,CV,SD,P & M & M,DC,VC & S, [Min] Distribution, inventory and wastage cost. & Greece & Two-Step Decomposition Algorithm \\
\hline
\citet{hovav2015network} & LP & $\times$ & S,IV,SD & M & M,DC,VC & S, [Min] Operational and shortage cost. & Israel & Commercial Solver \\
\hline
\citet{icsik2023optimizing} & RO & D,DR,S & S,CV,MD,P & M & DC,VC,WDC & M, [Max] Total vaccinations, [Min] holding, transhipment, waste management cost. & Izmir, Turkey & $\epsilon$-constraint Method  \\
\hline
\citet{jahani2022covid} & NLP,Q & S,LT & M,CV,SD & M & M,VC & M, [Min] Expected vaccine waiting time, [Min] holding and ordering cost. & Melbourne, Australia & Goal Attainment Algorithm, Non-dominated Sorting Genetic Algorithm \\
\hline
\citet{karakaya2023developing} & SP & S,LT & M,CV,MD & M & DC,PS & S, [Min] Deviations from the announced COVID-19 calendar. & Norway & Aggregation-based algorithm. \\
\hline
\citet{mak2022managing} & I & D & S,CV,MD & M & DC,VC & \textit{Three roll-out scenarios are compared:} (i) holding back second doses, (ii) releasing second doses, and (iii) stretching the lead time between doses. & Synthetic & Proof Based \\
\hline
\citet{mohammadi2022bi} & RO,Q & S & S,CV,MD & S & DC,VC,PS & M, [Min] Total expected number of deaths among the population, [Min] total distribution cost. & France & Scenario-based robust-stochastic optimization approach. \\
\hline
\citet{shah2022investigation} & I & VW & S,CV,SD,P & M & M,DC & S, [Min] Ordering, manufacturing, holding, inspection and carbon emission cost. & Synthetic & Exact Solution Algorithm \\
\hline
\citet{sinha2021strategies} & I,GT & D & S,CV,SD,P & M & DC,VC & [\textit{3}], [Min] Total cost, under no disruption scenario, [Min] total cost, under disruption, [Min] excess inventory cost. & Gorakhpur, India (Japanese Encephalitis Vaccine) & Heuristic Algorithm \\

\hline
\end{tabular}
\end{adjustbox}
\label{table:inventory-articles}
\end{center}
\end{table}

\begin{table}[!hbt] \setlength{\tabcolsep}{8pt}
\begin{center}
\caption{Detailed features of the optimization articles studying supply chain management problems.}
\begin{adjustbox}{width=\textwidth}
\renewcommand{\arraystretch}{1.25}
\begin{tabular}{p{1.5in}|p{0.8in}|p{0.3in}|p{0.7in}|p{0.2in}|p{1.2in}|p{3in}|p{1in}|p{3in}}
\hline
\rowcolor{Gray}
\rotatebox{270}{\textbf{Article}} & \rotatebox{270}{\textbf{Model}} & \rotatebox{270}{\textbf{Uncertainty}} & \rotatebox{270}{\textbf{Commodity}} & \rotatebox{270}{\textbf{Period}} & \rotatebox{270}{\textbf{Echelons}} & \rotatebox{270}{\textbf{Objective}} & \rotatebox{270}{\textbf{Case Study}} & \rotatebox{270}{\textbf{Solution Method  
 }} \\
\hline
\citet{abbasi2023designing} & MILP & $\times$ & S,CV,SD & S & M,DC,VC,PS,WDC & M, [Min] Economic cost, [Min] environmental cost. & Iran & Multi-objective grey wolf optimizer algorithm. \\
\hline
\citet{basciftci2023resource} & SP,RO & D,I & S,CV,SD & M & DC,PS & S, [Min] Total cost of locating DCs, installing their capacities, expected cost of shipping, holding inventory, back-orders. & Michigan, USA & Moment-based ambiguity set approach. \\
\hline
\citet{chowdhury2022modeling} & MILP & D,S & S,CV,SD,P & M & M,DC,VC & M, [Min] Operational cost, [Min] environmental cost, [Max] job opportunities. & Bangladesh & Two heuristic algorithms: multi-objective social engineering optimizer and multi-objective feasibility enhanced particle swarm optimization. \\
\hline
\citet{gilani2022data} & RO & VA & M,CV,SD & M & M,DC,VC & M, [Min] Economical cost, [Min] environmental cost, [Min] social cost. & Iran & Cutting Planes Approach \\
\hline
\citet{goodarzian2022sustainable} & MILP,CCP,S & D & M,CP & M & M,DC,VC & M, [Min] Operational cost [Min] environmental cost, [Max] social responsibility. & USA & Multi-Objective Teaching–learning-based Optimization Heuristic, Particle Swarm Optimization Heuristic, Genetic Algorithm Heuristic \\
\hline
\citet{goodarzian2022designing} & MILP & $\times$ & M,CV,SD & M & DC,VC & M, [Min] Operational, shortage, surplus cost, [Min] unmet demand, [Min] delivery time, [Min] environmental impact. & Iran & Gray Wolf Optimization Heuristic, Variable Neighborhood Search Heuristic \\
\hline
\citet{habibi2023designing} & SP,MC & D,LT & S,CV,SD & S & DC,VC & S, [Min] Operational cost. & Iran & Two heuristics: population-based simulated annealing and improved grey wolf optimizer. \\
\hline
\citet{kohneh2023optimization} & FP & D & M,CV,SD,P & M & DC,VC,PS & M, [Min] System cost, [Max] vaccination service level. & Mashhad, Iran & Model Linearization, Defuzzification \\
\hline
\citet{lai2021multi} & SP & D & M,CV,SD & M & VC,PS & S, [Min] Operational cost. & Synthetic & Benders decomposition-based heuristic algorithm. \\
\hline
\citet{li2021locate} & MINLP & $\times$ & S,CV,SD & S & VC,PS & M, [Min] Travel distance, cost, [Max] number of fully opened vaccination stations. & Shenzhen, China & Two-stage MILP Reformulation, $\epsilon$-Constraint Approach \\
\hline
\citet{lim2022redesign} & MILP & $\times$ & S,CV,SD & S & DC,VC & S, [Min] Operational, logistical cost. & Sub-Saharan Africa & Commercial Solver, Hybrid Heuristics for Large Scale \\
\hline
\citet{manupati2021multi} & MILP,DT & D,LT & S,CV,SD & M & M,DC,VC & S, [Min] Overall cost. & India & Commercial Solver \\
\hline
\citet{rahman2023optimising} & MILP,CCP & VW & M,CV,SD & M & DC,VC,PS & S, [Min] Operational cost. & Canberra, Australia & Heuristic Genetic Algorithm \\
\hline
\citet{rastegar2021inventory} & MILP & $\times$ & S,IV,SD & M & DC & S, [Max] Minimum vaccine distributed-to-demand ratio. & Iran & Commercial Solver, Linearization \\
\hline
\citet{sazvar2021capacity} & RO,FP & D & S,IV,SD & M & M,DC,PS & M, [Min] Economic cost, [Min] environmental cost, [Max] social responsibility, [Max] resilience. & Iran & Multi-choice goal programming with a utility function approach. \\
\hline
\citet{shiri2022equitable} & SP & D & M,CV,SD & M & DC,VC & S, [Min] Transportation, operation, shortage cost. & Iran & Scenario Reduction Algorithm \\
\hline
\citet{tang2022bi} & MILP & $\times$ & S,CV,SD & S & VC,PS & M, [Min] Operational, inventory cost, [Min] total travel distance. & Beijing City, China & Tailored Genetic Algorithm, DP Based Heuristic \\
\hline
\citet{tavana2021mathematical} & MILP & $\times$ & M,CV,SD & M & M,DC & S, [Max] Minimum delivery-to-demand ratio. & India & Commercial Solver \\
\hline
\citet{wang2023robust} & RO & S & S,CV,MD & M & DC,VC & M, [Max] vaccine benefits, [Min] cost. & Linshui, China & Column-and-Constraint Generation, Vertex Traversal and Dual Methods \\
\hline
\citet{xu2021hub} & MILP & $\times$ & S,CV,SD & S & DC,OC & [\textit{3}], [Max] Total coverage using three different coverage measures. & New Jersey, USA & Commercial Solver \\
\hline
\citet{yang2021optimizing} & MILP & $\times$ & S,ChV,SD & S & DC,VC & S, [Min] Annual hub facility, transportation costs. & Sub-Saharan Africa & Disaggregation and Merging Algorithm \\
\hline
\end{tabular}
\end{adjustbox}
\label{table:SCM-articles}
\end{center}
\end{table}

\FloatBarrier

\section{Appendix: Vehicle Routing Considerations for Vaccine Distribution}\label{sec: vehicle routing}

For the last-mile delivery of pandemic vaccines, routing of refrigerated trucks \citep{shamsi2021novel} or mobile vaccine clinics \citep{yang2021outreach,yucesoymobile} plays an important role. Routing problems are considered for the VC-PS part of the supply chain, under scenarios where vaccination center (VC) is mobile. All of the studies consider a single commodity, a single central depot and assume that the vehicles are equipped with all of the required equipment (i.e., vaccines, refrigerated storage units, and nurses to administer the vaccine). \citet{shamsi2021novel} propose a routing problem under an epidemiological model, after a disaster. Given the quantity of available vaccines and information about heterogeneous vehicles, the model distributes vaccines according to prioritization schemes by determining the time at which a vaccine is dropped at a population site (PS) (see \eqref{eq: vehicle-routing}):

\noindent \begin{align}
&\underbrace{\sum_{h} z_{r h} = 1, \hspace{0.25cm} \forall r}_{\substack{\text{A region can be} \\ \text{served by one vehicle.}}}; \quad \underbrace{\sum_{r} y_{r 0} \leq |\mathcal{H}|; \quad \sum_{r} y_{0 r} \leq |\mathcal{H}|}_{\substack{\text{All vehicles that depart the} \\ \text{depot should return back.}}}; \quad \underbrace{\sum_{r'} y_{r r'} = 1, \hspace{0.25cm} \forall r; \quad \sum_{r'} y_{r' r} = 1, \hspace{0.25cm} \forall r}_\text{Sequencing constraints.} \nonumber \\
& \underbrace{\sum_{r, r'} y_{r r'} \leq | \mathcal{R}^- | - 1, \hspace{0.25cm} \forall \mathcal{R}^- \subseteq \mathcal{R} \symbol{92} \{0\},  |\mathcal{R}^-| \geq 2}_\text{Sub-tour elimination constraint.} \label{eq: vehicle-routing} \\
& \underbrace{\sum_{r} d_{r g} z_{r h} \leq C^H_h, \hspace{0.25cm} \forall h}_\text{Vehicle capacity loading constraint.}; \quad \underbrace{T_{r'} \geq y_{r r'} \big( T_r + \tau_{rr'} + \tau^S_r\big), \hspace{0.25cm} \forall r,r', r \neq r'}_{\substack{\text{Consecutive vaccination start time} \\ \text{considering order.}}} \nonumber
\end{align}
This article uses time to serve a region and waiting time in each region to approximate the social cost per region in their objective function. On the other hand,  \citet{yang2021outreach} and \citet{yucesoymobile} assume homogeneous vehicles for routing vaccines to different priority groups (i.e., pregnant, elderly and disabled people) in different locations. \citet{yucesoymobile} assume that the demands at each node are known whereas \citet{yang2021outreach} use optimization under uncertainty techniques for handling demand uncertainty.

One major feature in both \citet{yang2021outreach} and \citet{yucesoymobile} is that the model selects population sites (PSs) to visit under demand information at hand, with the binary decision variable $Z^R_r$. Subsequently, the routing constraints are expressed as $\sum_{r'} y_{r r'} = Z^R_r, \hspace{0.25cm} \forall r$ and $\sum_{r'} y_{r' r} = Z^R_{r}, \hspace{0.25cm} \forall r$. Service time in each population site (PS) is treated as a decision variable and the total time budget constraint is formulated. Moreover, some articles consider routing within the entire supply chain design \citep{chowdhury2022modeling,habibi2023designing}. Further information on \textit{routing} articles can be found in Table~\ref{table:routing-articles}. 

\FloatBarrier

\begin{table}[!hbt] \setlength{\tabcolsep}{8pt}
\begin{center}
\caption{Detailed features of the optimization articles studying routing problems.}
\begin{adjustbox}{width=\textwidth}
\renewcommand{\arraystretch}{1.25}
\begin{tabular}{p{1.5in}|p{0.6in}|p{0.4in}|p{0.6in}|p{0.2in}|p{0.6in}|p{0.6in}|p{2.5in}|p{1in}|p{2.5in}}
\hline
\rowcolor{Gray}
\rotatebox{270}{\textbf{Article}} & \rotatebox{270}{\textbf{Model}} & \rotatebox{270}{\textbf{Uncertainty}} & \rotatebox{270}{\textbf{Commodity}} & \rotatebox{270}{\textbf{Period}} & \rotatebox{270}{\textbf{Vehicle}} & \rotatebox{270}{\textbf{Echelons}} & \rotatebox{270}{\textbf{Objective}} & \rotatebox{270}{\textbf{Case Study}} & \rotatebox{270}{\textbf{Solution Method  
 }} \\
\hline
\citet{shamsi2021novel} & VRP & D & S,GV,SD & S & RT,He & DC,PS & M, [Min] Social cost, [Min] routing cost. & Tehran, Iran & Dynamic Programming \\
\hline
\citet{yang2021outreach} & VRP,SCP & D,LT & S,ChV,SD & M & MVC,Ho & DC,PS & S, [Min] Vehicle and outreach trip cost. & Sub-Saharan Africa & MILP Commercial Solver, Multi-period Stochastic Modeling Approach \\
\hline
\citet{yucesoymobile} & LRP & $\times$ & S,CV,SD & S & MVC,Ho & DC,PS & S, [Max] Vaccination service levels. & Van, Turkey & Commercial Solver \\
\hline
\end{tabular}
\end{adjustbox}
\label{table:routing-articles}
\end{center}
\end{table}

\FloatBarrier

\section{Appendix: Additional Literature} \label{sec:sim-and-others}

We now review papers that do not use mathematical optimization  models but use  simulation approaches (\S \ref{subsec:sim}) to support decisions (\S \ref{subsec:other}) in vaccine management using descriptive statistics, observation, or other relevant methods in order to identify and conceptualize crucial steps. These papers highlight some additional, less studied aspects of vaccine management. Typically, these studies involve testing different scenarios using a simulation model and offering managerial insights for public health authorities, by observing the trends in the identified Key Performance Indicators (KPIs).

\subsection{Simulation-based Articles} \label{subsec:sim}

Table~\ref{table:simulation-based-classification} contains information about the decision levels, type of the simulation model, commodity, duration of the simulation, priority groups, and case studies used in \textit{simulation-based articles}. 

\paragraph{Prioritization:} A large number of the simulation papers analyzed questions on the prioritization of different population subgroups, due to vaccine scarcity at the early stages. The most common KPIs studied under a number of alternative prioritization schemes (e.g., by age groups, by profession, or by infection risk/immunity) are the number of infections, hospitalizations, and mortality. \citet{chen2020allocation} and \citet{rosenstrom2022can} also consider fairness by using the Gini index and race-based sub-population grouping, respectively. Papers that employ system dynamics models (i.e., SIR-type) addressed prioritization issues through the use of additional compartments that reflect differences between different population subgroups, but typically agent-based simulation models are deployed when analyzing the impact of interactions that may result in infections~\citep{chen2021prioritizing,oruc2021impact}.

The simulation based studies resulted in valuable insights shaping the COVID-19 vaccine deployment, such as prioritization of the elderly to minimize mortality, prioritization of younger and working individuals to minimize spread, benefits of dynamic prioritization (first vaccinate elderly, then switch to young), importance of equity objectives in reducing health disparities, etc. In addition, simulation approaches provided timelines for the duration of the pandemic as well as required number of vaccine doses to reach herd immunity \citep{chaturvedi2021predictive}.

\paragraph{Managing multiple vaccines:} COVID-19 vaccines have different, possibly time-dependent, efficacy, resource, storage and dosing requirements, and became available at different times during the pandemic. While most papers have considered a single type of vaccine, a few focused on issues related to managing multiple vaccines. For example, \citet{kim2021resource} provide the insight that a vaccine with fewer resource requirements and becoming available later can reduce the virus spread as much as a vaccine that becomes available earlier but requires more resources. In subsequent work, \citet{kim2022balancing} showed that if a vaccine with lower initial efficacy can be distributed fast, then it would be as effective as a vaccine with higher initial efficacy with a slower distribution rate, assuming that vaccine efficacy drops by time. In a similar vein, \citet{shim2021optimal} show that to minimize deaths under scarce supply of vaccines, COVID-19 vaccines with lower efficacy should be given to younger people and COVID-19 vaccines with higher efficacy should be given to older people, highlighting an important principle in prioritization with multiple vaccines. Finally in \citet{romero2021public} a simulation model was constructed to compare standard COVID-19 vaccination with a delayed second dose strategy, prioritizing the first dose, for vaccines with varying degrees of effectiveness. The results indicate that implementing a delayed second dose vaccination plan, specifically for individuals under the age of 65, could potentially lower the overall number of cumulative deaths under specific circumstances.

\paragraph{Management of vaccine inventory:} Vaccine vial size, recommended dosing, and population targeted for deployment are parameters that affect vaccine inventory management. For example, there are significant operational differences between storing and administering vaccines with 10-dose vials versus those with 2-dose vials. \citet{assi2011impact,assi2012influenza} build discrete-event simulation models where every vaccine with its storage location, refrigerator, freezer, and transport device is modeled as an entity to study the impact of vial size, vaccination period duration, and target vaccination population on supply chain bottlenecks.  

\subsection{Additional Considerations} \label{subsec:other}

\paragraph{Conceptual articles:} \citet{dai2021transforming} outline the COVID-19 supply chain research agenda under supply, demand, and supply-demand matching and highlight important areas such as cold chain management, demand eligibility, and vaccine hesitancy. \citet{forman2021covid} outline the obstacles to successful and coordinated international vaccine initiatives to prevent COVID-19, as well as proposed strategies to overcome them. \citet{schmidt2021equitable} use descriptive statistics to assess the equity of COVID-19 vaccine allocation within the US, and conclude that policymakers at the federal, state, and municipal levels should widely adopt the use of disadvantage indices and associated place-based indicators to ensure that fairness plays a prominent role in allocation strategies. \citet{alam2021challenges} use the decision-making trial and evaluation laboratory framework to examine the main challenges of the COVID-19 vaccine supply chain and come up with practical policy guidelines for the decision-makers. The article identifies the main bottlenecks of the {COVID-19} vaccine supply chain as (i) limited vaccine producers, (ii) improper collaboration with local organizations, (iii) lack of vaccine monitoring (iv) difficulties in regulating vaccine temperature, and (v) high cost of vaccinations. \citet{li2021strategies} describe the main methods used by a sizable hospital pharmacy department to overcome the challenges of preparing a large number of COVID-19 vaccine doses in a short period of time. \textcolor{black}{\citet{cano2023exploring} utilize social media data and natural language processing to demonstrate that vaccine supply chain challenges, especially in developing countries, significantly affect COVID-19 healthcare responses, emphasizing the critical need for improved coordination in vaccine supply chains to ensure equitable access.}

\paragraph{Vaccine hesitancy:} As mentioned in \citet{dai2021transforming}, a major obstacle to vaccination adoption is vaccine hesitancy, which the World Health Organization named as of the top threats to world health. Vaccine hesitancy is quantified using surveys. \citet{bogart2021covid} performed a questionnaire among HIV-infected African Americans about COVID-19 vaccine hesitancy and found that more than half of the participants supported at least one COVID-19 treatment reluctance or vaccination stance. \citet{silva2021covid} deployed a comparative questionnaire at a college campus to quantify COVID-19 and influenza vaccine hesitancy among college students. The results show that if the vaccinations are shown to be both safe and effective, college students are generally willing to get the vaccine.   \citet{sallam2021covid} reviewed the vaccine hesitancy surveys and came up with a vaccine hesitancy percentage for each country. \citet{troiano2021vaccine} performed a narrative review to find the underlying factors of vaccine hesitancy (i.e., age, gender, race, belief, education, work status) in the context of COVID-19. Finally, \citet{hegde2023two} proposed a framework with two distinct operating modes to examine service rates, which signify the capacity of the system to administer vaccinations under the first mode and the throughput affected by vaccine hesitancy under the second mode.

\paragraph{Waste management:} Mass COVID-19 vaccination practices increase medical waste dramatically. The  typical waste is vaccine vials, needles, syringes, and plastic vaccine administration equipment. \citet{hasija2022environmental} conceptualize and classify the types of vaccine waste and provide guidelines for vaccine waste management and treatment. Similarly, \citet{rayhan2022assessment} assess the COVID-19 vaccine-related waste management practices in Bangladesh. 

\FloatBarrier

\begin{table}[!hbt] \setlength{\tabcolsep}{8pt}
\begin{center}
\caption{Detailed features of the simulation-based articles.}
\begin{adjustbox}{width=\textwidth}
\renewcommand{\arraystretch}{1.25}
\begin{tabular}{l|l|l|l|l|l|l|l|l|l}
\hline
\rowcolor{Gray}
& \multicolumn{3}{c}{\textbf{Decision Level}} & \multicolumn{5}{c}{\textbf{Model Characteristics}} & \\
\rowcolor{Gray}
\textbf{Article} & \textbf{Inventory} & \textbf{Allocation} & \textbf{Distribution}  & \textbf{Type} & \textbf{SIR} & \textbf{Commodity} & \textbf{Period} & \textbf{Priority Groups} & \textbf{Case Study} \\

\hline
\rowcolor{Gray}
\textbf{Prioritization} \\

\hline

\citet{bubar2021model} & $\times$ & $\checkmark$ & $\times$ & \makecell[l]{SEIR \\ Simulation} & $\checkmark$ & \makecell[l]{Single COVID-19 \\ vaccine.} & One year. & \makecell[l]{$\checkmark$ Age based.} & \makecell[l]{USA} \\
\hline
\citet{chaturvedi2021predictive} & $\times$ & $\times$ & $\checkmark$ & \makecell[l]{SEIR \\ Simulation} & $\checkmark$ & \makecell[l]{Single COVID-19 \\ vaccine.} & Nine months. & \makecell[l]{$\times$} & \makecell[l]{India, Brazil, \\ USA} \\
\hline
\citet{chen2020allocation} & $\times$ & $\checkmark$ & $\times$ & \makecell[l]{SAPHIRE \\ Simulation} & $\checkmark$ & \makecell[l]{Single COVID-19 \\ vaccine.} & Three months. & \makecell[l]{$\checkmark$ Age based.} & \makecell[l]{New York, \\ USA} \\
\hline
\citet{chen2021prioritizing} & $\times$ & $\checkmark$ & $\times$ & \makecell[l]{Agent Based \\ Simulation} & $\checkmark$ & \makecell[l]{Single COVID-19 \\ vaccine.} & Four months. & \makecell[l]{$\checkmark$ Social \\ contact based.} & \makecell[l]{Virginia, \\ USA} \\
\hline
\citet{foy2021comparing} & $\times$ & $\checkmark$ & $\times$ & \makecell[l]{SEIR \\ Simulation} & $\checkmark$ & \makecell[l]{Single COVID-19 \\ vaccine.} & Five years. & \makecell[l]{$\checkmark$ Age based.} & \makecell[l]{India} \\
\hline
\citet{fujimoto2021significance} & $\times$ & $\checkmark$ & $\times$ & \makecell[l]{SEIR \\ Simulation} & $\checkmark$ & \makecell[l]{Single COVID-19 \\ vaccine.} & 550 days. & \makecell[l]{$\checkmark$ Antibody \\ testing based.} & \makecell[l]{Synthetic} \\
\hline
\citet{macintyre2022modelling} & $\times$ & $\checkmark$ & $\checkmark$ & \makecell[l]{SEIR \\ Simulation} & $\checkmark$ & \makecell[l]{Single COVID-19 \\ vaccine.} & 700 days. & \makecell[l]{$\checkmark$ Age, health \\ worker based.} & \makecell[l]{New South Wales, \\ Australia} \\
\hline
\citet{oruc2021impact} & $\times$ & $\checkmark$ & $\times$ & \makecell[l]{Agent Based \\ Simulation} & $\checkmark$ & \makecell[l]{Single COVID-19 \\ vaccine.} & 16 months. & \makecell[l]{$\checkmark$ Age, risk \\ based.} & \makecell[l]{Georgia, \\ USA} \\
\hline
\citet{rosenstrom2022can} & $\times$ & $\checkmark$ & $\times$ & \makecell[l]{Agent Based \\  Simulation} & $\checkmark$ & \makecell[l]{Single COVID-19 \\ vaccine.} & 17 months. & \makecell[l]{$\checkmark$ Income, \\ race based.} & \makecell[l]{North Carolina, \\ USA} \\
\hline
\citet{walker2022modeling} & $\times$ & $\checkmark$ & $\times$ & \makecell[l]{SEPIR \\ Simulation} & $\checkmark$ & \makecell[l]{Single COVID-19 \\ vaccine.} & 731 days. & \makecell[l]{$\checkmark$ Age, health \\ worker based.} & \makecell[l]{USA} \\

\hline
\rowcolor{Gray}
\textbf{Multiple Vaccine Management} \\

\hline

\citet{kim2021resource} & $\times$ & $\checkmark$ & $\times$ & \makecell[l]{SIRD \\ Simulation} & $\checkmark$ & \makecell[l]{Two different \\ COVID-19 vaccines.} & One year. & \makecell[l]{$\times$} & \makecell[l]{Synthetic} \\
\hline
\citet{kim2022balancing} & $\times$ & $\times$ & $\checkmark$ & \makecell[l]{SIRD \\ Simulation} & $\checkmark$ & \makecell[l]{COVID-19 vaccines \\ with different efficacy.} & One year. & \makecell[l]{$\times$} & \makecell[l]{Synthetic} \\
\hline
\citet{romero2021public} & $\times$ & $\checkmark$ & $\times$ & \makecell[l]{Agent Based \\  Simulation} & $\checkmark$ & \makecell[l]{COVID-19 vaccines \\ with different efficacy.} & 180 days. & \makecell[l]{$\checkmark$ Age, health \\ worker based.} & USA \\
\hline
\citet{shim2021optimal} & $\times$ & $\checkmark$ & $\times$ & \makecell[l]{SEIR \\ Simulation} & $\checkmark$ & \makecell[l]{COVID-19 vaccines \\ with different efficacy.} & One year. & \makecell[l]{$\checkmark$ Age based.} & \makecell[l]{South Korea} \\

\hline
\rowcolor{Gray}
\textbf{Management of Vaccine Inventory} \\

\hline

\citet{assi2011impact} & $\checkmark$ & $\checkmark$ & $\checkmark$ & \makecell[l]{Discrete Event \\ Simulation} & $\times$ & \makecell[l]{Multiple vial sized \\ measles vaccines.} & One year. & \makecell[l]{$\times$} & \makecell[l]{Niger} \\
\hline
\citet{assi2012influenza} & $\checkmark$ & $\checkmark$ & $\checkmark$ & \makecell[l]{Discrete Event \\ Simulation} & $\times$ & \makecell[l]{Ten different \\ influenza vaccines.} & One-six months. & \makecell[l]{$\checkmark$ Age based.} & \makecell[l]{Trang Province, \\ Thailand} \\
\hline
\citet{govindan2020decision} & $\times$ & $\times$ & $\times$ & \makecell[l]{Fuzzy Inference \\ System} & $\times$ & \makecell[l]{General demand \\ management.} & N/A. & \makecell[l]{$\checkmark$ Age, disease, \\ risk based.} & \makecell[l]{Synthetic} \\

\hline
\end{tabular}
\end{adjustbox}
\label{table:simulation-based-classification}
\end{center}
\end{table}

\end{APPENDIX}

\end{document}